\documentclass[12pt]{article}
\usepackage[a4paper,left=2cm,right=1.2cm,top=3cm,bottom=4cm]
{geometry}
\usepackage{amsmath,amstext, amsthm, empheq, extpfeil,  
amsbsy,amssymb,marvosym,fancyhdr,graphicx,amscd,amsfonts,latexsym,delarray,stackrel,lineno,color,cite,appendix,xcolor,url,hyperref}
\usepackage{wasysym}
\usepackage{subfig}

\usepackage{srcltx}
\usepackage{mathrsfs}
\usepackage{float} 
\usepackage[all]{xy}
\usepackage{t1enc}
\usepackage{mathrsfs}
\usepackage{pifont}

\usepackage{mathpple}
\usepackage[T1]{fontenc}

\definecolor{Green}{rgb}{0,1,0}
\definecolor{Blue}{RGB}{0,0,191}
\definecolor{mathmodecolor}{RGB}{0,102,0}
\definecolor{keywordcolor}{RGB}{0,51,151}
\definecolor{sourcebackgroundcolor}{RGB}{255,247,223}
\definecolor{unixagred}{RGB}{255,0,0}
\definecolor{lightgray}{RGB}{191,191,191}
\definecolor{green}{RGB}{1,191,191}

\newcommand*\patchAmsMathEnvironmentForLineno[1]{%
  \expandafter\let\csname old#1\expandafter\endcsname\csname #1\endcsname
  \expandafter\let\csname oldend#1\expandafter\endcsname\csname end#1\endcsname
  \renewenvironment{#1}%
     {\linenomath\csname old#1\endcsname}%
     {\csname oldend#1\endcsname\endlinenomath}}%
\newcommand*\patchBothAmsMathEnvironmentsForLineno[1]{%
  \patchAmsMathEnvironmentForLineno{#1}%
  \patchAmsMathEnvironmentForLineno{#1*}}%
\AtBeginDocument{%
\patchBothAmsMathEnvironmentsForLineno{equation}%
\patchBothAmsMathEnvironmentsForLineno{align}%
\patchBothAmsMathEnvironmentsForLineno{flalign}%
\patchBothAmsMathEnvironmentsForLineno{alignat}%
\patchBothAmsMathEnvironmentsForLineno{gather}%
\patchBothAmsMathEnvironmentsForLineno{multline}%
}

\newtheorem{thm}{Theorem}[section]
\newtheorem{prop}[thm]{Proposition}

\newtheorem{lem}[thm]{Lemma}

\newtheorem{rem}[thm]{Remark}

\newtheorem{hypo}[thm]{Hypothesis}
\newtheorem{fact}[thm]{Fact}

\newcommand{\dstar}{d^*u}

\def\GL{{\rm GL}}

\def\Qer{{\rm Ker}}

\def\A{{\mathbb A}}

\def\C{{\mathbb C}}
\def\F{{\mathbb F}}

\def\Q{{\mathbb Q}}
\def\R{{\mathbb R}}
\def\Z{{\mathbb Z}}

\def\cE{{\mathcal E}}

\def\cH{{\mathcal H}}

\def\cM{{\mathcal M}}

\def\cS{{\mathcal S}}

\def\rs{{\R_+^*}}

\def\sin{{{\rm sin}}}
\def\cos{{{\rm cos}}}

\def\dim{{\mbox{dim}}}
 \def\scal2{{\mathscr S}}
 
\def\fourier{\F}
\def\sr0{{\cS^{\rm ev}_0}}

\def\sar0{{\cS_0(\A_\Q)}}

\def\fourier{\F}

\newcommand{\Rad}{{\rm Rad}}

\newcommand{\Ps}{{\operatorname Ps}}
\newcommand{\PS}{{\operatorname PS}}

\newcommand{\fixme}[1]{\centerline{\color{red}\itshape{#1}}}
\newcommand{\nil}[1]{}

\title{The uncertainty operator, the Weil ground state and the zeros of zeta:\\ a numerical study and conditional positivity of the Weil quadratic form}

\author{Ronnie Andrews, Alain Connes, Antonio La Gatta}
\date{2026}

\begin{document}
\renewcommand{\fixme}[1]{\par\noindent{\color{red}\itshape #1}\par}
\maketitle

\begin{abstract}
Applying the Dirichlet variational principle to the Weil quadratic form $QW_\lambda$ restricted to $L^2([\lambda^{-1},\lambda],\dstar)$—which involves only the primes $p\le\lambda^2$ of the Euler product—yields a canonical finite-dimensional Hermitian matrix approximation $QW_\lambda^N$ whose ground state has a Mellin transform with zeros proved to be real in full generality whenever the lowest eigenvalue is simple and even under $u\mapsto u^{-1}$. While $\lambda^2=13$ in previous work \cite{CCM,CoJOMP} already recovered the first $50$ non-trivial zeros of $\zeta$, we show here that increasing $\lambda$ to $\lambda\sim16$ ($\lambda^2=250$, involving only the $53$ primes $p<250$) reproduces the first thousand zeros of $\zeta$, each with more than $500$ correct decimal places, the first one with $1341$, and all $2469$ zeros up to height $3000$, with a certified count showing that the finite transform has no other zero, real or complex, in this range.

To explain this phenomenon and test the strategy outlined in \cite{CoJOMP}, we establish a precise dictionary between the low-lying spectrum of $QW_\lambda$ and the spectrum of the time-frequency uncertainty operator $U_\lambda=2-P_\lambda-\widehat P_\lambda$ on even functions, where $P_\lambda$ is the projection on functions supported in $[-\lambda,\lambda]$ and $\widehat P_\lambda$ its Fourier conjugate. The latter is diagonalized by the prolate spheroidal wave functions, with minuscule eigenvalues $\delta_n(\lambda)=1-\vert\chi_n(\lambda)\vert$ governed by Fuchs' asymptotics $\delta_n(\lambda)\sim c_n\lambda^{2n+1}e^{-4\pi\lambda^2}$, and its two sectors (Fourier transform close to $\pm1$) correspond via the arithmetic map $\cE(f)(u)=u^{1/2}\sum_{n\ge1}f(nu)$ to the two sectors of $QW_\lambda$ under $u\mapsto u^{-1}$. While the upper bound $\epsilon_\lambda(j)\le C_j\delta_{4j}(\lambda)(\log\lambda)^3$ is a theorem obtained with Kastoris \cite{K,CK} after a Gevrey smoothing, we test the lower bound far beyond the range of previous computations (up to $\lambda^2=300$, where $\epsilon_\lambda(1)\sim 10^{-1618}$) and find a striking coincidence between the logarithms of the minuscule eigenvalues of $QW_\lambda$ and the Fuchs ladder. In particular, the lowest eigenvalue $\epsilon_\lambda(1)$ is confirmed to be simple, isolated, and even, separated from the next eigenvalues by factors $\asymp\lambda^4$ (odd) and $\asymp\lambda^8$ (even). Furthermore, the Gram--Schmidt vectors built from the $\cE$-images of zero-integral prolate functions reproduce all minuscule eigenvectors of $QW_\lambda$, and conditioning the form by them and by the pole functional removes the minuscule spectrum altogether, leaving a form bounded below by $\mu_\lambda\approx 2/\lambda$. The picture that emerges is that of a direct sum: $QW_\lambda$ decomposes into a finite-dimensional part, very close to the radical of the global Weil form and governed by the Fuchs ladder of $U_\lambda$, and an infinite-dimensional strictly positive part originating in the Sonin space $\operatorname{Ker}P_\lambda\cap\operatorname{Ker}\widehat P_\lambda$, where $U_\lambda$ takes its maximal value $2$. Under the Fuchs ladder hypothesis and the positivity of $QW_\lambda$, we obtain quantitative estimates for the convergence of the Weil ground state to an explicit target whose Mellin transform is $\tfrac14\xi(s+\tfrac12)$, setting up the approach to the Riemann Hypothesis via Hurwitz' theorem. Finally, in Section 8 we use the semilocal trace formula in its quantized-differential form—writing $QW_\lambda$ as the trace of the scaling action against the uncertainty operator at scale $1$—to obtain a first analytical justification of the direct sum picture: we establish the semilocal uncertainty principle from the non-vanishing of the semilocal Sonin space by Coburn's lemma and prove, under one quantitative hypothesis on the semilocal cutoff projections holding in the archimedean case, that $QW_\lambda$ is bounded below by an explicit positive constant on the orthogonal complement of an explicit finite-dimensional space built from the near-self-dual modes. Combined with the monotonicity of $\epsilon_\lambda$ in $\lambda$ and uniform Kastoris estimates, this reduces the positivity of $QW_\lambda$ to a two-sided estimate on the Schur complement of the form on the finitely many modes at the edge of the Fuchs ladder and on the two pole directions.


\paragraph{Mathematics Subject Classification 2020.}
Primary
\href{http://www.ams.org/mathscinet/msc/msc2020.html?t=11M06&btn=Current}{11M06},
\href{http://www.ams.org/mathscinet/msc/msc2020.html?t=11M26&btn=Current}{11M26},
\href{http://www.ams.org/mathscinet/msc/msc2020.html?t=11M55&btn=Current}{11M55};
Secondary
\href{http://www.ams.org/mathscinet/msc/msc2020.html?t=58B34&btn=Current}{58B34},
\href{http://www.ams.org/mathscinet/msc/msc2020.html?t=33E10&btn=Current}{33E10},
\href{http://www.ams.org/mathscinet/msc/msc2020.html?t=47B35&btn=Current}{47B35},
\href{http://www.ams.org/mathscinet/msc/msc2020.html?t=47A75&btn=Current}{47A75},
\href{http://www.ams.org/mathscinet/msc/msc2020.html?t=65F15&btn=Current}{65F15}.
\end{abstract}

\tableofcontents

\section{Introduction}

Let $QW$ be the Weil quadratic form and, for $\lambda > 1$, let $QW_\lambda$ be its restriction to functions supported in $[\lambda^{-1},\lambda]$, viewed as a quadratic form on $L^2([\lambda^{-1},\lambda],\dstar)$, $\dstar = du/u$. By Weil's criterion, the positivity of $QW_\lambda$ for all $\lambda$ is equivalent to the Riemann Hypothesis. A remarkable feature of $QW_\lambda$ is that its evaluation on the window $[\lambda^{-1},\lambda]$ involves only the finite set of prime powers $p^m \le \lambda^2$, and therefore uses only finitely many factors of the Euler product of $\zeta$. In \cite{CCM} and in the ``Letter to Riemann'' \cite{CoJOMP}, a variational route to the zeros of $\zeta$ was proposed in the spirit of Riemann's Dirichlet principle: one minimizes $QW_\lambda$ on the unit sphere of $L^2([\lambda^{-1},\lambda],\dstar)$ (or of its canonical finite-dimensional Galerkin truncation $QW_\lambda^N$ to the first $2N+1$ log-Fourier modes) and takes the Mellin transform $\mathcal{M}\eta(s)=\int_{\lambda^{-1}}^\lambda \eta(u)u^s\dstar$ of the resulting ground state $\eta$. By a general theorem proved in \cite{CvS} (recalled as Theorem \ref{thm:realzeros} below), which relies on the Galerkin structure of the truncated matrix and requires no unproved hypothesis on $\zeta$, as soon as the lowest eigenvalue $\epsilon_\lambda^N(1)$ of $QW_\lambda^N$ is simple and its eigenvector $\eta$ is even under the inversion $\iota(u)=u^{-1}$, \emph{all} the zeros of the entire function $t\mapsto \mathcal{M}\eta(it)$ are real—that is, all the zeros of $\mathcal{M}\eta(s)$ lie exactly on the critical line $\Re s=0$ (corresponding to $\Re w = \tfrac12$ for $w=s+\tfrac12$).

In \cite{CCM,CoJOMP}, taking $\lambda^2=13$ (which involves only the six primes $2, 3, 5, 7, 11, 13$) and $N=120$ already reproduced the first $50$ non-trivial zeros of $\zeta$, the first one to within $2.4\cdot 10^{-55}$. The first major finding of the present paper is the dramatic scaling of this phenomenon when $\lambda$ increases and more primes are included in the window: for $\lambda\sim 16$ ($\lambda^2=250$, involving only the $53$ primes $p<250$, with $N=4000$), the zeros of the Mellin transform of the ground state of $QW_\lambda^N$—all proved real a priori—reproduce the \emph{first thousand zeros of $\zeta$, each with more than five hundred correct decimal places}, the first zero with $1341$ correct decimal places, and all the $2469$ zeros of height at most $3000$ with at least twelve correct digits. Moreover, the error on each zero is proportional to the lowest eigenvalue $\epsilon_\lambda^N(1)\approx 1.1\cdot 10^{-1346}$ of the truncated form, and a certified ball-arithmetic sign count shows that the finite transform has no other zero in this range, real or complex (Section \ref{sec:zeros}, Figure \ref{fig:mellindigits}, and Remark \ref{rem:countzeros}).

To understand why the Dirichlet principle on a finite-prime window captures the zeros of $\zeta$ with such precision, and to advance the convergence strategy outlined in \cite{CoJOMP} after the Letter to Riemann, we investigate the low-lying spectrum of $QW_\lambda$. It was observed in  \cite{CC,CCM}  that the lowest eigenvalue $\epsilon_\lambda$ is minuscule, of the order of $e^{-4\pi\lambda^2}$, and that the Weil ground state is remarkably close to the image under the arithmetic map
\begin{equation}
\cE(f)(u) := u^{1/2}\sum_{n\ge 1} f(nu)
\end{equation}
of a two-mode combination of prolate spheroidal wave functions on $[-\lambda,\lambda]$. Here we explain this phenomenon by establishing a precise dictionary between $QW_\lambda$ and an exactly solvable operator on $L^2(\R)_{\mathrm{ev}}$ arising in the ``time-frequency'' uncertainty analysis of Slepian, Pollak, and Landau \cite{Sle}, and we combine high-precision numerical experiments (up to $\lambda^2=300$, where $\epsilon_\lambda(1)\sim 10^{-1618}$ requires arithmetic with thousands of digits) with two analytical tools—the Gevrey approximant method developed with Kastoris \cite{K,CK} and the semilocal trace formula \cite{Co99, CM}—which provide the first elements of proof for the emerging spectral picture.

The exactly solvable operator is the time-frequency uncertainty operator
\begin{equation}\label{eq:unc}
U_\lambda := 2 - P_\lambda - \widehat P_\lambda, \qquad \widehat P_\lambda = \fourier P_\lambda \fourier,
\end{equation}
where $P_\lambda$ is the orthogonal projection of $L^2(\R)_{\mathrm{ev}}$ on the functions supported in $[-\lambda,\lambda]$ and $\fourier$ is the Fourier transform for the character $e^{2\pi ixy}$ (so that $\fourier^2=1$ on $L^2(\R)_{\mathrm{ev}}$). Thus $\widehat P_\lambda$ is the projection on the functions band-limited to $[-\lambda,\lambda]$, and $U_\lambda=(1-P_\lambda)+(1-\widehat P_\lambda)$ measures the failure of a function to be simultaneously time- and band-limited. By construction $U_\lambda$ is equal to $2$ on the infinite-dimensional Sonin space \cite{CC} of even functions in the kernel of $P_\lambda$ and of $\widehat P_\lambda$. By the theory of two projections, $U_\lambda$ is diagonalized by the prolate spheroidal wave functions $h_{n,\lambda}$: its eigenvectors are $h_{n,\lambda}\pm \hat h_{n,\lambda}$ with eigenvalues $1\mp\vert\chi_n(\lambda)\vert$, where $\chi_n(\lambda)$ is the eigenvalue of the finite Fourier transform on $h_{n,\lambda}$ (Theorem \ref{thm:unc}). The minuscule eigenvalues $\delta_n(\lambda)=1-\vert\chi_n(\lambda)\vert$, $n$ even, form a ladder described by Fuchs' asymptotics \cite{Fu},
\begin{equation}\label{eq:fuchsintro}
\delta_n(\lambda) = \frac{2\sqrt{\pi}\,8^n(2\pi)^{n+1/2}}{n!}\,\lambda^{2n+1}e^{-4\pi\lambda^2}\bigl(1+O(\lambda^{-2})\bigr),
\end{equation}
and they split into two interlacing sectors: $n\equiv 0\pmod 4$, where the Fourier transform is close to the identity on $h_{n,\lambda}$, and $n\equiv 2\pmod 4$, where it is close to minus the identity.

The dictionary between $QW_\lambda$ and $U_\lambda$ rests on two facts. First, by Poisson summation the map $\cE$ intertwines the Fourier transform with the inversion $\iota(u)=u^{-1}$: $\cE(\hat f)=\cE(f)\circ\iota$ for $f$ even with $f(0)=\hat f(0)=0$. Hence a function in $\operatorname{Range}P_\lambda$ has $\cE(f)$ supported in $u\le\lambda$, a function in $\operatorname{Range}\widehat P_\lambda$ has $\cE(f)$ supported in $u\ge\lambda^{-1}$, and the window $[\lambda^{-1},\lambda]$ is the image of $\operatorname{Range}P_\lambda\cap\operatorname{Range}\widehat P_\lambda$, a space which is $\{0\}$ exactly by the uncertainty principle, but is approximately spanned by the prolate functions with $\vert\chi_n\vert\approx 1$. The two sectors of $U_\lambda$ become the two sectors of $QW_\lambda$ under $\iota$, with which $QW_\lambda$ commutes, while the zero-integral constraint $\hat f(0)=\int f=0$ required for $\cE(f)$ to lie in $L^2(\dstar)$ removes the lowest rung $\delta_0$ (resp.\ $\delta_2$) in each sector. Second, $\cE$ maps the even Schwartz functions with $f(0)=\hat f(0)=0$ into the radical of the global form $QW$ \cite{CC}, so that the only obstruction to $\cE(f)$ being a null vector of $QW_\lambda$ is its mass outside the window, which is precisely the frequency leakage measured by $U_\lambda$. Made quantitative via a Gevrey smoothing at the endpoints $\pm\lambda$ which makes the construction legitimate, this provides the first rigorous element of proof for the dictionary—the upper bound theorem obtained with Kastoris (\cite{K}, sharpened in \cite{CK}):
\begin{equation}
\epsilon_\lambda \le C\,\delta_4(\lambda)(\log\lambda)^3,
\end{equation}
and the same mechanism yields uniform upper bounds $\epsilon_\lambda(j)\le C_j\delta_{4j}(\lambda)(\log\lambda)^3$ for the $j$-th eigenvalue in the even sector, and $\delta_{4j+2}(\lambda)$ in the odd one (Theorem \ref{thm:upperK}).

What is not yet a theorem is the lower bound direction: that nothing on the Weil side is smaller than the uncertainty operator allows. We formulate this as Hypothesis \ref{hyp:LA} (logarithmic agreement with the Fuchs ladder) and test it numerically in Section \ref{sec:numspec} over a vast range ($5\le\lambda^2\le 300$). The second main numerical discovery of the paper is the remarkable coincidence of the logarithms of the minuscule numbers on both sides:
\begin{enumerate}
\item Across more than $1600$ orders of magnitude, the lowest eigenvalue $\epsilon_\lambda(1)$ tracks the exact prolate rung $\delta_4(\lambda)$, with the ratio $\epsilon_\lambda(1)/(\delta_4(\lambda)\log\lambda)$ remaining bounded between $4.9$ and $8.8$ without upward trend (Fact \ref{fact:logbound} and Table \ref{tab:logbound}).
\item The logarithmic gap $\log\epsilon_\lambda(2)-\log\epsilon_\lambda(1)$ agrees with the Fuchs gap $\log(\delta_8/\delta_4)=8\log\lambda+8.2425+o(1)$ to within a nearly constant offset of $0.5$--$0.8$ over the entire range $5\le\lambda^2\le 300$, and likewise $\log\epsilon_\lambda^{\mathrm{odd}}(1)-\log\epsilon_\lambda(1)$ agrees with $\log(\delta_6/\delta_4)=4\log\lambda+4.4342+o(1)$ to within $0.3$--$0.4$.
\item The even and odd eigenvalues strictly interlace, $\epsilon_\lambda(1)<\epsilon_\lambda^{\mathrm{odd}}(1)<\epsilon_\lambda(2)<\epsilon_\lambda^{\mathrm{odd}}(2)<\cdots$, for every minuscule eigenvalue, and the even and odd sectors contain $\lambda^2-1$ and $\lambda^2-2$ minuscule eigenvalues respectively, one fewer than the corresponding classes of $U_\lambda$ ($\lambda^2$ and $\lambda^2-1$ rungs, in agreement with Slepian's count $4\lambda^2$ split into four classes), the missing rung being the one removed by $\int f=0$ (Fact \ref{fact:count}).
\end{enumerate}
Crucially, these computations show that the lowest eigenvalue $\epsilon_\lambda(1)$ of $QW_\lambda$ is \emph{simple, isolated, and $\iota$-even}, separated from $\epsilon_\lambda^{\mathrm{odd}}(1)$ by a factor $\asymp\lambda^4$ and from $\epsilon_\lambda(2)$ by $\asymp\lambda^8$—thus verifying the key condition required by Theorem \ref{thm:realzeros} to guarantee that all zeros of the Mellin transform of the ground state lie on the critical line.

Section \ref{sec:conditioning} completes this spectral picture on the side of the eigenvectors and reveals a fundamental structural decomposition of $QW_\lambda$. The Gram--Schmidt orthogonalization of the $\cE$-images of the zero-integral prolate vectors reproduces \emph{every} minuscule eigenvector of $QW_\lambda$ (with $\sin\angle\le 2.6\cdot 10^{-2}$ up to the plunge, and down to $1.9\cdot 10^{-7}$ for the ground state at $\lambda^2=150$, Fact \ref{fact:gs}). Moreover, when we condition $QW_\lambda$ orthogonally to these $\lambda^2-1$ (even) and $\lambda^2-2$ (odd) prolate vectors together with the pole functional $\ell(f)=\int f(u)u^{1/2}\dstar$, \emph{all minuscule eigenvalues disappear}: the lowest remaining eigenvalue $\mu_\lambda$ is of order $2/\lambda$ (between $0.78$ at $\lambda^2=5$ and $0.19$ at $\lambda^2=150$, Fact \ref{fact:cond}). The picture which emerges is therefore that of a \textbf{direct sum}: on the window $[\lambda^{-1},\lambda]$, the Weil quadratic form is the sum of
\begin{itemize}
\item a \emph{finite-dimensional near-radical part} (of dimension $\approx 2\lambda^2$), very close to the radical of the global Weil form and governed by the low-lying Fuchs ladder of the uncertainty operator $U_\lambda$; and
\item an \emph{infinite-dimensional strictly positive part}, on which the form is bounded below by a positive constant and whose origin lies in the Sonin space $\operatorname{Ker}P_\lambda\cap\operatorname{Ker}\widehat P_\lambda$, where the uncertainty operator $U_\lambda$ takes its maximal value $2$.
\end{itemize}

Under Hypothesis \ref{hyp:LA} and the positivity of $QW_\lambda$, we prove in Section \ref{sec:ground} that the ground state $\eta_\lambda$ is simple and $\iota$-even, that it lies within $O(\lambda^{-4}(\log\lambda)^{3/2})$ in $L^2$ of the explicit Kastoris test vector, and we obtain quantitative bounds for its distance to the prolate state $\tau_\lambda=\cE(h_\lambda)/\cE(h_\lambda)(1)$ (at rate $O(\lambda^{-7/2}(\log\lambda)^{3/2})$, confirmed numerically as $\approx C\lambda^{-3.5}$ in Table \ref{tab:dist}) and to the explicit Hermite target $\tau(u)=\theta(u)/\theta(1)$,
\begin{equation}\label{eq:targetintro}
\theta(u) = \sqrt{u}\sum_{n\ge 1}\frac{\pi}{2}\,n^2u^2\bigl(2\pi n^2u^2 - 3\bigr)e^{-\pi n^2u^2},
\end{equation}
at rate $O(\lambda^{-2})$ in the norm $\|g\|_1=\int_1^\lambda\vert g(u)\vert u^{-1/2}du$ which controls the Mellin transform uniformly on the closed critical strip $\vert\Re s\vert\le\tfrac12$. Since the Mellin transform of $\theta$ is exactly $\tfrac14\xi(s+\tfrac12)$ (Theorem \ref{thm:mellintheta}) and the zeros of the Mellin transform of the ground state are all real (Theorem \ref{thm:realzeros}), uniform convergence of the normalized ground state to $\tau$ on compact subsets of the strip implies the Riemann Hypothesis by Hurwitz' theorem (Proposition \ref{prop:hurwitz}). We make the logical structure of this implication explicit in Section \ref{sec:zeros}—in particular the fact that the positivity of $QW_\lambda$ assumed in Section \ref{sec:ground} cannot be used in an unconditional proof—before presenting the high-precision zero computations for $\lambda^2=250$ and $\lambda^2=256$.

Finally, Section \ref{sec:condpos} provides a first analytical justification for the direct sum picture and turns the dictionary into a tool for proving conditional positivity. The semilocal trace formula of \cite{Co99}, in the quantized-differential form of \cite{CM}, Chapter 2, writes the Weil form of a test function $f$ supported in the window, with vanishing pole terms $\mathcal{M}f(\pm\tfrac12)=0$, as a trace,
\[
QW_\lambda(f,f) = \operatorname{Tr}\bigl(\theta(f)(1-P-\widehat P)\theta(f)^*\bigr),
\]
where $\theta(f)$ is the scaling convolution action and $1-P-\widehat P$ is the uncertainty operator at scale $1$ on the semilocal ad\`ele class space, shifted by $-1$ (Theorem \ref{thm:TF}). Three steps then explain the prolate conditioning of Section \ref{sec:conditioning} and lead to a rigorous lower bound:
\begin{enumerate}
\item On the semilocal space, the operator $1-P-\widehat P$ equals $+1$ on the Sonin space and acts on each prolate pair by $\pm\kappa_n$, so that its negative part is carried solely by the near-self-dual directions, with eigenvalues $-\kappa_n\uparrow 0$. The absence of any eigenvalue $-1$ (i.e.\ of exactly self-dual directions) is the \emph{semilocal uncertainty principle}, which we prove in Lemma \ref{lem:H1} from the infinite-dimensionality of the semilocal Sonin space \cite{CCqi,CCMsl} by applying Coburn's lemma on Toeplitz operators.
\item Because $f$ is supported in $[\lambda^{-1},\lambda]$, the scaling action $\theta(f)$ transports the cutoff projections at scale $\sigma$ into those at scale $\lambda\sigma$ (Lemma \ref{lem:transport}). Consequently, the negative part at scale $\sigma$ can only act through the inner products $\langle\theta(f)e_n, e_m^\lambda\rangle=\langle f, e_m^\lambda * e_n^*\rangle$ of $f$ with the convolved near-self-dual modes $e_m^\lambda$ at scale $\lambda\sigma$ that are more self-dual than $e_n$ (Proposition \ref{prop:pairexp}); orthogonality to these explicit vectors is precisely the conditioning.
\item Shifting the cutoffs from scale $1$ to scale $\sigma>1$ extracts from the projections an explicit positive margin $2\log\sigma\,\|f\|^2$ (Lemma \ref{lem:shift}).
\end{enumerate}
The outcome is Theorem \ref{thm:condpos}: under a  quantitative hypothesis ($\mathrm{H2}'$) on the high-frequency tail of the semilocal pair expansion—which holds unconditionally in the archimedean case and is reduced to the localization estimate \eqref{eq:H2ii} on smoothed spaces—one has $QW_\lambda(f,f)\ge\tfrac12\log\sigma\,\|f\|^2$ on the orthogonal complement of an explicit finite-dimensional subspace built from the near-self-dual modes (of codimension $O(\sigma^4\lambda^2)$ in the archimedean case). Finally, since the lowest eigenvalue $\epsilon_\lambda$ is non-increasing in $\lambda$, a minuscule lower bound $\epsilon_\lambda\ge -\eta(\lambda)$ with $\eta(\lambda)\to 0$ suffices to imply $\epsilon_\lambda\ge 0$ for all $\lambda$ (Proposition \ref{prop:monotone}). Combined with two-sided Kastoris estimates on the deeply minuscule rungs of the conditioning space, this reduces the Riemann Hypothesis to a two-sided estimate on the Schur complement of $QW_\lambda$ on the finitely many transition modes at the edge of the Fuchs ladder and on the two pole directions (Proposition \ref{prop:reduction}). We emphasize, however, that what really matters in Section \ref{sec:condpos} is the description of the method and the interplay of its main ingredients—the two-parameter shift in the trace formula, the role of the Sonin space, the scale transport by $\theta(f)$, and the Schur--monotonicity reduction—rather than the precise technical statements of Theorem \ref{thm:condpos} and Proposition \ref{prop:reduction}, which are only a first step and not yet sharp compared with the numerical conditioning of Section \ref{sec:conditioning}.

\paragraph{Earlier works.} The minuscule nature of the eigenvalues of $QW_\lambda$ as well as the prolate approximation of the eigenvectors were discovered in \cite{CC}. The approximation of the zeros of $\zeta$ by those of the Mellin transform of the ground state was discovered in \cite{CCM} (see also \cite{CoJOMP}), with $N=120$ and $\lambda^2\le 14$. A. Groskin \cite{Gr} extended these computations to $\lambda^2=100$, obtaining the first ten zeros with $307$ to $329$ correct decimals, and R.A. carried them to $\lambda=16$ and the first $200$ zeros, the starting point of the computations of \S\ref{ss:zeroscomp}. In connection with a conjecture of D. \'Sliwi\'nski \cite{Sl}, R.A. \cite{An} also studied the error over the whole finite spectrum along the scaling $\lambda=N$, finding that it grows while the low-lying zeros become more accurate, in agreement with the breakdown of the agreement described in \S\ref{ss:zeroscomp}. The construction of prolate approximants in the radical of the Weil quadratic form, and the first upper bound for the lowest eigenvalue in terms of the Fuchs estimate, are due to D. Kastoris \cite{K}; his bound carries an unspecified power of $\lambda$ in front of $\delta_4(\lambda)$, and the improvement of \cite{CK}, which replaces it by $(\log\lambda)^3$ as in \eqref{eq:upper}, is what makes the bound usable here, since the rungs of the ladder are separated by powers of $\lambda$. The constant $C_*$ of \eqref{eq:Cstar} is also due to him.

Throughout, $\lambda$ is the half-length of the window; the parameter $c=\lambda^2$ of \cite{CCM} is never used, and $\gamma := 2\pi\lambda^2$ is the prolate parameter.

\paragraph{Acknowledgements.} The computations of A.L.G. run on a cluster whose power is equivalent to that of one hundred A100 GPUs. We are indebted to Claude Shannon, both for creating information theory, which led to the work of Slepian and his collaborators on which this paper rests, and for the AI assistant bearing his name, which we used in the elaboration of this paper.

\section{The two operators}\label{sec:two}

\subsection{The Weil quadratic form on the window}

For real functions $g,h$ on $\rs$ put $F(x)=\int_0^\infty g(u)h(ux)\dstar$. The Weil quadratic form is $QW(g,h)=\Psi(F)$ where, in the normalization of \cite{CCM},
\begin{equation}\label{eq:Psi}
\Psi=W_{0,2}-W_\R-\sum_pW_p ,\qquad W_{0,2}(F)=\int_0^\infty F(x)(x^{1/2}+x^{-1/2})\,d^*x,
\end{equation}
\begin{equation}\label{eq:Wp}
W_p(F)=\log p\sum_{m\ge1}p^{-m/2}\big(F(p^m)+F(p^{-m})\big),
\end{equation}
\begin{equation}\label{eq:WR}
W_\R(F)=(\log4\pi+\gamma_E)F(1)+\int_1^\infty\Big(F(x)+F(x^{-1})-2x^{-1/2}F(1)\Big)\frac{x^{1/2}}{x-x^{-1}}\,d^*x .
\end{equation}
$QW_\lambda$ is the restriction of $QW$ to functions supported in $[\lambda^{-1},\lambda]$. Each term of \eqref{eq:Psi} is invariant under $F\mapsto F\circ\iota$, $\iota(x)=x^{-1}$, and $g\mapsto g\circ\iota$ is unitary on $L^2(\rs,\dstar)$ and preserves the window; hence $QW_\lambda$ commutes with $\iota$ and splits as $QW_\lambda^{\rm ev}\oplus QW_\lambda^{\rm odd}$ on the $\pm1$ eigenspaces of $\iota$. In the variable $t=\log u$ the archimedean term is, up to bounded operators, the Fourier multiplier $\log\vert\tau\vert$, so that $QW_\lambda$ is a bounded-below form with compact resolvent, whose eigenvalues we denote $\epsilon_\lambda(1)\le\epsilon_\lambda(2)\le\cdots$ in the even sector, $\epsilon^{\rm odd}_\lambda(1)\le\cdots$ in the odd one, and $\epsilon_\lambda=\epsilon_\lambda(1)$ (the numerics confirm that the lowest eigenvalue is even). Numerically $QW_\lambda$ is computed in the basis of the $2N+1$ lowest periodic log-Fourier modes on $[\lambda^{-1},\lambda]$, see Appendix \ref{app:matrix}.

\subsection{The arithmetic map and the radical}

Let $\sr0=\{f\in\cS(\R)^{ev}:f(0)=\hat f(0)=0\}$ with $\hat f(y)=\int f(x)e^{2\pi ixy}dx$. The basic facts are \cite{CC}:
\begin{equation}\label{eq:radical}
\cE(\sr0)\subset\Rad(QW),\qquad \cE(\hat f)=\cE(f)\circ\iota\quad(f\in\sr0),
\end{equation}
the second being Poisson summation $\sum_{n\in\Z}f(nu)=u^{-1}\sum_{k\in\Z}\hat f(k/u)$ with the terms $n=0$, $k=0$ removed. If $f\in\sr0$ is supported in $[-\lambda,\lambda]$ then $\cE(f)(u)=0$ for $u>\lambda$ and $\cE(f)(u)=u^{1/2}f(u)$ for $u>\lambda/2$; if $\hat f$ is supported in $[-\lambda,\lambda]$ then $\cE(f)(u)=0$ for $u<\lambda^{-1}$. In the Mellin variable, with $\cM g(s)=\int_0^\infty g(u)u^s\dstar$,
\begin{equation}\label{eq:mellinE}
\cM\cE(f)(s)=\zeta(s+\tfrac12)\int_0^\infty f(x)x^{s+1/2}\,d^*x ,
\end{equation}
so that the critical strip of $\zeta$ corresponds to $\vert\Re s\vert\le\tfrac12$ and the Hilbert space $L^2(\rs,\dstar)$ to the critical line.

\subsection{Prolate spheroidal wave functions}\label{ss:prolate}

Let $\gamma=2\pi\lambda^2$. The finite Fourier transform $F_\lambda=\mathbf 1_{[-\lambda,\lambda]}\F\mathbf 1_{[-\lambda,\lambda]}$ is a compact operator on $L^2[-\lambda,\lambda]$, symmetric on even functions. With $x=\lambda z$, $y=\lambda w$ one has $e^{2\pi ixy}=e^{i\gamma zw}$, and its eigenfunctions are the angular prolate spheroidal functions of Meixner--Sch\"afke \cite{MS}, $z\mapsto\Ps_n(z;\gamma^2)$, which are the functions $\PS(n,0,\gamma,z)$ of Mathematica, normalized by
\begin{equation}\label{eq:MSnorm}
\int_{-1}^1\Ps_n(z;\gamma^2)^2\,dz=\frac{2}{2n+1}\quad\text{for all }\gamma,\qquad \Ps_n(z;0)=P_n(z).
\end{equation}
For even $n$ we define
\begin{equation}\label{eq:hnl}
h_{n,\lambda}(x)=\Big(\frac{2n+1}{2}\Big)^{1/2}\lambda^{-1/2}\,\Ps_n\!\Big(\frac x\lambda;\gamma^2\Big)\ \ (\vert x\vert\le\lambda),\qquad h_{n,\lambda}(x)=0\ \ (\vert x\vert>\lambda),
\end{equation}
so that the $h_{n,\lambda}$ form an orthonormal basis of $L^2[-\lambda,\lambda]^{ev}$, with $h_{n,\lambda}(0)>0$ for $n\equiv0$ (mod $4$) and $h_{n,\lambda}(0)<0$ for $n\equiv2$ (mod $4$) (the sign of $P_n(0)$). We write
\begin{equation}\label{eq:chi}
\hat h_{n,\lambda}=\chi_n(\lambda)\,h_{n,\lambda}\quad\text{on }[-\lambda,\lambda],\qquad \delta_n(\lambda):=1-\vert\chi_n(\lambda)\vert ,
\end{equation}
and note that $\chi_n(\lambda)>0$ for $n\equiv0$ (mod $4$) and $\chi_n(\lambda)<0$ for $n\equiv2$ (mod $4$), since the Hermite limit $h_n$ of $h_{n,\lambda}$ has Fourier eigenvalue $i^n$. The operator $F_\lambda^*F_\lambda$ is the sinc kernel $\sin(\gamma(z-\tau))/(\pi(z-\tau))$ on $[-1,1]$, whose eigenvalues are $\lambda_n(\gamma)=\chi_n(\lambda)^2$; they are also the eigenvalues of $P_\lambda\widehat P_\lambda P_\lambda$ on ${\rm Range}\,P_\lambda^{ev}$. Two exact identities will be used throughout:
\begin{align}
&\int_{\vert y\vert>\lambda}\vert\hat h_{n,\lambda}(y)\vert^2dy=1-\chi_n(\lambda)^2,\qquad \int_{\vert y\vert>\lambda}\hat h_{n,\lambda}\hat h_{m,\lambda}=0\ \ (n\ne m),\label{eq:leak}\\
&\int_\R h_{n,\lambda}=\hat h_{n,\lambda}(0)=\chi_n(\lambda)h_{n,\lambda}(0).\label{eq:intlambda}
\end{align}
The first follows from Plancherel and $\int_{\vert y\vert\le\lambda}\hat h_{n,\lambda}\hat h_{m,\lambda}=\chi_n\chi_m\langle h_{n,\lambda},h_{m,\lambda}\rangle$: it says that $1-\chi_n^2$ is the $L^2$ mass of $\hat h_{n,\lambda}$ outside the band, the \emph{leakage}. Numerical values of $\chi_0,\chi_4$ are given in Appendix \ref{app:prolate}; for instance at $\lambda^2=50$, $1-\chi_4(\lambda)\approx1.36\cdot10^{-259}$.

\subsection{The uncertainty operator}

\begin{thm}\label{thm:unc}
Let $U_\lambda=2-P_\lambda-\widehat P_\lambda$ on $L^2(\R)^{ev}$. For every even $n$ the two-dimensional space $V_n={\rm span}\{h_{n,\lambda},\hat h_{n,\lambda}\}$ is invariant, and
\begin{equation}\label{eq:unceig}
\begin{aligned}U_\lambda\big(h_{n,\lambda}+{\rm sgn}(\chi_n)\hat h_{n,\lambda}\big)&=\delta_n(\lambda)\big(h_{n,\lambda}+{\rm sgn}(\chi_n)\hat h_{n,\lambda}\big),\\
U_\lambda\big(h_{n,\lambda}-{\rm sgn}(\chi_n)\hat h_{n,\lambda}\big)&=(2-\delta_n(\lambda))\big(h_{n,\lambda}-{\rm sgn}(\chi_n)\hat h_{n,\lambda}\big).\end{aligned}
\end{equation}
The orthogonal complement of $\bigoplus_nV_n$ is the Sonin space $\Qer P_\lambda\cap\Qer\widehat P_\lambda$, on which $U_\lambda=2$. Hence the spectrum of $U_\lambda$ below $1$ consists exactly of the simple eigenvalues $\delta_n(\lambda)$, $n$ even, with eigenvectors the self-dual vectors $h_{n,\lambda}+\hat h_{n,\lambda}$ for $n\equiv0$ (mod $4$) and the anti-self-dual vectors $h_{n,\lambda}-\hat h_{n,\lambda}$ for $n\equiv2$ (mod $4$).
\end{thm}
\begin{proof}
${\rm Range}\,P_\lambda\cap{\rm Range}\,\widehat P_\lambda=0$ since a band-limited function is entire. ${\rm Range}\,P_\lambda\cap\Qer\widehat P_\lambda=0$ since $\widehat P_\lambda f=0$ means that $\hat f$ vanishes on $[-\lambda,\lambda]$, while $\hat f$ is entire for compactly supported $f$; by symmetry $\Qer P_\lambda\cap{\rm Range}\,\widehat P_\lambda=0$. By Halmos' theorem on two projections, $L^2(\R)^{ev}=(\Qer P_\lambda\cap\Qer\widehat P_\lambda)\oplus\cH_{\rm gen}$ where on $\cH_{\rm gen}$ the pair is in generic position, and $\cH_{\rm gen}$ is the orthogonal sum of the two-dimensional spaces spanned by an eigenvector $e$ of $P\widehat PP$ in ${\rm Range}\,P$ and the unit vector $f=\widehat Pe/\Vert\widehat Pe\Vert$, with $\langle e,f\rangle=\cos\theta$, $\cos^2\theta$ the eigenvalue. Here $e=h_{n,\lambda}$, $\widehat P_\lambda h_{n,\lambda}=\F(\mathbf 1_{[-\lambda,\lambda]}\hat h_{n,\lambda})=\chi_n\hat h_{n,\lambda}$, so $f={\rm sgn}(\chi_n)\hat h_{n,\lambda}$ and $\cos\theta=\vert\chi_n\vert$. On such a pair $Pf=\cos\theta\,e$ and $\widehat Pe=\cos\theta\,f$, hence $(P+\widehat P)(e\pm f)=(1\pm\cos\theta)(e\pm f)$, which is \eqref{eq:unceig}. On the Sonin space $\Qer P_\lambda\cap\Qer\widehat P_\lambda$ both projections vanish, so $U_\lambda=2$ there.
\end{proof}

\begin{rem}
The Rayleigh quotient of $U_\lambda$ on a function $f\in{\rm Range}\,P_\lambda$ is its leakage: $\langle U_\lambda f,f\rangle=\langle(1-\widehat P_\lambda)f,f\rangle=\int_{\vert y\vert>\lambda}\vert\hat f\vert^2$. The vectors $h_{n,\lambda}+\hat h_{n,\lambda}$ have norm $\sqrt{2+2\chi_n}\approx2$ and are not time-limited, but their restriction to $[-\lambda,\lambda]$ is $(1+\chi_n)h_{n,\lambda}$.
\end{rem}

\subsection{The Fuchs ladder}

For fixed $n$ and $\gamma\to\infty$, Fuchs \cite{Fu} (see also \cite{Sle}) proved
\begin{equation}\label{eq:fuchs}
1-\lambda_n(\gamma)=A_n\gamma^{n+1/2}e^{-2\gamma}\big(1+O(\gamma^{-1})\big),\qquad A_n=\frac{4\sqrt\pi\,8^n}{n!}.
\end{equation}
Since $1-\lambda_n=(1-\vert\chi_n\vert)(1+\vert\chi_n\vert)=2\delta_n-\delta_n^2$, this gives \eqref{eq:fuchsintro}, and we record the values used below:
\begin{equation}\label{eq:deltas}
\delta_4=\frac{2^{14}\sqrt2\,\pi^5}{3}\lambda^9e^{-4\pi\lambda^2},\quad
\delta_6=\frac{2^{21}\sqrt2\,\pi^7}{45}\lambda^{13}e^{-4\pi\lambda^2},\quad
\delta_8=\frac{2^{26}\sqrt2\,\pi^9}{315}\lambda^{17}e^{-4\pi\lambda^2},
\end{equation}
\begin{equation}\label{eq:ratio}
\frac{\delta_{n+2}}{\delta_n}=\frac{256\pi^2\lambda^4}{(n+1)(n+2)},
\end{equation}
all up to the factor $1+O(\lambda^{-2})$. In logarithmic scale the ladder is $-\log\delta_n(\lambda)=4\pi\lambda^2-(2n+1)\log\lambda-\ell_n$, with $\ell_4=\log(2^{29/2}\pi^5/3)=14.6757$; the rungs are separated by $\log(\delta_{n+2}/\delta_n)=4\log\lambda+\log\frac{256\pi^2}{(n+1)(n+2)}$, so that for $n\equiv0$ (mod $4$)
\begin{equation}\label{eq:gaps}
\log\frac{\delta_{4j+2}}{\delta_{4j}} = 4\log\lambda + \log\frac{256\pi^2}{(4j+1)(4j+2)},\quad
\log\frac{\delta_{4j+4}}{\delta_{4j}} = 8\log\lambda + \log\frac{2^{16}\pi^4}{\prod_{k=1}^4(4j+k)};
\end{equation}
for $j=1$ the second constant is $\log(4096\pi^4/105)=8.2425$ and the first is $\log(128\pi^2/15)=4.4342$. The ladder is exact in the sense of Theorem \ref{thm:unc}; Fuchs' formula only describes its rungs for fixed $n$. For $n$ growing with $\lambda$, with $\sigma_n\in(0,1)$ the turning point of the prolate equation (fixed to leading order by $2n+1=\frac{2\gamma}{\pi}J(\sigma_n)$, $J(\sigma)=2[E(\sigma)-(1-\sigma^2)K(\sigma)]$), one has $\log(1/\delta_n)=2\gamma I(\sigma_n)+O(\log\gamma)$ with $I(\sigma)=E(\sigma')-\sigma^2K(\sigma')$, $\sigma'=\sqrt{1-\sigma^2}$, uniformly for $\sigma_n\le\sigma_0<1$ \cite{Du}; the rungs stay exponentially small up to $n\approx2\gamma/\pi=4\lambda^2$, Slepian's count of the eigenvalues of the sinc kernel close to $1$, so that each of the two sectors contains about $\lambda^2$ minuscule eigenvalues (exactly $\lambda^2$ and $\lambda^2-1$ below $10^{-2}$ in the computed range, Fact \ref{fact:count}). The rest of the spectrum of $U_\lambda$ consists of the eigenvalues $2-\delta_n(\lambda)$, accumulating at $2$ from below, and of the eigenvalue $2$ itself, of infinite multiplicity on the Sonin space.

\subsection{The zero-integral constraint and the Hermite limit}\label{ss:hermite}

The map $\cE$ requires $\hat f(0)=\int f=0$, otherwise $\cE(f)(u)\sim u^{-1/2}\int_0^\infty f$ as $u\to0$ is not in $L^2(\dstar)$; this constraint has no counterpart on the prolate side and is what removes the lowest rung $\delta_0$ from the Weil spectrum. The relevant vectors are therefore the zero-integral combinations of prolate functions of a given sector. In the even sector the first one is the two-mode state of \cite{CCM},
\begin{equation}\label{eq:twomode}
h_\lambda=\Big(\int_\R h_{0,\lambda}\Big)h_{4,\lambda}-\Big(\int_\R h_{4,\lambda}\Big)h_{0,\lambda}=\chi_0h_{0,\lambda}(0)\,h_{4,\lambda}-\chi_4h_{4,\lambda}(0)\,h_{0,\lambda},
\end{equation}
by \eqref{eq:intlambda}; it is bilinear in $(h_{0,\lambda},h_{4,\lambda})$, so that no sign or normalization convention for the prolate functions enters in the ratio $h_\lambda/h_\lambda(1)$. Its normalized form is $a_\lambda h_{0,\lambda}+b_\lambda h_{4,\lambda}$ with $a_\lambda^2+b_\lambda^2=1$. For fixed $n$, Meixner--Sch\"afke (\cite{MS}, Satz 9, p.~243; \cite{CCM}, Lemma 7.2) give $\max_{\vert x\vert\le\lambda}\vert h_{n,\lambda}-h_n\vert\le C_n\lambda^{-2}$, $h_n$ the normalized Hermite function, so that $a_\lambda\to-3/\sqrt{33}$, $b_\lambda\to\sqrt{8/11}$ and $h_\lambda/\Vert h_\lambda\Vert$ converges to the vector $h$ of \cite{CCM}, Lemma 7.1, which is a multiple of
\begin{equation}\label{eq:h}
h(x)=\frac\pi2\,x^2(2\pi x^2-3)\,e^{-\pi x^2}:
\end{equation}
indeed $h$ is self-dual, $\hat h=h$, so that $\int h=\hat h(0)=h(0)$ and the zero-integral constraint is the absence of a constant term. Two consequences of \eqref{eq:leak}, \eqref{eq:intlambda}:
\begin{equation}\label{eq:h0leak}
h_\lambda(0)=a_\lambda\delta_0h_{0,\lambda}(0)+b_\lambda\delta_4h_{4,\lambda}(0)=O(\delta_4),\qquad \int_{y>\lambda}\vert\hat h_\lambda\vert^2=a_\lambda^2\delta_0+b_\lambda^2\delta_4+O(\delta_4^2)=\frac8{11}\delta_4\big(1+O(\lambda^{-2})\big),
\end{equation}
using $\delta_0/\delta_4=\frac{3}{2^9(2\pi)^4}\lambda^{-8}$. Thus the one-sided leakage of the zero-integral two-mode state is $\frac8{11}\delta_4$, and its Rayleigh quotient for $U_\lambda$, which is the full leakage $\int_{\vert y\vert>\lambda}\vert\hat h_\lambda\vert^2=2\int_{y>\lambda}\vert\hat h_\lambda\vert^2$, is $\frac{16}{11}\delta_4$, against $2\delta_4$ for $h_{4,\lambda}$ alone and $\delta_4$ for the self-dual pair vector $h_{4,\lambda}+\hat h_{4,\lambda}$ (the one-sided coefficient $\frac8{11}$ is the $n=1$ term of the constant $C_*$ of \eqref{eq:Cstar} below): the compression of $U_\lambda$ to the hyperplane $\{\int f=0\}$ of each sector is the operator whose spectrum the Weil form should be compared with.

More generally let $H_K(\lambda)=\{\sum_{k\le K}c_kh_{4k,\lambda}:\int h=0\}$, a space of dimension $K$ spanned by the two-mode states $\xi(\lambda,k)$ built on $(h_{0,\lambda},h_{4k,\lambda})$, and $H_K^{\rm Herm}$ its Hermite limit. The Mellin transform of $\cE$ on $H_K^{\rm Herm}$ is computable:

\begin{prop}\label{prop:xipoly}
For $h\in H_K^{\rm Herm}$, $\cM\cE(h)(s)=\xi(s+\tfrac12)\,q_h(s^2)$ with $q_h$ a polynomial of degree $\le K-1$, and $h\mapsto q_h$ is a linear bijection of $H_K^{\rm Herm}$ onto the polynomials of degree $<K$. Consequently
\begin{equation}\label{eq:gramherm}
\Vert\cE(h)\Vert^2_{L^2(\rs,\dstar)}=\frac1{2\pi}\int_\R\Xi(t)^2\,\vert q_h(-t^2)\vert^2dt ,
\end{equation}
and $g_K^{\rm Herm}:=\inf\{\Vert\cE(h)\Vert^2:h\in H_K^{\rm Herm},\Vert h\Vert_2=1\}>0$.
\end{prop}
\begin{proof}
With $\Gamma_\R(w)=\pi^{-w/2}\Gamma(w/2)$ and $\int_0^\infty x^{2j}e^{-\pi x^2}x^{w-1}dx=\frac12\Gamma_\R(w+2j)=\frac12\pi^{-j}\Gamma_\R(w)\prod_{i<j}(\frac w2+i)$, the Mellin transform $\int_0^\infty h_n(x)x^{w-1}dx$ of a Hermite function is $\Gamma_\R(w)P_n(w)$ with $P_n$ a polynomial of degree $n/2$, so by \eqref{eq:mellinE} $\cM\cE(h)(s)=\Lambda(w)P_h(w)$, $w=s+\frac12$, $\Lambda=\Gamma_\R\zeta$, $\deg P_h\le2K$. Since $\hat h=h$ and the Mellin transform of $\hat h$ at $w$ is $\frac{\Gamma_\R(1-w)}{\Gamma_\R(w)}$ times that of $h$ at $1-w$, the functional equation of $\Lambda$ forces $P_h(w)=P_h(1-w)$; the constraint $\int h=0$ is $P_h(1)=0$, hence also $P_h(0)=0$, and $P_h(w)=\frac12w(w-1)Q_h(w)$ with $Q_h$ symmetric under $w\mapsto1-w$, i.e. a polynomial $q_h$ in $(w-\frac12)^2=s^2$ of degree $\le K-1$. Then $\cM\cE(h)=\xi(w)Q_h(w)$ with $\xi(w)=\frac12w(w-1)\Lambda(w)$. The map $h\mapsto q_h$ is injective since $\xi\not\equiv0$, hence bijective by dimension; \eqref{eq:gramherm} is Mellin--Plancherel with $\xi(\frac12+it)=\Xi(t)$, and positivity follows since $\Xi^2>0$ off a discrete set and $q_h\ne0$, by compactness of the unit sphere.
\end{proof}

\begin{figure}[ht]\centering
\includegraphics[width=0.8\textwidth]{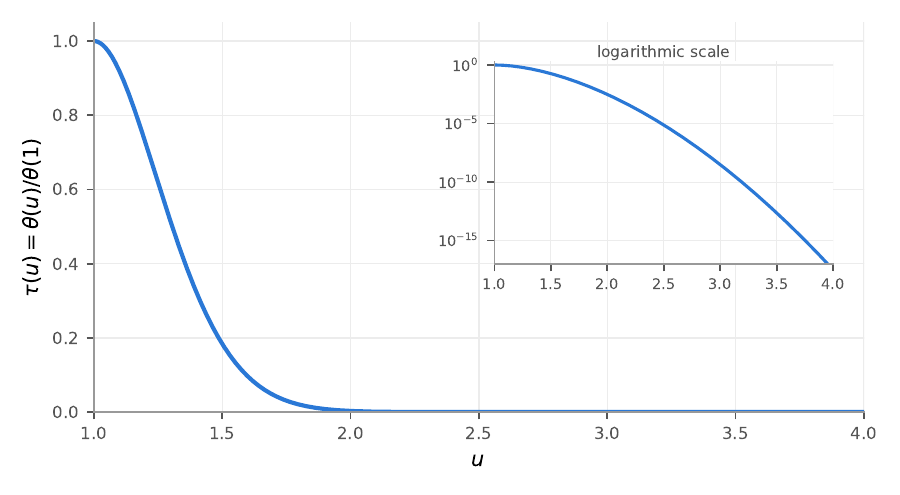}
\caption{The target $\tau=\theta/\theta(1)$ on $[1,4]$ ($\tau(u)=\tau(1/u)$), with $\theta(1)=0.2233484\ldots$; inset: the same on a logarithmic scale, showing the Gaussian decay $\tau(u)\approx\theta(1)^{-1}u^{1/2}h(u)$ for $u\gtrsim2$.}\label{fig:target}
\end{figure}

In the basis $s^{2j}$ the Gram matrix \eqref{eq:gramherm} is the Hankel matrix of the even moments $\frac1{2\pi}\int\Xi(t)^2t^{2j}dt$. For $K=1$ one gets the target:

\begin{thm}\label{thm:mellintheta}
Let $\theta=\cE(h)$ with $h$ as in \eqref{eq:h}, i.e. $\theta$ is the function \eqref{eq:targetintro}. Then
\begin{equation}\label{eq:mellintheta}
\cM\theta(s)=\int_0^\infty\theta(u)u^{s}\dstar=\frac14\,\xi\Big(s+\frac12\Big)\qquad(s\in\C),
\end{equation}
$\theta\circ\iota=\theta$ exactly, $\theta>0$ on $(0,\infty)$ with its maximum at $u=1$, and
$$\theta(1)=0.223348\ldots,\qquad\int_0^\infty\theta\,\dstar=\frac14\xi(\tfrac12)=0.124280\ldots,\qquad\Vert\theta\Vert^2_{L^2(\rs,\dstar)}=0.0199845\ldots$$
\end{thm}
\begin{proof}
$\int_0^\infty h(x)x^{w-1}dx=\frac\pi2\big[2\pi\cdot\frac12\pi^{-(w+4)/2}\Gamma(\frac w2+2)-3\cdot\frac12\pi^{-(w+2)/2}\Gamma(\frac w2+1)\big]=\frac14\pi^{-w/2}\Gamma(\frac w2)\frac w2\big[(w+2)-3\big]=\frac18w(w-1)\Gamma_\R(w)$, and \eqref{eq:mellinE} gives $\frac18w(w-1)\Gamma_\R(w)\zeta(w)=\frac14\xi(w)$. The symmetry is \eqref{eq:radical} and $\hat h=h$. Positivity: by \eqref{eq:h}, $h(x)>0$ for $x>(3/2\pi)^{1/2}=0.691\ldots$; for $u\ge1$ every argument $nu$ in $\theta(u)=u^{1/2}\sum_{n\ge1}h(nu)$ is $\ge1$, so all the terms are positive and $\theta(u)>0$, and for $u<1$ one uses $\theta(u)=\theta(1/u)$. The maximum at $u=1$, where $\theta'(1)=0$ by the symmetry, and the decrease of $\theta$ on $[1,\infty)$ are verified numerically (Figure \ref{fig:target}). The value of $\int\theta\dstar$ is $\frac14\xi(\frac12)=-\frac1{32}\pi^{-1/4}\Gamma(\frac14)\zeta(\frac12)$.
\end{proof}

\section{The dictionary}\label{sec:dict}

\subsection{Sectors and the window}

Combining \eqref{eq:radical} with Theorem \ref{thm:unc}: if $f\in\sr0$ is supported in $[-\lambda,\lambda]$ and close to self-dual, $\hat f\approx f$, then $\cE(f)$ vanishes for $u>\lambda$ and $\cE(f)=\cE(\hat f)\circ\iota\approx\cE(f)\circ\iota$ vanishes approximately for $u<\lambda^{-1}$: it is approximately supported in the window and approximately $\iota$-even. If instead $\hat f\approx-f$, then $\cE(f)$ is approximately $\iota$-odd. The correspondence is therefore
$$\text{Fourier}\approx+1\ (n\equiv0\bmod4)\ \longleftrightarrow\ \iota=+1,\qquad \text{Fourier}\approx-1\ (n\equiv2\bmod4)\ \longleftrightarrow\ \iota=-1,$$
and both are the same symmetry transported by $\cE$. The precise statement is the Poisson defect: for $f\in\sr0$ supported in $[-\lambda,\lambda]$,
\begin{equation}\label{eq:poissondefect}
\cE(f)(u)-\cE(f)(1/u)=\cE(\hat f-f)(1/u),
\end{equation}
so that the departure from $\iota$-symmetry on the window is measured by $\hat f-f$, i.e. by $\delta_n$ inside the band and by the leakage outside.

\subsection{The test vector and the upper bounds}

A function supported in $[-\lambda,\lambda]$ is not in $\sr0$ (it is not smooth at $\pm\lambda$), and the prolate functions do not vanish at $\pm\lambda$: $h_{n,\lambda}(\lambda)^2=4\pi\lambda\delta_n(\lambda)(1+O(\lambda^{-2}))$, an exponentially small but non-zero value with logarithmic derivative $-2\pi^2\lambda^3$ (the regular singular point of the prolate equation). The construction of \cite{K,CK} repairs this: $h_\lambda$ is multiplied by an even Gevrey cutoff of class $G^s$, $s>1$, equal to $1$ except on a layer of width $\lambda^{-3}$ at the endpoints, and corrected by two fixed Gevrey bumps so as to restore $f(0)=\int f=0$; the result $f_\lambda\in\sr0$ satisfies $\Vert f_\lambda-h_\lambda\Vert_2\le C\lambda^{-1}\delta_4^{1/2}$, and $K_\lambda:=\cE(f_\lambda)\in\Rad(QW)$ exactly. The test vector is
\begin{equation}\label{eq:testvector}
k_\lambda:=\chi_\lambda K_\lambda,\qquad r_\lambda:=(1-\chi_\lambda)K_\lambda ,
\end{equation}
with $\chi_\lambda$ a multiplicative Gevrey cutoff vanishing for $u\le\lambda^{-1}$ and equal to $1$ for $u\ge\lambda^{-1}e^{\lambda^{-4}}$, so that $k_\lambda$ is smooth and supported in $[\lambda^{-1},\lambda]$. The mechanism is the identity
\begin{equation}\label{eq:bridge}
QW_\lambda(k_\lambda,k_\lambda)=QW(r_\lambda,r_\lambda),
\end{equation}
which follows from $K_\lambda\in\Rad(QW)$ and the symmetry of the form: the energy of the test vector is the energy of the lower tail, which is the Poisson defect \eqref{eq:poissondefect} read on $u<\lambda^{-1}$, i.e. the alias sum $u^{-1/2}\sum_{n\ge1}\hat f_\lambda(n/u)$. Its $n=1$ term has $L^2$ norm equal to the leakage $\int_{y>\lambda}\vert\hat f_\lambda\vert^2\approx\frac8{11}\delta_4$, the Rayleigh quotient of $U_\lambda$ on $h_\lambda$, and the full tail has
\begin{equation}\label{eq:Cstar}
\Vert r_\lambda\Vert^2_{L^2(\dstar)}=(C_*+o(1))\,\delta_4(\lambda),\qquad C_*=\frac{16}{11\pi}\sum_{n\ge1}\frac{\arcsin(1/n)}{n}=1.03293\ldots
\end{equation}
(\cite{CK}; the higher aliases contribute $42\%$ of the mass and the interference terms are of lower order). Estimating the three parts of the explicit formula on $r_\lambda$ by absolute values gives:

\begin{thm}[\cite{K,CK}]\label{thm:upper}
For every $s>1$ there are $C_s$, $c_0>0$ such that for $\lambda\ge\lambda_0(s)$: $\Vert k_\lambda\Vert\ge c_0$, $\Vert k_\lambda-\mathbf 1_{[\lambda^{-1},\lambda]}\cE(h_\lambda)\Vert\le C_s\lambda^{-1}\delta_4^{1/2}\log\lambda$, and
\begin{equation}\label{eq:upper}
\frac{\vert QW_\lambda(k_\lambda,k_\lambda)\vert}{\Vert k_\lambda\Vert^2}\le C_s\,\delta_4(\lambda)(\log\lambda)^3 .
\end{equation}
In particular $\epsilon_\lambda\le C_s\delta_4(\lambda)(\log\lambda)^3$. The same bound holds for the Rayleigh quotient of the finite-dimensional truncation of $QW_\lambda$ to the first $N(\lambda)=O(\lambda^{4+2s}\log\lambda)$ log-Fourier modes.
\end{thm}

The construction is linear in the prolate vector, and all the estimates behind \eqref{eq:upper} only use the leakage, the endpoint values and the uniform asymptotics of \cite{Du}, which hold uniformly for $\sigma_n\le\sigma_0<1$. The method of \cite{CK} therefore gives the following uniform version; it is not stated in \cite{CK} in this form, and we indicate below the only point where the argument has to be followed for several prolate functions at once.

\begin{thm}[after \cite{CK}]\label{thm:upperK}
Let $E_\lambda:f\mapsto\chi_\lambda\cE(f^{\rm Gev})$ denote the construction above (Gevrey repair, $\cE$, cutoff), which is linear. For $0<\sigma_0<1$ there is $C=C(s,\sigma_0)$ such that for all $K$ with $\sigma_{4K}\le\sigma_0$ and all $f\in H_K(\lambda)$,
\begin{equation}\label{eq:formcomp}
\vert QW_\lambda(E_\lambda f,E_\lambda f)\vert\le C\,K\,(\log\lambda)^3\,\langle U_\lambda f,f\rangle\le C\,K\,\delta_{4K}(\lambda)(\log\lambda)^3\,\Vert f\Vert_2^2 ,
\end{equation}
and for fixed $K$, with $g_K^{\rm Herm}$ as in Proposition \ref{prop:xipoly}, $\epsilon_\lambda(j)\le C_K\,\delta_{4j}(\lambda)(\log\lambda)^3$ for $j\le K$ and $\lambda\ge\lambda_0(K)$. The analogous statements hold in the odd sector with $h_{4k+2,\lambda}$ and $\delta_{4j+2}$.
\end{thm}
Indeed, in \cite{CK} the estimates are carried out for the two-mode vector $h_\lambda$, but each of them (endpoint values, layer, Fourier tail) is expressed through the single quantity $S(f)=\sum_k\vert c_k\vert\delta_{4k}^{1/2}$, and $S(f)^2\le K\sum_kc_k^2\delta_{4k}\le K\langle U_\lambda f,f\rangle$ since $\langle U_\lambda f,f\rangle=\sum_kc_k^2(1-\chi_{4k}^2)$ on $H_K(\lambda)$, so that they apply verbatim on $H_K(\lambda)$, with the constant multiplied by $K$.

\begin{rem}
The exponent $3$ of the logarithm is expected to improve to $1$: the archimedean term of the explicit formula contributes one genuine $\log\lambda$ (it is the multiplier $\log\vert\tau\vert$ acting on a tail oscillating at log-frequency $\lambda^2$), the prime sum contributes $\log(Y_\lambda/\lambda)$ with $Y_\lambda\asymp\lambda^{3+2s}$ the Gevrey scale, and the pole term is $O(\lambda^{-4}\delta_4)$ once the oscillation of $\hat f_\lambda$ is taken into account \cite{CK}.
\end{rem}

\subsection{The hypothesis}

\begin{hypo}[logarithmic agreement with the Fuchs ladder]\label{hyp:LA}
There are constants $0<c\le C$ such that, for all $\lambda\ge\lambda_0$ and $1\le j\le J$,
\begin{equation}\label{eq:LA}
c\,\delta_{4j}(\lambda)\le\epsilon_\lambda(j)\le C\,\delta_{4j}(\lambda)\log\lambda,\qquad c\,\delta_{4j+2}(\lambda)\le\epsilon^{\rm odd}_\lambda(j)\le C\,\delta_{4j+2}(\lambda)\log\lambda .
\end{equation}
\end{hypo}

The upper halves are Theorem \ref{thm:upperK} (with the exponent $3$ at present); the lower halves are the content of the hypothesis, which says that nothing in the Weil spectrum is smaller than the uncertainty operator allows. It is a statement about constants, not only about exponents: the estimates of Section \ref{sec:ground} use the gap $\epsilon_\lambda(2)/\epsilon_\lambda(1)\asymp\delta_8/\delta_4\asymp\lambda^8$ against polynomial losses, and a weaker hypothesis allowing factors $\lambda^{O(1)}$ in \eqref{eq:LA} would not suffice. The sharp numerical test is therefore the difference
\begin{equation}\label{eq:gaptest}
\log\epsilon_\lambda(2)-\log\epsilon_\lambda(1)\ \text{ against }\ \log\frac{\delta_8}{\delta_4}=8\log\lambda+8.2425 ,
\end{equation}
and likewise $\log\epsilon^{\rm odd}_\lambda(1)-\log\epsilon_\lambda(1)$ against $4\log\lambda+4.4342$. Two further predictions: the interlacing $\epsilon_\lambda(j)<\epsilon^{\rm odd}_\lambda(j)<\epsilon_\lambda(j+1)$, and the count of about $\lambda^2$ minuscule eigenvalues in each sector.

\begin{rem}[why a lower bound is of a different nature]
An upper bound for $\epsilon_\lambda$ is produced by exhibiting a test vector; a lower bound $\epsilon_\lambda\ge c\delta_4$ implies in particular $QW_\lambda\ge0$, which for all $\lambda$ is equivalent to the Riemann Hypothesis. Hypothesis \ref{hyp:LA} is thus strictly stronger than Weil positivity on the window, and it is tested, not proved, in this paper.
\end{rem}

\section{Numerical evidence for the spectrum}\label{sec:numspec}

\subsection{Method}\label{ss:method}
The matrix of $QW_\lambda$ is computed in the orthonormal basis $b_j(u)=\xi_j(\log u)$, $\vert j\vert\le N$, of $E_N\subset L^2([\lambda^{-1},\lambda],\dstar)$ described in Appendix \ref{app:matrix}; the even sector is spanned by $b_0,\dots,b_N$ and the odd sector by $b_{-1},\dots,b_{-N}$. All entries are given in closed form (digamma and trigamma values and series in $\lambda^{-4}$), so that the matrix can be produced at any $N$ and any precision; at $\lambda^2=13$, $N=120$ it agrees with an independent 3386-bit computation to $8\cdot10^{-109}$ in every entry. The values below are for
$$\lambda^2\in\{5,9,22,24,28,32,36,40,50\},\qquad N=960 ,$$
with $N=240$ and $N=480$ as convergence checks (and $N$ up to $1920$ for $\lambda^2=13$), together with the larger values $\lambda^2=100$, computed at $N=1000$ (all quantities) and $N=2000$ (lowest eigenvalues and conditioning), with $N=500$ as a check, and $\lambda^2=150$, computed at $N=1500$ (all quantities) and $N=2500$ (lowest eigenvalues, ground state and conditioning), with $N=750$ as a check.

\emph{Precision.} Since the minuscule eigenvalues span hundreds of orders of magnitude, the relevant condition is not the size of $\epsilon_\lambda(1)$ alone: an eigenpair $\epsilon_\lambda(k)$ is obtained accurately by inverse iteration only if the working precision $2^{-b}$ satisfies $2^{-b}\ll\epsilon_\lambda(1)^2/(\Vert QW^N_\lambda\Vert\,\epsilon_\lambda(k))$. We used $b\approx3.33\,(2\vert\log_{10}\epsilon_\lambda(1)\vert+40)$ bits, from $300$ bits at $\lambda^2=5$ to $1880$ bits at $\lambda^2=50$ ($\epsilon_\lambda(1)\approx10^{-258}$) $3700$ bits at $\lambda^2=100$ ($\epsilon_\lambda(1)\approx10^{-529}$) and $5500$ bits at $\lambda^2=150$ ($\epsilon_\lambda(1)\approx10^{-801}$). When only the $k$ lowest eigenpairs are wanted, the condition involves the eigenvalue $\epsilon_\lambda(k+6)$ reached by the iteration block instead of $10^{-8}$, and $2400$ bits suffice at $\lambda^2=100$, $3400$ at $\lambda^2=150$; at $N=1000$ they reproduce the $3700$-bit values to $30$ digits.

\emph{Eigensolver.} All eigenpairs below $10^{-8}$ are computed in ball arithmetic (Arb) by block inverse subspace iteration with a Rayleigh--Ritz step on the matrix itself; the larger eigenvalues, and the counts of \S\ref{ss:counting}, come from the full spectrum in double precision. The relative residuals $\Vert QW\,v-\epsilon v\Vert/\epsilon$ are below $10^{-20}$ except in the even sector at $\lambda^2=40,50$, where they are at most $10^{-11}$ and $3\cdot10^{-4}$; since the Rayleigh quotient error is bounded by $\Vert r\Vert^2/{\rm gap}$, every eigenvalue reported is still correct to at least fifteen digits. At $\lambda^2=100$ the inverse of the sector matrix has condition number $\approx10^{529}$ and the ball radii produced by Arb for it are meaningless; the iteration is then carried out on midpoints, the accuracy being controlled by the precision rule and the residuals: the iteration converges to a relative change of $10^{-34}$--$10^{-40}$ in four or five steps, the relative residuals that can be evaluated in double precision are below $10^{-20}$, and those of the lowest eigenpairs at $N=2000$ are below $10^{-63}$. The midpoint version reproduces the earlier values to $30$ digits (checked at $\lambda^2=13$ and $50$). At $\lambda^2=150$ the explicit inverse is replaced by an approximate LU solve at each step, which halves the memory and gives the same eigenvalues. For all the minuscule eigenvalues at $\lambda^2=150$, $N=2500$ the required precision ($5500$ bits) exceeds the available memory, and the eigenvalues are computed in batches $[a,b]$ of consecutive indices by block inverse iteration with the shifted matrix $(QW^N_\lambda+\sigma_a)^{-1}$, $\sigma_a\approx\delta_{4a}(\lambda)$, the eigenvectors of the previous batches being projected out. The eigenvalues below the batch are then mapped to $\approx\sigma_a^{-1}$ instead of being amplified, so that each batch needs only $2^{-b}\ll\epsilon_\lambda(a)^2/(\Vert QW^N_\lambda\Vert\,\epsilon_\lambda(b+6))$; $3200$ bits suffice for seven batches. At $\lambda^2=50$, $N=960$ this reproduces all the eigenvalues to $4\cdot10^{-15}$, and at $\lambda^2=150$, $N=2500$ the two lowest ones to all the digits computed before.

\emph{Convergence in $N$.} The eigenvalues decrease with $N$ (nested subspaces). The truncation is the only source of error at the precision used, and it is visible:
\begin{table}[ht]\centering\small
\begin{tabular}{c|ccc|ccc}
$\lambda^2$ & \multicolumn{3}{c|}{$-\log\epsilon_\lambda(1)$} & \multicolumn{3}{c}{$-\log\epsilon^{\rm odd}_\lambda(1)$}\\
 & $N=240$ & $N=480$ & $N=960$ & $N=240$ & $N=480$ & $N=960$\\ \hline
5 & 39.219 & 39.229 & 39.231 & 32.214 & 32.240 & 32.256 \\
9 & 86.488 & 86.553 & 86.587 & 78.284 & 78.297 & 78.307 \\
22 & 245.456 & 245.653 & 245.717 & 235.223 & 235.415 & 235.508 \\
24 & 270.101 & 270.338 & 270.438 & 259.719 & 259.923 & 260.010 \\
28 & 319.513 & 319.833 & 319.955 & 308.743 & 309.102 & 309.225 \\
32 & 368.795 & 369.331 & 369.480 & 357.820 & 358.328 & 358.493 \\
36 & 417.416 & 419.091 & 419.276 & 406.293 & 407.802 & 408.009 \\
40 & 460.038 & 468.788 & 469.002 & 449.022 & 457.347 & 457.571 \\
50 & 542.045 & 593.242 & 593.666 & 531.107 & 581.322 & 581.736 \\
$100^\ast$ & 1066.141 & 1218.065 & 1218.627 & 1053.997 & 1204.777 & 1205.361 \\
$150^\dagger$ & 1562.160 & 1844.115 & 1845.018 & 1549.356 & 1829.997 & 1830.872 \\
\end{tabular}
\caption{Convergence of the lowest even and odd eigenvalues in the truncation $N$. $^\ast$For $\lambda^2=100$ the three columns are $N=500,1000,2000$; $^\dagger$for $\lambda^2=150$, $N=750,1500,2500$.}\label{tab:conv}
\end{table}
The increment from $N=480$ to $N=960$ grows from $0.002$ ($\lambda^2=5$) to $0.42$ ($\lambda^2=50$) and decreases by a factor $0.3$--$0.4$ per doubling of $N$; for $\lambda^2\ge40$ the value $N=240$ is far from sufficient, and $N$ should exceed about $10\lambda^2$. The eigenvalues at $N=960$ are thus slight overestimates; extrapolating the increments geometrically, $-\log\epsilon_\lambda(1)$ should still increase by about $0.2$--$0.3$ beyond $N=960$ at $\lambda^2=50$, and by less for smaller $\lambda$. At $\lambda^2=100$ the same pattern appears at twice the scale: the step $N=500\to1000$ raises $-\log\epsilon_\lambda(1)$ by $152$, the step $N=1000\to2000$ by $0.56$, and about $0.3$ should remain beyond $N=2000$. (At $\lambda^2=100$, $N=250$ still gives $-\log\epsilon_\lambda(1)=767.7$, against $1218.6$ at $N=2000$: well below $N\approx10\lambda^2$ the eigenvalue is far from converged, although the ground state is already close to the prolate vector, with $\sin\angle=5\cdot10^{-3}$.) At $\lambda^2=150$ the steps $N=750\to1500\to2500$ raise $-\log\epsilon_\lambda(1)$ by $282$ and $0.90$. For $\lambda^2=100$ the tables below use $N=2000$ for the lowest eigenvalues and $N=1000$ for the counts, the ratios over all $j$, and the prolate angles; for $\lambda^2=150$, $N=2500$ and $N=1500$ respectively.

\subsection{The even sector}\label{ss:evensector}
Writing $D_j(\lambda)=-\log\delta_{4j}(\lambda)$ for the exact rungs, computed from the prolate functions as in Appendix \ref{app:prolate}, and $D_j^{\rm Fu}$ for Fuchs' asymptotic value \eqref{eq:fuchsintro}, one finds (rows marked ${}^\ast$: eigenvalues computed by A.L.G. at $N=3200$, $4000$, $4800=16\lambda^2$ for $\lambda^2=200$, $250$, $300$, with $\epsilon^{\rm odd}_\lambda(2)$ at $N=3000$ for $\lambda^2=250$; the exact rungs are ours):
\begin{equation}\label{eq:table2}\setlength{\arraycolsep}{3.4pt}
\begin{array}{c|ccccc|ccc}
\lambda^2 & D_1 & -\log\epsilon_\lambda(1) & D_2 & D_2^{\rm Fu} & -\log\epsilon_\lambda(2) & \log\frac{\epsilon_\lambda(2)}{\epsilon_\lambda(1)} & \log\frac{\delta_8}{\delta_4} & 8\log\lambda+8.2425\\ \hline
5 & 41.185 & 39.231 & 27.235 & 26.233 & 26.052 & 13.179 & 13.950 & 14.680 \\
9 & 88.680 & 86.587 & 72.023 & 71.503 & 70.762 & 15.824 & 16.657 & 17.031 \\
22 & 247.933 & 245.717 & 227.472 & 227.268 & 225.832 & 19.885 & 20.461 & 20.607 \\
24 & 272.669 & 270.438 & 251.848 & 251.661 & 250.230 & 20.208 & 20.822 & 20.955 \\
28 & 322.233 & 319.955 & 300.776 & 300.616 & 299.088 & 20.866 & 21.458 & 21.571 \\
32 & 371.892 & 369.480 & 349.886 & 349.747 & 348.110 & 21.370 & 22.007 & 22.105 \\
36 & 421.623 & 419.276 & 399.134 & 399.011 & 397.314 & 21.962 & 22.489 & 22.577 \\
40 & 471.411 & 469.002 & 448.492 & 448.381 & 446.766 & 22.236 & 22.919 & 22.998 \\
50 & 596.064 & 593.666 & 572.236 & 572.148 & 570.413 & 23.254 & 23.828 & 23.891 \\
100 & 1221.251 & 1218.627 & 1194.619 & 1194.575 & 1192.694 & 25.933 & 26.632 & 26.663 \\
150 & 1847.741 & 1845.018 & 1819.476 & 1819.447 & 1817.359 & 27.660 & 28.264 & 28.285 \\
200^\ast & 2474.762 & 2471.868 & 2445.342 & 2445.320 & 2443.068 & 28.800 & 29.420 & 29.436 \\
250^\ast & 3102.075 & 3099.148 & 3071.759 & 3071.742 & 3069.451 & 29.697 & 30.316 & 30.328 \\
300^\ast & 3729.573 & 3726.627 & 3698.525 & 3698.511 & 3696.170 & 30.457 & 31.047 & 31.058 \\
\end{array}
\end{equation}
\begin{enumerate}
\item The exact rungs $D_j$ and Fuchs' asymptotics differ by $O(\lambda^{-2})$ in relative terms but by $1.0$ ($\lambda^2=5$), $0.09$ ($\lambda^2=50$), $0.04$ ($\lambda^2=100$), $0.03$ ($\lambda^2=150$) and $0.02$--$0.015$ ($\lambda^2=200$--$300$) in absolute terms for $D_2$; the comparison with the Weil spectrum must use the exact $D_j$.
\item The second eigenvalue exceeds $\delta_8(\lambda)$ by a factor between $3.3$ ($\lambda^2=5$) and $6.2$ ($\lambda^2=50$), $6.9$ at $\lambda^2=100$, $8.3$ at $\lambda^2=150$, and $9.7$, $10.1$, $10.5$ at $\lambda^2=200$, $250$, $300$; the lowest one exceeds $\delta_4(\lambda)$ by a factor between $7.1$ and $11.2$ for $\lambda^2\le50$, $13.8$ at $\lambda^2=100$, $15.2$ at $\lambda^2=150$, and $18.1$, $18.7$, $19.0$ at $\lambda^2=200$, $250$, $300$ ($16.1$ at $\lambda^2=250$, $N=5000$). All these factors still decrease with $N$ (Fact \ref{fact:logbound}).
\item The test \eqref{eq:gaptest}: $\log(\epsilon_\lambda(2)/\epsilon_\lambda(1))$ falls short of the exact $\log(\delta_8/\delta_4)$ by $0.53$ to $0.83$, with no trend in $\lambda$ between $\lambda^2=5$ and $\lambda^2=300$ ($0.70$ at $\lambda^2=100$, $0.60$ at $\lambda^2=150$, $0.62$, $0.62$, $0.59$ at $\lambda^2=200$, $250$, $300$); against Fuchs' value $8\log\lambda+8.2425$ the shortfall is $0.60$--$1.50$, the larger values coming from the inaccuracy of Fuchs' asymptotics at small $\lambda$.
\end{enumerate}

\begin{fact}\label{fact:logbound}
The upper half of Hypothesis \ref{hyp:LA} holds for $j=1$ with the first power of $\log\lambda$ and a small constant: for every $\lambda$ of Table \ref{tab:logbound} ($5\le\lambda^2\le300$),
\begin{equation}\label{eq:logbound}
\epsilon_\lambda(1)\le 8.8\,\delta_4(\lambda)\log\lambda,\qquad\text{and}\qquad \epsilon_\lambda(1)\le 6.9\,\delta_4(\lambda)\log\lambda\quad(\lambda^2\ge22).
\end{equation}
Since $QW^N_\lambda$ is the compression of $QW_\lambda$ to $E_N$, the min-max principle gives $\epsilon^N_\lambda(1)\ge\epsilon_\lambda(1)$ for the lowest eigenvalue $\epsilon^N_\lambda(1)$ of $QW^N_\lambda$, so that every computed ratio $\epsilon^N_\lambda(1)/\delta_4(\lambda)$ is an upper bound for the true one, and \eqref{eq:logbound} involves no extrapolation in $N$. Table \ref{tab:logbound} gives the ratios at the largest $N$ computed.
\begin{enumerate}
\item The ratio $\epsilon_\lambda(1)/(\delta_4(\lambda)\log\lambda)$ stays between $4.9$ and $8.8$, without upward trend; the values at large $\lambda$ are the least converged in $N$ and hence the most overestimated.
\item The ratio $\epsilon_\lambda(1)/(\delta_4(\lambda)(\log\lambda)^3)$ decreases from $13.5$ to $0.8$: the exponent $3$ of Theorem \ref{thm:upperK} is far from sharp in this range.
\item Comparing all $\lambda$ at the same ratio $N/\lambda^2$ (interpolating in $\log N$), a least-squares fit $\epsilon^N_\lambda(1)/\delta_4(\lambda)\approx C(\log\lambda)^p$ gives $p=0.73$, $0.54$, $0.37$ for $N/\lambda^2=16$, $20$, $24$, the exponent decreasing as the truncation effect diminishes. Within the computed range the ratio thus grows more slowly than $\log\lambda$, which is consistent with the factor $\log\lambda$ in \eqref{eq:LA} being an upper bound rather than the actual behavior.
\end{enumerate}
\end{fact}
\begin{table}[ht]\centering\small
\begin{tabular}{c|c|ccc}
$\lambda^2$ & $N$ & $\dfrac{\epsilon^N_\lambda(1)}{\delta_4(\lambda)}$ & $\dfrac{\epsilon^N_\lambda(1)}{\delta_4(\lambda)\log\lambda}$ & $\dfrac{\epsilon^N_\lambda(1)}{\delta_4(\lambda)(\log\lambda)^3}$\\[6pt] \hline
5 & 960 & 7.06 & 8.77 & 13.55 \\
9 & 960 & 8.11 & 7.38 & 6.12 \\
13 & 1920 & 9.49 & 7.40 & 4.50 \\
22 & 960 & 9.17 & 5.94 & 2.49 \\
24 & 960 & 9.31 & 5.86 & 2.32 \\
28 & 960 & 9.76 & 5.86 & 2.11 \\
32 & 960 & 11.16 & 6.44 & 2.14 \\
36 & 960 & 10.46 & 5.84 & 1.82 \\
40 & 960 & 11.12 & 6.03 & 1.77 \\
50 & 1920 & 9.66 & 4.94 & 1.29 \\
$100^\ast$ & 3000 & 12.27 & 5.33 & 1.01 \\
$150^\ast$ & 3000 & 13.83 & 5.52 & 0.88 \\
$200^\ast$ & 3200 & 18.07 & 6.82 & 0.97 \\
$250^\ast$ & 5000 & 16.08 & 5.82 & 0.76 \\
$300^\ast$ & 4800 & 19.02 & 6.67 & 0.82 \\
\end{tabular}
\caption{The lowest eigenvalue against the rung $\delta_4(\lambda)$ (exact, as in \eqref{eq:table2}) at the largest $N$ computed; by the min-max principle each entry is an upper bound for the corresponding ratio of $QW_\lambda$. At $\lambda^2=50$ the value $N=1920$ is used (the ratio is $11.00$ at $N=960$). ${}^\ast$Eigenvalues computed by A.L.G.; the exact rungs $\delta_4$ are ours. At $\lambda^2=100$ and $150$ his values at $N=1000,2000$ and $N=1500$ coincide with ours, and the ratio still decreases with $N$ (from $13.79$ at $N=2000$ and $15.22$ at $N=2500$); at $\lambda^2=200$ the ratio is $20.37$ at $N=2800$ and $18.07$ at $N=3200$, and at $\lambda^2=250$ it is $28.27$, $18.68$, $16.08$ at $N=3000$, $4000$, $5000$.}\label{tab:logbound}
\end{table}

\subsection{The odd sector and interlacing}\label{ss:oddsector}
\begin{equation}\label{eq:table3}
\begin{array}{c|cccc|ccc}
\lambda^2 & -\log\delta_6 & -\log\epsilon^{\rm odd}_\lambda(1) & -\log\delta_{10} & -\log\epsilon^{\rm odd}_\lambda(2) & \log\frac{\epsilon^{\rm odd}_\lambda(1)}{\epsilon_\lambda(1)} & \log\frac{\delta_6}{\delta_4} & 4\log\lambda+4.4342\\ \hline
5 & 33.833 & 32.256 & 21.261 & 20.175 & 6.975 & 7.353 & 7.653 \\
9 & 80.008 & 78.307 & 64.578 & 63.392 & 8.280 & 8.672 & 8.829 \\
22 & 237.379 & 235.508 & 218.063 & 216.573 & 10.208 & 10.554 & 10.616 \\
24 & 261.936 & 260.010 & 242.255 & 240.787 & 10.428 & 10.733 & 10.790 \\
28 & 311.183 & 309.225 & 290.860 & 289.284 & 10.729 & 11.050 & 11.099 \\
32 & 360.569 & 358.493 & 339.693 & 338.116 & 10.987 & 11.323 & 11.366 \\
36 & 410.060 & 408.009 & 388.697 & 387.040 & 11.267 & 11.563 & 11.601 \\
40 & 459.633 & 457.571 & 437.837 & 436.230 & 11.431 & 11.778 & 11.812 \\
50 & 583.833 & 581.736 & 561.124 & 559.476 & 11.930 & 12.231 & 12.258 \\
100 & 1207.620 & 1205.361 & 1182.096 & 1180.315 & 13.266 & 13.631 & 13.645 \\
150 & 1833.295 & 1830.872 & 1806.135 & 1804.191 & 14.146 & 14.446 & 14.455 \\
200^\ast & 2459.739 & 2457.160 & 2431.422 & 2429.312 & 14.708 & 15.023 & 15.031 \\
250^\ast & 3086.604 & 3083.988 & 3057.391 & 3054.833 & 15.160 & 15.471 & 15.477 \\
300^\ast & 3713.736 & 3711.079 & 3683.790 & \text{--} & 15.548 & 15.837 & 15.842 \\
\end{array}
\end{equation}
The difference $\log(\epsilon^{\rm odd}_\lambda(1)/\epsilon_\lambda(1))$ falls short of the exact $\log(\delta_6/\delta_4)$ by $0.29$--$0.39$, again without trend ($0.36$ at $\lambda^2=100$, $0.30$ at $\lambda^2=150$, $0.32$, $0.31$, $0.29$ at $\lambda^2=200$, $250$, $300$). The interlacing
$$\epsilon_\lambda(1)<\epsilon^{\rm odd}_\lambda(1)<\epsilon_\lambda(2)<\epsilon^{\rm odd}_\lambda(2)<\cdots$$
holds for every minuscule eigenvalue (eigenvalues below $10^{-2}$) and every $\lambda$ in the table, i.e. for $3$ pairs at $\lambda^2=5$ up to $47$ pairs at $\lambda^2=50$ $96$ pairs at $\lambda^2=100$ and $146$ pairs at $\lambda^2=150$; for $\lambda^2=200$, $250$ it is checked for the first two pairs, the only ones computed.

\subsection{Counting}\label{ss:counting}
\begin{figure}[ht]\centering
\includegraphics[width=0.8\textwidth]{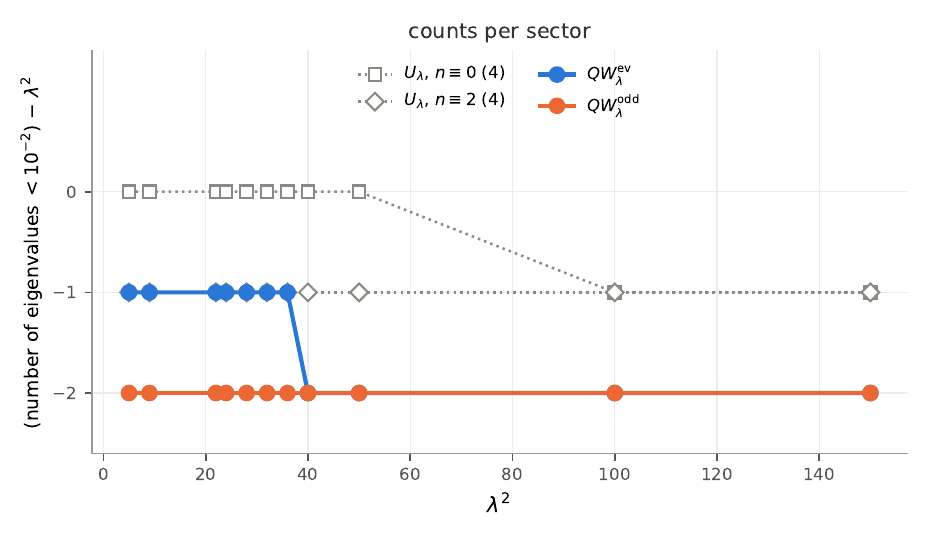}
\caption{Number of eigenvalues below $10^{-2}$ minus $\lambda^2$, for $U_\lambda$ (rungs $\delta_n$, $n\equiv0$ and $n\equiv2$ mod $4$) and for $QW^{\rm ev}_\lambda$, $QW^{\rm odd}_\lambda$.}\label{fig:counting}
\end{figure}

\begin{figure}[ht]\centering
\includegraphics[width=0.8\textwidth]{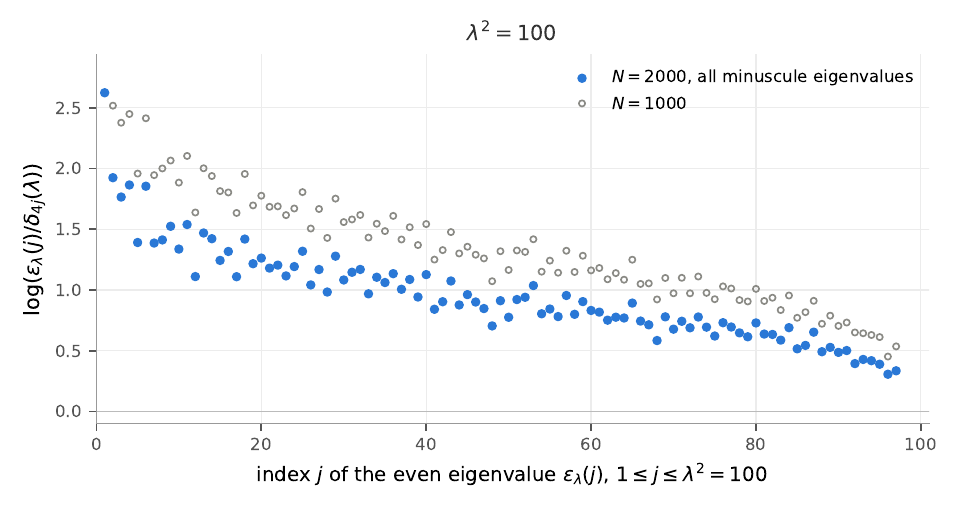}
\caption{$\lambda^2=100$: $\log\big(\epsilon_\lambda(j)/\delta_{4j}(\lambda)\big)$ against the index $j$ of the even eigenvalue, $1\le j\le\lambda^2$, where $\epsilon_\lambda(j)$ is the $j$-th eigenvalue of $QW^{\rm ev}_\lambda$ and $\delta_{4j}(\lambda)=1-\vert\chi_{4j}(\lambda)\vert$ the exact rung of $U_\lambda$ (not Fuchs' asymptotic value, which is only valid for fixed $j$). Filled dots: all $97$ minuscule even eigenvalues at $N=2000=20\lambda^2$; open dots: the same at $N=1000=10\lambda^2$. The truncation lowers the curve by about $0.5$ between the two, a little less near the plunge.}\label{fig:ratio100}
\end{figure}

\begin{figure}[ht]\centering
\includegraphics[width=0.8\textwidth]{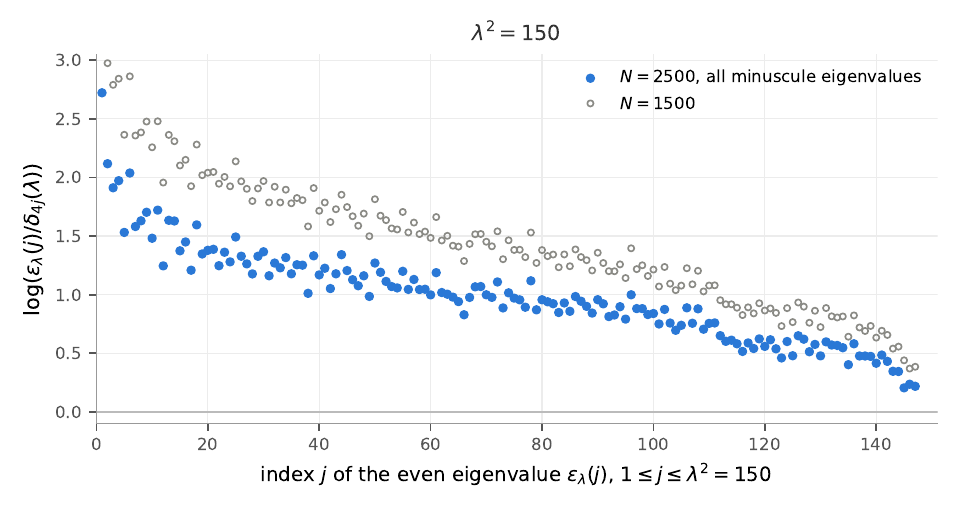}
\caption{The same as Figure \ref{fig:ratio100} for $\lambda^2=150$, $1\le j\le150$. Filled dots: all $147$ minuscule even eigenvalues at $N=2500\approx17\lambda^2$, computed by the batched shifted inverse iteration of \S\ref{ss:method}; open dots: $N=1500=10\lambda^2$. Between $N=1750$ and $N=2500$ the curve still moves down by $0.36$ at $j=1$ and $0.13$ at $j=140$.}\label{fig:ratio150}
\end{figure}

By Theorem \ref{thm:unc} the eigenvalues of $U_\lambda$ below $1$ are exactly the $\delta_n(\lambda)$, $n$ even, and we computed all of them from the prolate functions. Table \ref{tab:count} compares, per sector, the number of eigenvalues below a threshold.
\begin{table}[ht]\centering\small
\begin{tabular}{c|cc|cc|cc}
 & \multicolumn{2}{c|}{$<10^{-10}$} & \multicolumn{2}{c|}{$<10^{-5}$} & \multicolumn{2}{c}{$<10^{-2}$}\\
$\lambda^2$ & $QW_\lambda$ ev/odd & $U_\lambda$ $0/2\ (4)$ & $QW_\lambda$ ev/odd & $U_\lambda$ $0/2\ (4)$ & $QW_\lambda$ ev/odd & $U_\lambda$ $0/2\ (4)$\\ \hline
5 & 2/1 & 3/2 & 3/2 & 4/3 & 4/3 & 5/4 \\
9 & 5/5 & 6/6 & 7/6 & 8/7 & 8/7 & 9/8 \\
22 & 18/17 & 19/18 & 19/19 & 21/20 & 21/20 & 22/21 \\
24 & 20/19 & 21/20 & 21/21 & 23/22 & 23/22 & 24/23 \\
28 & 24/23 & 25/24 & 25/25 & 26/26 & 27/26 & 28/27 \\
32 & 28/27 & 29/28 & 29/29 & 30/30 & 31/30 & 32/31 \\
36 & 32/31 & 33/32 & 33/33 & 34/34 & 35/34 & 36/35 \\
40 & 36/35 & 37/36 & 37/37 & 38/38 & 38/38 & 40/39 \\
50 & 45/45 & 46/46 & 47/47 & 48/48 & 48/48 & 50/49 \\
100 & 95/95 & 96/96 & 97/97 & 98/98 & 98/98 & 99/99 \\
150 & 145/144 & 146/145 & 147/146 & 148/147 & 148/148 & 149/149 \\
\end{tabular}
\caption{Number of eigenvalues below a threshold: $QW^{\rm ev}_\lambda/QW^{\rm odd}_\lambda$ against the rungs $\delta_n$ with $n\equiv0/n\equiv2$ (mod $4$).}\label{tab:count}
\end{table}
\begin{fact}\label{fact:count}
For every $\lambda$ in the table:
\begin{enumerate}
\item $U_\lambda$ has exactly $\lambda^2$ eigenvalues below $10^{-2}$ in the sector $n\equiv0$ (mod $4$) and $\lambda^2-1$ in the sector $n\equiv2$ (mod $4$), in agreement with Slepian's count $2\gamma/\pi=4\lambda^2$ split into four classes; the exceptions are $\lambda^2=100$ and $150$, where the rung $\delta_{4\lambda^2-4}$ ($0.0116$ at $\lambda^2=100$) lies just above the threshold (there are $\lambda^2$ eigenvalues below $0.1$ in each class), which shows that the count at a fixed threshold is only meaningful up to the one or two rungs of the plunge;
\item below $10^{-10}$, each sector of $QW_\lambda$ has exactly one eigenvalue fewer than the corresponding sector of $U_\lambda$; the missing eigenvalue is the rung $\delta_0$ (resp. $\delta_2$) removed by the zero-integral constraint (\S\ref{ss:hermite});
\item at the thresholds $10^{-5}$, $10^{-2}$ the difference is occasionally $2$ (e.g. the even sector at $\lambda^2=40,50$), because each Weil eigenvalue lies above its rung and $\epsilon_\lambda(\lambda^2-1)\approx1$--$2\cdot10^{-2}$ sits on the threshold;
\item pairing $\epsilon_\lambda(j)$ with $\delta_{4j}$ and $\epsilon^{\rm odd}_\lambda(j)$ with $\delta_{4j+2}$, the ratios satisfy, for all minuscule $j$,
\begin{equation}\label{eq:ratiorange}
1.40\le\frac{\epsilon_\lambda(j)}{\delta_{4j}(\lambda)}\le11.2,\qquad 1.51\le\frac{\epsilon^{\rm odd}_\lambda(j)}{\delta_{4j+2}(\lambda)}\le8.2 \qquad(\lambda^2\le50),
\end{equation}
the largest ratio occurring at $j=1$ and the ratio decreasing steadily to about $1.5$ at the plunge, as essentially the same function of $j/\lambda^2$ for all $\lambda$ (Figures \ref{fig:ratio100} and \ref{fig:ratio150} for $\lambda^2=100$, $150$; the counts are in Figure \ref{fig:counting}). At $\lambda^2=100$ the lower ends are $1.36$ (even sector, $N=2000$; $1.57$ at $N=1000$) and $1.58$, while the ratios at $j=1$ are $13.8$ and $9.6$ at $N=2000$; the extrapolated convergence in $N$ (Table \ref{tab:conv}) would bring them to about $10$ and $7$, within \eqref{eq:ratiorange}. At $\lambda^2=150$ the even ratios at $N=2500$ lie between $1.23$ ($j=145$) and $15.2$ ($j=1$), and in the odd sector the ratio at $j=1$ is $11.3$ ($N=2500$) and the lower end $1.48$ ($N=1500$); at $N=1500$ the ratios for $j\lesssim10$ are still inflated by the truncation (up to $37.6$ at $j=1$). The lower ends decrease slowly with $\lambda$ and with $N$ ($1.40$, $1.36$, $1.23$ at $\lambda^2=50,100,150$), so that the lower bound in \eqref{eq:ratiorange} is specific to $\lambda^2\le50$.
\end{enumerate}
\end{fact}
Thus each sector of $QW_\lambda$ contains about $\lambda^2-1$ minuscule eigenvalues, one less than the corresponding sector of $U_\lambda$, and \eqref{eq:ratiorange} is the logarithmic agreement of Hypothesis \ref{hyp:LA} with explicit constants in the whole range $1\le j\le\lambda^2-1$ (subject to the convergence in $N$ of Table \ref{tab:conv}, which can only lower the ratios).
\begin{table}[ht]\centering\small
\begin{tabular}{c|ccc|ccc}
$\lambda^2$ & \# even & $j=1$ & range & \# odd & $j=1$ & range\\ \hline
5 & 4 & 7.06 & 1.77--7.06 & 3 & 4.84 & 2.92--4.84 \\
9 & 8 & 8.11 & 1.71--8.11 & 7 & 5.48 & 1.94--5.48 \\
22 & 19 & 9.17 & 1.65--9.17 & 19 & 6.49 & 1.82--6.49 \\
24 & 21 & 9.31 & 1.77--9.31 & 21 & 6.86 & 1.67--6.86 \\
28 & 25 & 9.76 & 1.64--9.76 & 25 & 7.09 & 1.63--7.09 \\
32 & 29 & 11.16 & 1.56--11.16 & 29 & 7.97 & 1.63--7.97 \\
36 & 33 & 10.46 & 1.61--10.46 & 33 & 7.77 & 1.68--7.77 \\
40 & 37 & 11.12 & 1.71--11.12 & 37 & 7.86 & 1.51--7.86 \\
50 & 47 & 11.00 & 1.40--11.00 & 47 & 8.14 & 1.51--8.14 \\
100 & 97 & 13.79 & 1.36--13.79 & 96 & 9.58 & 1.58--9.58 \\
150 & 147 & 15.22 & 1.23--15.22 & 146 & 11.27 & 1.48--17.64 \\
\end{tabular}
\caption{The ratios $\epsilon_\lambda(j)/\delta_{4j}(\lambda)$ and $\epsilon^{\rm odd}_\lambda(j)/\delta_{4j+2}(\lambda)$ over the minuscule eigenvalues. For $\lambda^2=100$ the even ratios are taken at $N=2000$ and the odd ones at $N=1000$ ($j=1$ at $N=2000$); for $\lambda^2=150$ the even ratios at $N=2500$, and the odd ones at $N=1500$ ($j=1,2$ at $N=2500$), where the ratios for small $j$ are still inflated by the truncation (the upper end $17.6$ occurs at $j\approx6$).}\label{tab:ratio}
\end{table}

\section{Eliminating the minuscule eigenvalues by prolate conditioning}\label{sec:conditioning}

The dictionary predicts not only the minuscule eigenvalues but their eigenvectors: the images under $\cE$ of the zero-integral prolate combinations. We construct these vectors explicitly, compare them with the eigenvectors of $QW_\lambda$, and show that the orthogonal complement of their span, together with the vanishing of the pole term, carries no minuscule eigenvalue at all.

\subsection{The prolate vectors}\label{ss:prolatevectors}
For $k\ge1$ let
\begin{equation}\label{eq:fk}
f^{\rm ev}_k=I_0\,h_{4k,\lambda}-I_{4k}\,h_{0,\lambda},\qquad
f^{\rm odd}_k=I_2\,h_{4k+2,\lambda}-I_{4k+2}\,h_{2,\lambda},
\end{equation}
where $I_n:=\int_\R h_{n,\lambda}(x)\,dx$, so that $\int f=0$ and $f^{\rm ev}_1=h_\lambda$ is the two-mode state \eqref{eq:twomode}. The function $\cE(f)(u)=u^{1/2}\sum_{1\le n<\lambda^2}f(nu)$ is restricted to the window and projected on $E_N$, the components $\langle\cE(f),b_j\rangle$ being computed with the inner product of $L^2([\lambda^{-1},\lambda],\dstar)$; only the integers $n<\lambda^2$ contribute, and the integral is split at the points $u=\lambda/n$ where a term leaves the support (Appendix \ref{app:matrix}). By \S\ref{sec:dict}, $\cE(f^{\rm ev}_k)$ is approximately $\iota$-even and $\cE(f^{\rm odd}_k)$ approximately $\iota$-odd; numerically the wrong-parity part has relative size $10^{-14}$ for small $k$ and grows only near the plunge. Projecting on the respective sector and applying Gram--Schmidt in the order $k=1,2,\dots$ gives orthonormal vectors
$$\mathfrak p^{\rm ev}_k\in E_N^{\rm ev},\qquad \mathfrak p^{\rm odd}_k\in E_N^{\rm odd}$$
(Gram--Schmidt in $E_N$ in the interleaved order $f^{\rm ev}_1,f^{\rm odd}_1,f^{\rm ev}_2,\dots$ gives the same vectors up to $10^{-15}$, except near the plunge). These are the natural candidates for the eigenvectors $\eta_\lambda(k)$, $\eta^{\rm odd}_\lambda(k)$ of $\epsilon_\lambda(k)$, $\epsilon^{\rm odd}_\lambda(k)$, as in \cite{CC}, \S3.

\begin{figure}[ht]\centering
\includegraphics[width=\textwidth]{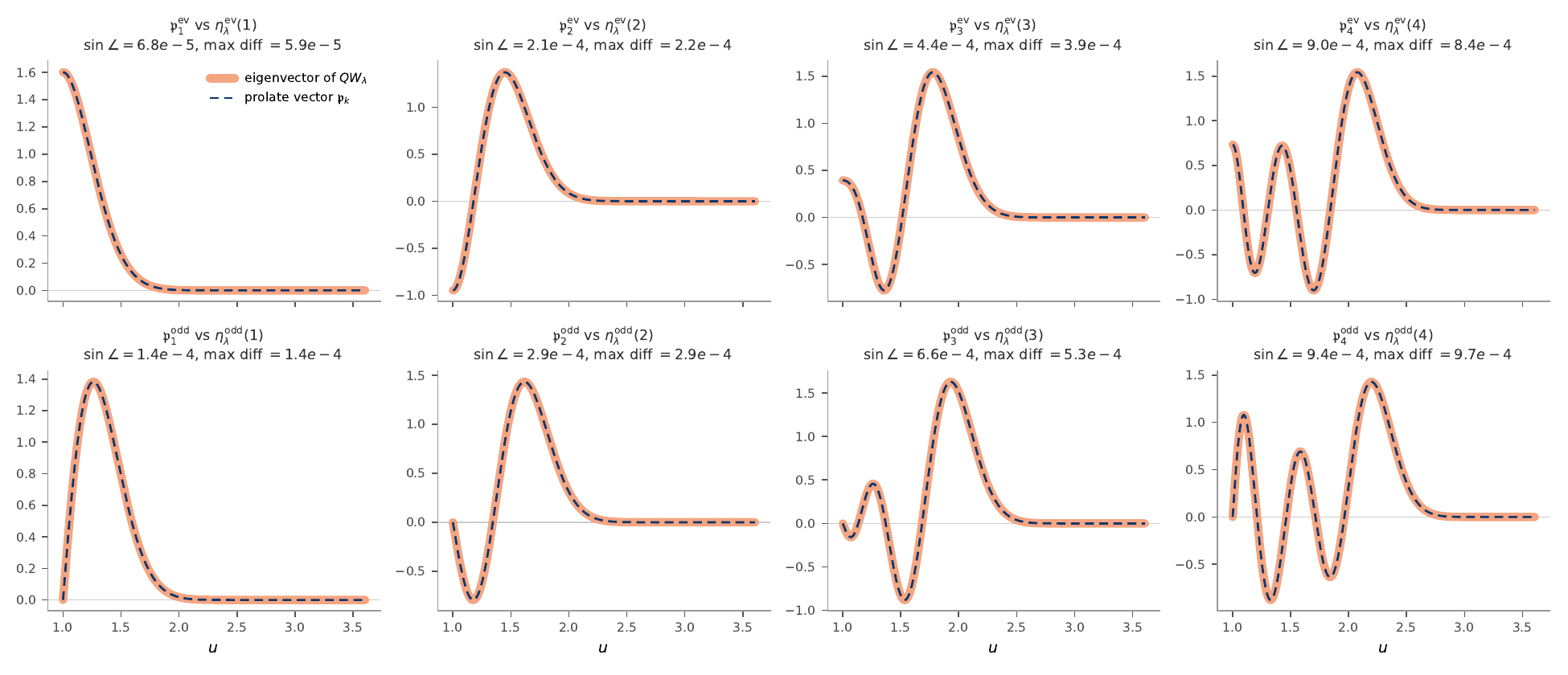}
\caption{$\lambda^2=13$, $N=120$: the prolate vectors $\mathfrak p^{\rm ev}_k$, $\mathfrak p^{\rm odd}_k$ (dark dashed line) and the eigenvectors $\eta_\lambda(k)$, $\eta^{\rm odd}_\lambda(k)$ of $QW_\lambda$ (wide light band), $k=1,\dots,4$, on $u\in[1,\lambda]$, both normalized in $L^2(\dstar)$. Each pair appears as a single graph because the two functions coincide to within the width of the line: the angle between them is $\sin\angle\le10^{-3}$, and the maximal pointwise difference, relative to $\max\vert\eta_\lambda(k)\vert$ and given above each panel, is between $6\cdot10^{-5}$ and $10^{-3}$. The functions on $[\lambda^{-1},1]$ follow by $\eta\circ\iota=\pm\eta$.}\label{fig:gs}
\end{figure}

\begin{fact}\label{fact:gs}
$\mathfrak p^{\rm ev}_k$ is close to $\pm\eta_\lambda(k)$ and $\mathfrak p^{\rm odd}_k$ to $\pm\eta^{\rm odd}_\lambda(k)$ for every minuscule eigenvalue, with $\sin\angle\le2.6\cdot10^{-2}$ for $\lambda^2\le50$ ($4.2\cdot10^{-2}$ at $\lambda^2=100$, $6.3\cdot10^{-2}$ at $\lambda^2=150$), the largest angles occurring at the plunge, and much smaller for small $k$:
\begin{table}[ht]\centering\small
\begin{tabular}{c|ccc|ccc|c}
$\lambda^2$ & $\mathfrak p^{\rm ev}_1$ & $\mathfrak p^{\rm ev}_2$ & $\mathfrak p^{\rm ev}_3$ & $\mathfrak p^{\rm odd}_1$ & $\mathfrak p^{\rm odd}_2$ & $\mathfrak p^{\rm odd}_3$ & max over minuscule\\ \hline
5 & 3.4e-04 & 9.5e-04 & 2.3e-03 & 7.1e-04 & 1.6e-03 & 8.3e-03 & 1.7e-02 \\
9 & 1.1e-04 & 3.0e-04 & 3.4e-04 & 2.5e-04 & 3.4e-04 & 8.1e-04 & 2.6e-02 \\
22 & 2.0e-05 & 3.9e-05 & 6.8e-05 & 2.6e-05 & 4.8e-05 & 1.1e-04 & 1.4e-02 \\
24 & 1.6e-05 & 3.6e-05 & 8.0e-05 & 2.5e-05 & 4.8e-05 & 1.2e-04 & 1.6e-02 \\
28 & 1.2e-05 & 2.3e-05 & 5.0e-05 & 1.8e-05 & 3.1e-05 & 6.8e-05 & 1.8e-02 \\
32 & 8.1e-06 & 2.0e-05 & 4.4e-05 & 1.3e-05 & 2.8e-05 & 6.6e-05 & 1.7e-02 \\
36 & 7.8e-06 & 1.6e-05 & 3.2e-05 & 1.0e-05 & 2.1e-05 & 5.3e-05 & 1.8e-02 \\
40 & 6.8e-06 & 1.5e-05 & 2.9e-05 & 1.1e-05 & 2.1e-05 & 3.9e-05 & 2.2e-02 \\
50 & 5.4e-06 & 1.2e-05 & 2.2e-05 & 8.2e-06 & 1.6e-05 & 3.0e-05 & 2.6e-02 \\
100 & 1.2e-06 & 2.8e-06 & 4.8e-06 & 2.0e-06 & 3.6e-06 & 6.9e-06 & 4.2e-02 \\
150 & 1.9e-07 & 4.7e-07 & 8.7e-07 & 3.3e-07 & 6.4e-07 & 1.3e-06 & 6.3e-02 \\
\end{tabular}
\caption{$\sin\angle(\mathfrak p_k,\eta_\lambda(k))$ at $N=960$ ($N=1000$ for $\lambda^2=100$, $N=1500$ for $\lambda^2=150$), first three vectors of each sector and maximum over all eigenvalues below $10^{-2}$.}\label{tab:gs}
\end{table}
\end{fact}
The angle of the first vector decreases with $\lambda$, from $3.4\cdot10^{-4}$ at $\lambda^2=5$ to $5.4\cdot10^{-6}$ at $\lambda^2=50$ $1.2\cdot10^{-6}$ at $\lambda^2=100$ ($1.3\cdot10^{-6}$ at $N=2000$) and $1.9\cdot10^{-7}$ at $\lambda^2=150$ ($6.4\cdot10^{-7}$ at $N=2500$). Part of it is truncation: at $\lambda^2=13$,
\begin{table}[ht]\centering\small
\begin{tabular}{c|cccc|cccc}
$N$ & $\mathfrak p^{\rm ev}_1$ & $\mathfrak p^{\rm ev}_2$ & $\mathfrak p^{\rm ev}_3$ & $\mathfrak p^{\rm ev}_4$ & $\mathfrak p^{\rm odd}_1$ & $\mathfrak p^{\rm odd}_2$ & $\mathfrak p^{\rm odd}_3$ & $\mathfrak p^{\rm odd}_4$\\ \hline
120 & 6.8e-05 & 2.1e-04 & 4.4e-04 & 9.0e-04 & 1.4e-04 & 2.9e-04 & 6.6e-04 & 9.4e-04 \\
240 & 4.9e-05 & 1.6e-04 & 3.3e-04 & 6.7e-04 & 9.5e-05 & 2.2e-04 & 5.3e-04 & 7.3e-04 \\
480 & 3.5e-05 & 1.3e-04 & 2.8e-04 & 6.1e-04 & 7.6e-05 & 1.8e-04 & 4.0e-04 & 5.6e-04 \\
960 & 2.9e-05 & 1.2e-04 & 2.4e-04 & 5.1e-04 & 6.8e-05 & 1.7e-04 & 3.9e-04 & 5.4e-04 \\
1920 & 2.3e-05 & 1.0e-04 & 2.2e-04 & 4.8e-04 & 6.7e-05 & 1.6e-04 & 3.8e-04 & 5.2e-04 \\
\end{tabular}
\caption{$\sin\angle(\mathfrak p_k,\eta_\lambda(k))$ at $\lambda^2=13$ as a function of $N$.}\label{tab:gsN}
\end{table}
All angles decrease with $N$; in the odd sector they level off near $N\approx1000$ at a value which is the genuine distance, at this $\lambda$, between the prolate vector and the eigenvector, while $\mathfrak p^{\rm ev}_1$ is still improving by a factor $0.8$ per doubling. The angles are the same in double precision and with $60$ digits. Figure \ref{fig:gs} shows the agreement in the style of \cite{CC}, Figures 26--37.

\subsection{The pole term}\label{ss:pole}
Let $\ell(g)=\int_{\lambda^{-1}}^\lambda g(u)\,u^{1/2}\dstar=\cM g(\tfrac12)$, so that $\ell(g\circ\iota)=\cM g(-\tfrac12)$.

\begin{lem}\label{lem:pole}
For real $g$ supported in the window, $W_{0,2}(F)=2\,\cM g(\tfrac12)\,\cM g(-\tfrac12)$ with $F$ as in \eqref{eq:Psi}. Hence the pole part of $QW_\lambda$ is the rank-one form $+2\,\ell(g)^2$ on $\iota$-even functions and $-2\,\ell(g)^2$ on $\iota$-odd functions.
\end{lem}
\begin{proof}
$W_{0,2}(F)=\int_0^\infty\!\!\int_0^\infty g(u)g(ux)(x^{1/2}+x^{-1/2})\dstar\,d^*x$, and the substitution $v=ux$ gives $\int g(u)u^{-1/2}\dstar\int g(v)v^{1/2}d^*v+\int g(u)u^{1/2}\dstar\int g(v)v^{-1/2}d^*v=2\cM g(\tfrac12)\cM g(-\tfrac12)$. For $g\circ\iota=\pm g$ one has $\cM g(-\tfrac12)=\pm\cM g(\tfrac12)$.
\end{proof}
In the basis $b_j$ (Appendix \ref{app:matrix}), with $\omega_j=2\pi j/L$, $L=2\log\lambda$,
\begin{equation}\label{eq:ellbasis}
\ell(b_0)=\frac{4\sinh(L/4)}{\sqrt L},\qquad \ell(b_j)=\sqrt{\frac2L}\,\frac{\sinh(L/4)}{\tfrac14+\omega_j^2},\qquad \ell(b_{-j})=\sqrt{\frac2L}\,\frac{2\omega_j\sinh(L/4)}{\tfrac14+\omega_j^2}\qquad(j\ge1),
\end{equation}
and the restriction of the pole part of the matrix to $E_N^{\rm ev}$, resp. $E_N^{\rm odd}$, is $+2\,\ell\ell^T$, resp. $-2\,\ell\ell^T$ (checked against \cite{CCM}, (4.2), to $10^{-52}$). Since $\iota$ exchanges $\cM g(\tfrac12)$ and $\cM g(-\tfrac12)$, the single condition $\ell(g)=0$ kills both pole values in each sector.

\subsection{Prolate conditioning}\label{ss:prolatecond}
In each sector we restrict $QW_\lambda$ to the orthogonal complement of the prolate vectors and of the pole form:
\begin{equation}\label{eq:mu}
\mu^{\rm ev}_\lambda=\min\Big\{\frac{QW_\lambda(g,g)}{\Vert g\Vert^2}:\ g\in E_N^{\rm ev},\ g\perp\mathfrak p^{\rm ev}_1,\dots,\mathfrak p^{\rm ev}_{\lambda^2-1},\ \ell(g)=0\Big\},
\end{equation}
and $\mu^{\rm odd}_\lambda$ likewise with $\mathfrak p^{\rm odd}_1,\dots,\mathfrak p^{\rm odd}_{\lambda^2-2}$. The numbers of prolate conditions, $\lambda^2-1$ and $\lambda^2-2$, are the numbers of minuscule eigenvalues found in \S\ref{ss:counting}; the pole condition replaces the rung removed by the zero-integral constraint. By Courant--Fischer, with $k$ conditions $\mu$ is at most the $(k+1)$-st eigenvalue of the sector, i.e. $\mu^{\rm ev}_\lambda\le\epsilon_\lambda(\lambda^2+1)$ and $\mu^{\rm odd}_\lambda\le\epsilon^{\rm odd}_\lambda(\lambda^2)$.
\begin{table}[ht]\centering\small
$$\begin{array}{c|cccc|cccc}
 & \multicolumn{4}{c|}{\text{even sector}} & \multicolumn{4}{c}{\text{odd sector}}\\
\lambda^2 & \text{prolates} & \mu^{\rm ev}_\lambda & \text{exact}+\ell & \epsilon_\lambda(\lambda^2+1) & \text{prolates} & \mu^{\rm odd}_\lambda & \text{exact}+\ell & \epsilon^{\rm odd}_\lambda(\lambda^2)\\ \hline
5 & 0.606 & \mathbf{0.782} & 0.784 & 0.853 & 0.080 & \mathbf{0.731} & 0.731 & 0.767 \\
9 & 0.520 & \mathbf{0.658} & 0.658 & 0.663 & 0.087 & \mathbf{0.633} & 0.634 & 0.651 \\
13 & 0.428 & \mathbf{0.492} & 0.487 & 0.582 & 0.086 & \mathbf{0.471} & 0.471 & 0.477 \\
22 & 0.375 & \mathbf{0.437} & 0.442 & 0.485 & 0.092 & \mathbf{0.438} & 0.440 & 0.450 \\
24 & 0.340 & \mathbf{0.399} & 0.401 & 0.435 & 0.096 & \mathbf{0.383} & 0.384 & 0.396 \\
28 & 0.319 & \mathbf{0.379} & 0.381 & 0.444 & 0.098 & \mathbf{0.361} & 0.362 & 0.384 \\
32 & 0.321 & \mathbf{0.404} & 0.404 & 0.407 & 0.092 & \mathbf{0.395} & 0.396 & 0.403 \\
36 & 0.320 & \mathbf{0.388} & 0.388 & 0.390 & 0.099 & \mathbf{0.377} & 0.381 & 0.389 \\
40 & 0.317 & \mathbf{0.335} & 0.335 & 0.366 & 0.102 & \mathbf{0.332} & 0.332 & 0.334 \\
50 & 0.262 & \mathbf{0.270} & 0.271 & 0.310 & 0.104 & \mathbf{0.270} & 0.270 & 0.273 \\
100 & 0.194 & \mathbf{0.195} & 0.195 & 0.234 & 0.093 & \mathbf{0.194} & 0.195 & 0.195 \\
150 & 0.184 & \mathbf{0.199} & 0.199 & 0.199 & 0.093 & \mathbf{0.194} & 0.194 & 0.198 \\
\end{array}$$
\caption{Lowest eigenvalue of $QW_\lambda$ ($N=960$; $N=2000$ for $\lambda^2=100$, $N=2500$ for $\lambda^2=150$) on the orthogonal complement of: the prolate vectors only; the prolate vectors and $\ell$ ($\mu_\lambda$); the exact minuscule eigenvectors and $\ell$; last column, the Courant--Fischer bound.}\label{tab:cond}
\end{table}

\begin{fact}\label{fact:cond}
For every $\lambda$ with $5\le\lambda^2\le50$, and for $\lambda^2=100$ and $150$:
\begin{enumerate}
\item No minuscule eigenvalue survives the conditioning: $0.27\le\mu^{\rm ev}_\lambda,\mu^{\rm odd}_\lambda\le0.78$ for $\lambda^2\le50$, $\mu^{\rm ev}_\lambda=\mu^{\rm odd}_\lambda=0.1945$ at $\lambda^2=100$, and $0.199/0.194$ at $\lambda^2=150$. In fact $1.6\le\lambda\,\mu_\lambda\le2.5$ throughout, with $\lambda\mu_\lambda=1.91$ at $\lambda^2=50$, $1.945$ at $\lambda^2=100$ and $2.43/2.37$ at $\lambda^2=150$: the conditioned minimum is of order $2/\lambda$, with fluctuations of $\pm20\%$ that follow the edge of the non-minuscule spectrum.
\item In the even sector the $\lambda^2-1$ prolate vectors alone already give the first non-minuscule eigenvalue $\epsilon_\lambda(\lambda^2)$ (to $4$ digits): they span the minuscule eigenspace. The pole condition raises the minimum further, by $29\%$ at $\lambda^2=5$, $3\%$ at $\lambda^2=50$, $0.3\%$ at $\lambda^2=100$ and $8\%$ at $\lambda^2=150$.
\item In the odd sector the $\lambda^2-2$ prolate vectors leave one transitional eigenvalue $\epsilon^{\rm odd}_\lambda(\lambda^2-1)\approx0.08$--$0.10$, the partner of the plunge rung $\delta_{4\lambda^2-2}$; the pole condition removes it, multiplying the minimum by a factor between $2.1$ ($\lambda^2=100,150$) and $9$ ($\lambda^2=5$) and bringing it within $0.2$--$6\%$ of the Courant--Fischer bound.
\item Replacing the prolate vectors by the exact minuscule eigenvectors changes $\mu_\lambda$ by at most $1.3\%$.
\item For $\lambda^2\ge22$ the two sectors give nearly the same value ($0.437/0.438$ at $\lambda^2=22$, $0.270/0.270$ at $\lambda^2=50$, $0.1945/0.1945$ at $\lambda^2=100$, $0.199/0.194$ at $\lambda^2=150$), which decreases like $2/\lambda$ together with the edge of the non-minuscule spectrum.
\end{enumerate}
The values change by at most $0.5\%$ between $N=480$ and $N=960$, and by $1\%$ between $N=120$ and $N=960$ at $\lambda^2=13$; at $\lambda^2=100$ they change by $0.7\%$ between $N=500$ and $N=1000$ and by less than $0.05\%$ between $N=1000$ and $N=2000$; at $\lambda^2=150$ by $3\%$ (even) and $1\%$ (odd) between $N=1500$ and $N=2500$.
\end{fact}
In other words, once the $\cE$-images of the prolate functions with $\vert\chi_n\vert\approx1$ and the pole term are factored out, $QW_\lambda$ is bounded below on what remains by about $2/\lambda$, the edge of the non-minuscule spectrum ($0.27$ at $\lambda^2=50$, $0.19$ at $\lambda^2=100$, $0.19$--$0.20$ at $\lambda^2=150$): the minuscule part of the Weil spectrum is entirely accounted for by the prolate side of the dictionary, with one prolate condition per rung of $U_\lambda$ except the rung removed by $\int f=0$, which is replaced by the pole condition.
\begin{figure}[ht]\centering
\includegraphics[width=\textwidth]{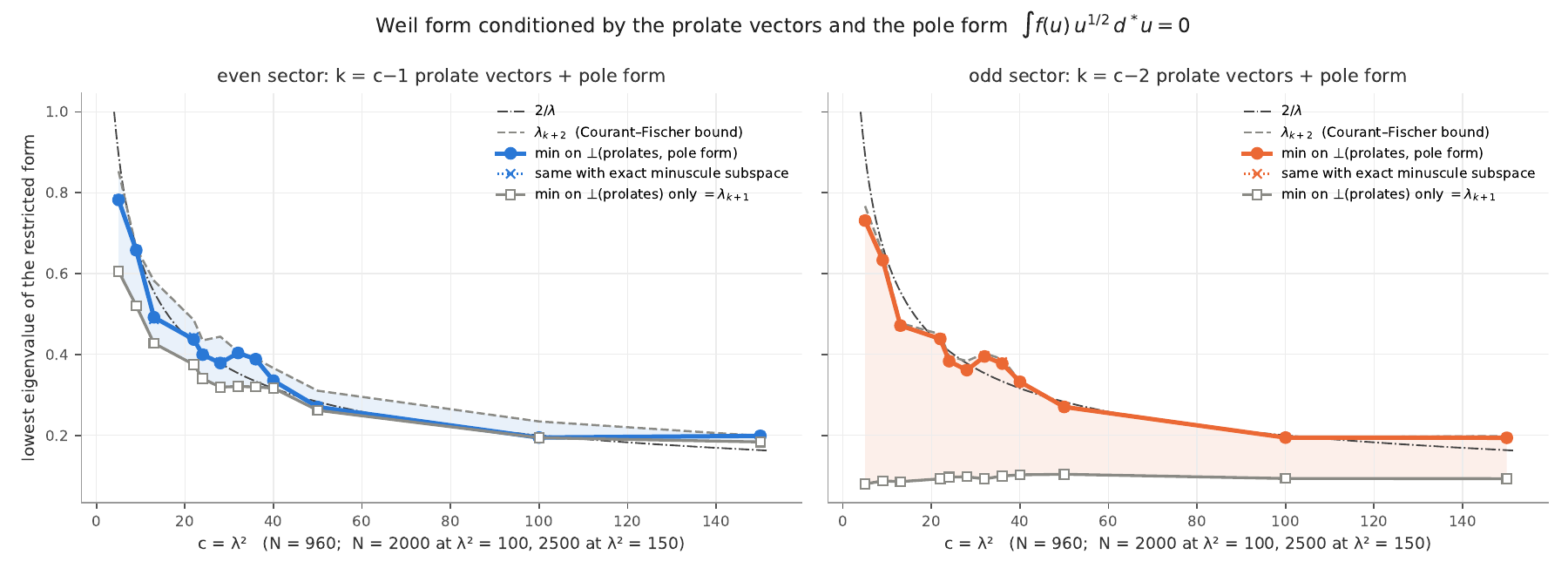}
\caption{$\mu^{\rm ev}_\lambda$, $\mu^{\rm odd}_\lambda$ (thick), the same with the exact minuscule eigenspace (dotted), the minimum with the prolate vectors only (squares), the Courant--Fischer bound (dashed) and the curve $2/\lambda$ (dash-dotted), against $\lambda^2$; $N=960$, $N=2000$ at $\lambda^2=100$ and $N=2500$ at $\lambda^2=150$.}\label{fig:conditioning}
\end{figure}

\section{The ground state}\label{sec:ground}

In this section $QW_\lambda\ge0$ is assumed (Weil positivity on the window), so that $\epsilon_\lambda\ge0$, together with the gap part of Hypothesis \ref{hyp:LA}. We write $\delta=\delta_4(\lambda)$, $\langle\cdot,\cdot\rangle$ and $\Vert\cdot\Vert$ for $L^2([\lambda^{-1},\lambda],\dstar)$, and we normalize $\Vert k_\lambda\Vert=1$.

\subsection{Simplicity, evenness, and the projection on the ground state}

Let $\eta_\lambda$ be a ground state, $\xi_\lambda$ the orthogonal projection of $k_\lambda$ on the ground eigenspace and $\omega_\lambda=k_\lambda-\xi_\lambda$ its projection on the spectral subspace of the remaining eigenvalues, so that $\Vert\xi_\lambda\Vert^2+\Vert\omega_\lambda\Vert^2=1$ and $QW_\lambda(\xi_\lambda,\omega_\lambda)=0$.

\begin{lem}\label{lem:gap}
Assume $QW_\lambda\ge0$ and that, the eigenvalues being counted with multiplicity, every eigenvalue of $QW_\lambda$ other than the lowest one $\epsilon_\lambda$ is $\ge g_\lambda>0$ (so that $\epsilon_\lambda$ is simple as soon as $g_\lambda>\epsilon_\lambda$). Then
\begin{equation}\label{eq:omega}
\Vert\omega_\lambda\Vert^2\le\frac{C\,\delta_4(\lambda)(\log\lambda)^3}{g_\lambda},\qquad
\Big\Vert k_\lambda-\frac{\eta_\lambda}{\Vert\eta_\lambda\Vert}\Big\Vert\le\Vert\omega_\lambda\Vert(1+\Vert\omega_\lambda\Vert^2)^{1/2},
\end{equation}
the sign of $\eta_\lambda$ being chosen so that $\langle k_\lambda,\eta_\lambda\rangle>0$. If moreover $\Vert\omega_\lambda\Vert<1/2$, the ground eigenvalue is simple and the ground state is $\iota$-even.
\end{lem}
\begin{proof}
By orthogonality $QW_\lambda(k_\lambda,k_\lambda)=\epsilon_\lambda\Vert\xi_\lambda\Vert^2+QW_\lambda(\omega_\lambda,\omega_\lambda)\ge g_\lambda\Vert\omega_\lambda\Vert^2$, using $\epsilon_\lambda\ge0$; then \eqref{eq:upper}. Next $1-\Vert\xi_\lambda\Vert=\Vert\omega_\lambda\Vert^2/(1+\Vert\xi_\lambda\Vert)\le\Vert\omega_\lambda\Vert^2$ and $\Vert k_\lambda-\eta_\lambda/\Vert\eta_\lambda\Vert\Vert^2=\Vert\omega_\lambda\Vert^2+(1-\Vert\xi_\lambda\Vert)^2$ when the ground eigenspace is a line. Simplicity: the gap hypothesis as stated makes the ground eigenspace one-dimensional. Evenness: $QW_\lambda$ commutes with $\iota$, so the ground line is $\iota$-invariant and $\iota\eta_\lambda=\pm\eta_\lambda$; by the Poisson defect \eqref{eq:poissondefect} and the estimates of \cite{CK}, $\Vert k_\lambda-k_\lambda\circ\iota\Vert\le C\lambda^{1/2}\delta^{1/2}\log\lambda$, so the $\iota$-odd part of $k_\lambda$ has norm $o(1)$, while $\langle k_\lambda,\eta_\lambda\rangle/\Vert\eta_\lambda\Vert=\Vert\xi_\lambda\Vert>\sqrt{3}/2$; the sign $-$ is impossible.
\end{proof}

The positivity $\epsilon_\lambda\ge0$ is used: without it one only gets $g_\lambda\Vert\omega_\lambda\Vert^2\le C\delta(\log\lambda)^3+\vert\epsilon_\lambda\vert$, and the construction gives no lower bound on $\epsilon_\lambda$. Two versions of the gap are relevant, with the ratios \eqref{eq:gaps}:
\begin{itemize}
\item \emph{full gap} $g_\lambda=c\,\delta_6(\lambda)$, i.e. the lower half of \eqref{eq:LA} for $\epsilon^{\rm odd}_\lambda(1)$ and $\epsilon_\lambda(2)$: $\delta_4/\delta_6=\frac{15}{128\pi^2\lambda^4}$ and $\Vert\omega_\lambda\Vert\le C\lambda^{-2}(\log\lambda)^{3/2}$;
\item \emph{even gap} $g_\lambda=c\,\delta_8(\lambda)$ for the second even eigenvalue only: since the $\iota$-odd component of $k_\lambda$ is $O(\lambda^{1/2}\delta^{1/2}\log\lambda)$ and $QW_\lambda(k_\lambda,k_\lambda)$ dominates each of its two sector contributions (both $\ge0$), the argument can be run inside the even sector, $\delta_4/\delta_8=\frac{105}{4096\pi^4\lambda^8}$ and $\Vert\omega_\lambda\Vert\le C\lambda^{-4}(\log\lambda)^{3/2}$.
\end{itemize}
The nested version applies verbatim to the spans $H_j(\lambda)$: under \eqref{eq:LA} the $j$-th eigenvector lies within $O(\lambda^{-4}(\log\lambda)^{3/2})$ of ${\rm span}\{E_\lambda\xi(\lambda,k)\}_{k\le j}$, and the Gram--Schmidt orthogonalization of the $E_\lambda\xi(\lambda,k)$ in the order $k=1,2,\ldots$ approximates the successive eigenvectors with that accuracy (Davis--Kahan). This is confirmed numerically for every minuscule eigenvalue, not only the first few: see Fact \ref{fact:gs} and Table \ref{tab:gs} in Section \ref{sec:conditioning}.

\subsection{From $L^2$ to the norm $\Vert\cdot\Vert_1$}\label{ss:L1}

Let $\rho(u)=u^{1/2}$ for $1\le u\le\lambda$, $\rho(u)=\rho(1/u)=u^{-1/2}$ for $\lambda^{-1}\le u\le1$. Then $\Vert\rho\Vert^2=2(\lambda-1)$ and $\rho\,\dstar=u^{-1/2}du$ on $[1,\lambda]$, so for $g$ on the window
\begin{equation}\label{eq:schwarzrho}
\Vert g\Vert_1:=\int_1^\lambda\vert g(u)\vert u^{-1/2}du\le\langle\vert g\vert,\rho\rangle\le\sqrt{2(\lambda-1)}\,\Vert g\Vert,\qquad \Vert g\Vert_1\le\sqrt{\lambda-1}\,\Vert g\mathbf 1_{[1,\lambda]}\Vert .
\end{equation}
The pairing with $\rho$ is the pairing with the pole vector $u^{\pm1/2}$ of the explicit formula, which is why this is the natural $L^1$ norm; what it controls is explained in \S\ref{ss:strip}. Applied to Lemma \ref{lem:gap},
\begin{equation}\label{eq:L1unnorm}
\Big\Vert k_\lambda-\frac{\eta_\lambda}{\Vert\eta_\lambda\Vert}\Big\Vert_1\le\sqrt{\lambda-1}\,\Vert\omega_\lambda\Vert(1+o(1))\le
\begin{cases}C\lambda^{-3/2}(\log\lambda)^{3/2} & \text{(full gap)},\\[2pt] C\lambda^{-7/2}(\log\lambda)^{3/2} & \text{(even gap)}.\end{cases}
\end{equation}
Cauchy--Schwarz discards the structure of $\omega_\lambda$; since $\rho$ is far from an eigenvector nothing better is available without information on the spectral measure of $\rho$, but the factor $\sqrt\lambda$ is paid once.

\subsection{Normalization}\label{ss:normalization}

The estimates above concern $L^2$-normalized vectors, whereas in the computations the ground state and the comparison functions are normalized at $u=1$. Point evaluation is not bounded on $L^2(\dstar)$, so the $L^2$ bound on $\omega_\lambda$ says nothing about $\omega_\lambda(1)$, and the estimates transfer to a normalization by a bounded linear functional. The target of Theorem \ref{thm:mellintheta} provides the natural one: with $\tau=\theta/\theta(1)$, define the positive averaging functional
\begin{equation}\label{eq:P}
P(g):=\frac{\langle g,\theta\rangle}{\langle\tau,\theta\rangle}=\frac{\theta(1)}{\Vert\theta\Vert^2}\,\langle g,\theta\rangle ,\qquad P(\tau)=1 ,
\end{equation}
the mean of $g/\tau$ for the weight $\theta^2\dstar$, a window of logarithmic width $\approx0.2$ around $u=1$ ($99.2\%$ of the mass of $\theta^2\dstar$ lies in $[\frac23,\frac32]$, all but $10^{-6}$ in $[\frac12,2]$); normalizing by $P(g)=1$ is the same as choosing the multiple of $g$ closest in $L^2$ to $\tau$. Put $\tilde k_\lambda=k_\lambda/P(k_\lambda)$, $\tilde\eta_\lambda=\eta_\lambda/P(\eta_\lambda)$. Since $k_\lambda\to\theta/\Vert\theta\Vert$ in $L^2$ of the window (Theorem \ref{thm:upper} and the Hermite limit, at the rate $O(\lambda^{-2})$), $a:=P(k_\lambda)=\theta(1)/\Vert\theta\Vert+O(\lambda^{-2})=1.580+O(\lambda^{-2})$, $b:=P(\eta_\lambda/\Vert\eta_\lambda\Vert)$ satisfies $\vert a-b\vert\le\frac{\theta(1)}{\Vert\theta\Vert}\Vert k_\lambda-\eta_\lambda/\Vert\eta_\lambda\Vert\Vert$, and
\begin{equation}\label{eq:normalized}
\begin{aligned}\Vert\tilde k_\lambda-\tilde\eta_\lambda\Vert&\le\frac{\Vert\omega_\lambda\Vert}{a}\Big(1+\frac{\theta(1)}{\Vert\theta\Vert\,b}\Big)(1+o(1))=(1.27+o(1))\Vert\omega_\lambda\Vert ,\\
\Vert\tilde k_\lambda-\tilde\eta_\lambda\Vert_1&\le(1.27+o(1))\sqrt{\lambda-1}\,\Vert\omega_\lambda\Vert .\end{aligned}
\end{equation}
The exponents of \eqref{eq:L1unnorm} are thus preserved for the normalization $P(\cdot)=1$, and for any fixed bounded functional non-degenerate on $\theta$.

For the strict normalization $\eta_\lambda(1)=1$ one needs in addition $\vert\omega_\lambda(1)\vert$, i.e. a modulus of continuity of the ground state at $u=1$. Near $u=1$ the test vector is as regular as $h$ ($u\partial_u\cE(f_\lambda)=\frac12\cE(f_\lambda)+u^{1/2}\sum_nnu\,f_\lambda'(nu)$, a rapidly convergent sum with $O(1)$ terms), so it is Lipschitz in $t=\log u$ on $[\frac12,2]$ with constant $O(1)$. If the ground state is Lipschitz there with constant $L_\lambda$, averaging against a bump of logarithmic width $r$ gives $\vert\omega_\lambda(1)\vert\le r^{-1/2}\Vert\omega_\lambda\Vert+(C+L_\lambda)r$, hence
\begin{equation}\label{eq:pointnorm}
\vert\omega_\lambda(1)\vert\le C\Vert\omega_\lambda\Vert^{2/3}(1+L_\lambda)^{1/3},\qquad
\Big\Vert\frac{k_\lambda}{k_\lambda(1)}-\frac{\eta_\lambda}{\eta_\lambda(1)}\Big\Vert_1\le C\lambda^{-13/6}(\log\lambda)(1+L_\lambda)^{1/3}
\end{equation}
under the even gap ($k_\lambda(1)\to\theta(1)/\Vert\theta\Vert\ne0$). Since the eigenvalue equation only yields logarithmic Sobolev regularity, a bound on $L_\lambda$ is at present a numerical observation. In the computations we report both normalizations and the diagnostic ratio $P(\eta_\lambda)/\eta_\lambda(1)$, which equals $1$ for $\tau$ itself.

\subsection{The three distances}

Let $\theta_\lambda=\cE(h_\lambda)$ with $h_\lambda$ as in \eqref{eq:twomode}, $\tau_\lambda=\theta_\lambda/\theta_\lambda(1)$ (or its $P$-normalized version $\tilde\tau_\lambda$), and $\tau=\theta/\theta(1)$ the target. The comparison of the ground state with the target splits into three pieces of very different sizes:
\begin{enumerate}
\item \emph{ground state versus test vector}: $\Vert\tilde\eta_\lambda-\tilde k_\lambda\Vert_1\le C\lambda^{-7/2}(\log\lambda)^{3/2}$ (even gap; $\lambda^{-3/2}$ for the full gap), conditional on $\epsilon_\lambda\ge0$ and on the gap;
\item \emph{test vector versus $\theta_\lambda$}: $\Vert\tilde k_\lambda-\tilde\tau_\lambda\Vert_1\le C\lambda^{-1/2}\delta^{1/2}\log\lambda$ by Theorem \ref{thm:upper}: the Gevrey repair and the lower cutoff are invisible at the scale of the target;
\item \emph{$\theta_\lambda$ versus $\theta$} (prolate versus Hermite): $h_\lambda-h$ is $\lambda^{-2}$ times a Hermite-type function (e.g. $h_{0,\lambda}=h_0-\frac{\sqrt6}{32\pi}\lambda^{-2}h_4+O(\lambda^{-4})$ uniformly on $[-\lambda,\lambda]$), so $\vert\theta_\lambda(u)-\theta(u)\vert\le C\lambda^{-2}u^{-1/2}$ for $u\ge1$ up to the tail of $\theta$, and directly in the $L^1$ norm, without Schwarz loss,
$$\Vert\tilde\tau_\lambda-\tau\Vert_1\le C\lambda^{-2}\log\lambda,\qquad\Vert\tau_\lambda-\tau\Vert_1\le C\lambda^{-2}\log\lambda .$$
\end{enumerate}
Consequently the distance to target $\Vert\eta_\lambda-\tau\Vert_1$ is dominated by the third piece: it should decay like $\lambda^{-2}\log\lambda$, a property of the Hermite limit of the prolate functions and not of the Weil form. The information carried by Theorem \ref{thm:upper} and the gap is only visible if the ground state is compared with $\tau_\lambda$ (or with $\cE$ of the second-order Meixner--Sch\"afke expansion of $h_\lambda$) rather than with $\tau$, and with the normalization $P(\cdot)=1$. \smallskip\noindent\emph{Numerical values.} Table \ref{tab:dist} gives the quantities of Appendix \ref{app:report} for $5\le\lambda^2\le250$, with $\eta_\lambda$ the ground state of $QW^N_\lambda$ normalized by $\eta_\lambda(1)=1$, $g_\lambda=\eta_\lambda-\tau_\lambda$, and tildes for the normalization $P(\cdot)=1$; $\theta_\lambda$ is evaluated by the direct sum for $u\ge1$ and by the dual expression \eqref{eq:dual} for $u<1$, and the integrals by the trapezoidal rule in $t=\log u$ with $2\cdot10^4$ and $4\cdot10^4$ intervals on $[0,\log\lambda]$, which agree to all digits shown. At $\lambda^2=13$, $N=120$ the values coincide with an independent computation by R.A. to ten digits.
\begin{table}[ht]\centering\small\setlength{\tabcolsep}{3pt}
$$\begin{array}{cc|cc|cc|ccc}
\lambda^2 & N & \Vert g_\lambda\Vert_1 & \Vert\tilde\eta_\lambda-\tilde\tau_\lambda\Vert_1 & \Vert\tau_\lambda-\tau\Vert_1 & \Vert\tilde\tau_\lambda-\tau\Vert_1 & P(\eta_\lambda) & 1-\rho_\lambda^2 & \max\vert\cM g_\lambda\vert\\ \hline
5 & 960 & 1.35\cdot10^{-4} & 1.28\cdot10^{-4} & 2.68\cdot10^{-2} & 2.46\cdot10^{-2} & 0.96129 & 1.18\cdot10^{-7} & 2.30\cdot10^{-4}\\
9 & 960 & 4.63\cdot10^{-5} & 4.26\cdot10^{-5} & 1.44\cdot10^{-2} & 1.30\cdot10^{-2} & 0.97976 & 1.22\cdot10^{-8} & 7.89\cdot10^{-5}\\
13 & 1920 & 1.00\cdot10^{-5} & 9.09\cdot10^{-6} & 9.86\cdot10^{-3} & 8.85\cdot10^{-3} & 0.98632 & 5.45\cdot10^{-10} & 1.70\cdot10^{-5}\\
22 & 960 & 8.92\cdot10^{-6} & 8.02\cdot10^{-6} & 5.76\cdot10^{-3} & 5.15\cdot10^{-3} & 0.99207 & 4.17\cdot10^{-10} & 1.51\cdot10^{-5}\\
24 & 960 & 7.04\cdot10^{-6} & 6.33\cdot10^{-6} & 5.27\cdot10^{-3} & 4.71\cdot10^{-3} & 0.99275 & 2.59\cdot10^{-10} & 1.20\cdot10^{-5}\\
28 & 960 & 5.36\cdot10^{-6} & 4.81\cdot10^{-6} & 4.51\cdot10^{-3} & 4.02\cdot10^{-3} & 0.99381 & 1.49\cdot10^{-10} & 9.10\cdot10^{-6}\\
32 & 960 & 3.58\cdot10^{-6} & 3.21\cdot10^{-6} & 3.94\cdot10^{-3} & 3.51\cdot10^{-3} & 0.99460 & 6.62\cdot10^{-11} & 6.08\cdot10^{-6}\\
36 & 960 & 3.46\cdot10^{-6} & 3.10\cdot10^{-6} & 3.50\cdot10^{-3} & 3.12\cdot10^{-3} & 0.99521 & 6.16\cdot10^{-11} & 5.88\cdot10^{-6}\\
40 & 960 & 3.00\cdot10^{-6} & 2.68\cdot10^{-6} & 3.14\cdot10^{-3} & 2.80\cdot10^{-3} & 0.99569 & 4.61\cdot10^{-11} & 5.09\cdot10^{-6}\\
50 & 1920 & 2.01\cdot10^{-6} & 1.79\cdot10^{-6} & 2.51\cdot10^{-3} & 2.24\cdot10^{-3} & 0.99657 & 2.05\cdot10^{-11} & 3.40\cdot10^{-6}\\
100 & 2000 & 5.78\cdot10^{-7} & 5.14\cdot10^{-7} & 1.25\cdot10^{-3} & 1.11\cdot10^{-3} & 0.99829 & 1.68\cdot10^{-12} & 9.79\cdot10^{-7}\\
150 & 2500 & 2.87\cdot10^{-7} & 2.56\cdot10^{-7} & 8.33\cdot10^{-4} & 7.41\cdot10^{-4} & 0.99887 & 4.14\cdot10^{-13} & 4.87\cdot10^{-7}\\
250^\ast & 5000 & 9.76\cdot10^{-8} & 8.68\cdot10^{-8} & 5.00\cdot10^{-4} & 4.44\cdot10^{-4} & 0.99932 & 4.80\cdot10^{-14} & 1.65\cdot10^{-7}\\
\end{array}$$
\caption{The three distances in the norm $\Vert\cdot\Vert_1$, with $\eta_\lambda(1)=\tau_\lambda(1)=1$ (tildes: $P(\cdot)=1$); $P(\eta_\lambda)$ is the ratio $P(\eta_\lambda)/\eta_\lambda(1)$, $\rho_\lambda$ the overlap of $\theta_\lambda$ and $\eta_\lambda$, and $\max\vert\cM g_\lambda\vert$ the maximum on $\Re s=\pm\frac12$, $\vert\Im s\vert\le100$. ${}^\ast$Ground state computed by A.L.G.; the distances are computed with the same program as the other rows.}\label{tab:dist}
\end{table}
\begin{enumerate}
\item \emph{The third piece is exactly of order $\lambda^{-2}$.} $\lambda^2\Vert\tau_\lambda-\tau\Vert_1$ decreases from $0.134$ ($\lambda^2=5$) to $0.1251$ ($\lambda^2=100$), $0.1250$ ($\lambda^2=150$) and $0.1249$ ($\lambda^2=250$), and is fitted by $0.1247+0.044\lambda^{-2}$ for $\lambda^2\ge22$; likewise $\lambda^2\Vert\tilde\tau_\lambda-\tau\Vert_1\to0.111$. The logarithm in the bound $C\lambda^{-2}\log\lambda$ is therefore not present. This piece dominates the distance to the target: $\Vert\eta_\lambda-\tau\Vert_1$ agrees with $\Vert\tau_\lambda-\tau\Vert_1$ to three digits for every $\lambda$.
\item \emph{Ground state versus $\tau_\lambda$.} $\Vert g_\lambda\Vert_1$ is two to three orders of magnitude smaller, from $1.3\cdot10^{-4}$ to $9.8\cdot10^{-8}$, and a least-squares fit over the thirteen values gives $\Vert g_\lambda\Vert_1\approx C\lambda^{-3.5}$ (exponent $-3.54$, and $-3.50$ without $\lambda^2=250$), i.e. the exponent of the even-gap estimate \eqref{eq:L1unnorm}, without visible logarithm. These values still depend on the truncation: at $\lambda^2=13$, $\Vert g_\lambda\Vert_1=2.94,\,2.12,\,1.50,\,1.23,\,1.00\cdot10^{-5}$ for $N=120,\dots,1920$; at $\lambda^2=50$, $2.64,\,2.41,\,2.01\cdot10^{-6}$ for $N=480,960,1920$; at $\lambda^2=100$, $5.38$ and $5.78\cdot10^{-7}$ for $N=1000,2000$; at $\lambda^2=150$, $0.87$ and $2.87\cdot10^{-7}$ for $N=1500,2500$, the value at $N=2500$ being the one on the $\lambda^{-3.5}$ fit; at $\lambda^2=250$, $0.78$, $1.02$ and $0.98\cdot10^{-7}$ for $N=3000$, $4000$, $5000$ (A.L.G.), the last two within $4\%$ of each other and $15\%$ below the fit over $\lambda^2\le150$. The order of magnitude and the exponent are stable, the individual values only to $10$--$20\%$; the distances computed at $N\le5\lambda^2$ (e.g. $3.3\cdot10^{-4}$ at $\lambda^2=100$, $N=500$) measure the truncation, not $g_\lambda$.
\item \emph{Shape of the deviation.} $g_\lambda$ is proportional to the next target $\theta_1(b,\cdot)$ of \S\ref{ss:nexttarget} up to a relative $L^2(\dstar)$ residual $0.53\lambda^{-2}$ ($12\%$ at $\lambda^2=5$, $1.1\%$ at $\lambda^2=50$, $0.53\%$ at $\lambda^2=100$, $0.35\%$ at $\lambda^2=150$), and this residual is independent of $N$ to three digits even where the amplitude is not: the ground state deviates from $\tau_\lambda$ in the direction of the second eigenvector, as the Davis--Kahan argument predicts.
\item \emph{Normalization.} $1-P(\eta_\lambda)/\eta_\lambda(1)\approx0.17\lambda^{-2}$ ($0.19\lambda^{-2}$ at $\lambda^2=5$), and $P(\tau_\lambda)/\tau_\lambda(1)$ is the same to $10^{-6}$: the difference between the two normalizations is a property of $\tau_\lambda$, not of the ground state. Accordingly $\Vert\tilde\eta_\lambda-\tilde\tau_\lambda\Vert_1/\Vert g_\lambda\Vert_1$ lies between $0.89$ and $0.95$. The ratios $1/\Vert\eta_\lambda\Vert$ and $\theta_\lambda(1)/\Vert\theta_\lambda\Vert$ agree to $3\cdot10^{-4}$ at $\lambda^2=5$, $2\cdot10^{-5}$ at $\lambda^2=22$ and $10^{-6}$ at $\lambda^2=100$, and both equal $\theta(1)/\Vert\theta\Vert+0.27\lambda^{-2}$, with $\theta(1)/\Vert\theta\Vert=1.5799$.
\item \emph{Mellin transform and pointwise size.} $\max_{\Re s=\pm1/2}\vert\cM g_\lambda(s)\vert=1.70\Vert g_\lambda\Vert_1$ for every $\lambda$, below the bound $2\Vert g_\lambda\Vert_1$ of \eqref{eq:stripbound}; $\sup\vert g_\lambda\vert=(2.3$--$2.6)\Vert g_\lambda\Vert_1$, attained near $u=1.41$; the contribution of $[\lambda/2,\lambda]$ to $\Vert g_\lambda\Vert_1$ is $5\%$ at $\lambda^2=13$ and below $10^{-9}$ of it for $\lambda^2\ge22$. The symmetry defect $\sup\vert\tau_\lambda-\tau_\lambda\circ\iota\vert$ is at the level of the double-precision rounding ($\le6\cdot10^{-15}$), and $\eta_\lambda\circ\iota=\eta_\lambda$ exactly, the even sector being used.
\item \emph{Overlaps.} $1-\rho_\lambda^2$ is the square of the angle between $\eta_\lambda$ and the prolate vector $\mathfrak p^{\rm ev}_1$ of Table \ref{tab:gs}, from $1.2\cdot10^{-7}$ to $4\cdot10^{-13}$ ($\lambda^2=150$) and $5\cdot10^{-14}$ ($\lambda^2=250$).
\end{enumerate}

\subsection{The next target}\label{ss:nexttarget}
For the second even eigenvector the Hermite limit lies in the zero-integral part of ${\rm span}\{h_0,h_4,h_8\}$, hence $\cE$ of a function of the form $u^2e^{-\pi u^2}\times$(polynomial of degree $6$ in $u^2$), and by Proposition \ref{prop:xipoly} its Mellin transform is $\xi(s+\frac12)$ times a polynomial of degree one in $s^2$. A convenient representative, normalized so as to vanish at the normalization point $u=1$ of the ground state, is $\theta_1(b,u)=\sqrt u\sum_{n\ge1}d(b,nu)$ with
\begin{align*}d(b,u)=e^{-\pi u^2}u^2\Big(&\frac14\sqrt{\tfrac32}\big(32\pi^3u^6-224\pi^2u^4+350\pi u^2-105\big)\,b\\&+16\pi^3u^6-112\pi^2u^4+210\pi u^2-105\Big)\end{align*}
and $b$ determined by $\theta_1(b,1)=0$, $b=-0.75243767525394408195519\ldots$ Up to scale, $d(b,\cdot)$ is the only linear combination of the Hermite functions $h_{4j}$, $j=0,1,2$, which vanishes at $0$ and is such that $\theta_1(b,1)=0$: the two functions multiplying $b$ and $1$ in $d(b,u)$ are invariant under the Fourier transform and vanish at $0$ (hence have vanishing integral), and they span the combinations of $h_0,h_4,h_8$ vanishing at $0$.

\section{Mellin transforms and the zeros of $\zeta$}\label{sec:zeros}

\subsection{What the norm $\Vert\cdot\Vert_1$ controls}\label{ss:strip}

For $g$ supported in $[\lambda^{-1},\lambda]$, $\cM g(s)=\int g(u)u^s\dstar$ is entire, and for $\sigma=\Re s\in[-\frac12,\frac12]$
\begin{equation}\label{eq:stripbound}
\vert\cM g(s)\vert\le\int_1^\lambda\vert g\vert u^{1/2}\dstar+\int_{\lambda^{-1}}^1\vert g\vert u^{-1/2}\dstar=\Vert g\Vert_1+\Vert g\circ\iota\Vert_1 ,
\end{equation}
the second integral becoming $\Vert g\circ\iota\Vert_1$ under $u\mapsto1/u$. Thus the $L^1$ distance of \S\ref{ss:L1} controls the Mellin transforms uniformly on the closed strip $\vert\Re s\vert\le\frac12$, which is the critical strip of $\zeta(s+\frac12)$ by \eqref{eq:mellinE}. By Phragm\'en--Lindel\"of the supremum on the strip is attained on the two lines $\Re s=\pm\frac12$, i.e. it is the larger of the two suprema of the Fourier transforms in $t$ of $g(e^t)e^{\pm t/2}$, computable by two FFTs; computing it next to the right-hand side of \eqref{eq:stripbound} shows the slack of the $L^1$ bound. For $\iota$-even $g$ the two terms coincide.

\subsection{Zeros of the Mellin transform of the ground state}

\begin{thm}[{\cite[Theorem 5.6]{CvS}}]\label{thm:realzeros}
Let $N\ge1$ and assume that the lowest eigenvalue $\epsilon^N_\lambda(1)$ of $QW^N_\lambda$ is simple, with eigenvector $\eta=\sum_{\vert j\vert\le N}\xi_jV_j$ even ($\xi_{-j}=\xi_j$). Then all the zeros of the entire function $t\mapsto\cM\eta(it)$ are real: they are the points $t=2\pi j/L$, $\vert j\vert>N$, and the $2N$ points $t=2\pi s/L$, $s$ a root of the polynomial
\begin{equation}\label{eq:Pcvs}
P(s)=\sum_{\vert k\vert\le N}\xi_k\prod_{\vert j\vert\le N,\,j\neq k}(j-s),
\end{equation}
all of whose roots are real.
\end{thm}
Indeed $Q=QW^N_\lambda-\epsilon^N_\lambda(1)$ is a real symmetric matrix of the form considered in \cite{CvS}: by Appendix \ref{app:matrix} (see also \cite[Proposition 4.3]{CCM}) its off-diagonal entries are $(b_n-b_m)/(n-m)$ with $b_{-n}=-b_n$, and its diagonal is even in $n$; it is positive, with kernel $\C\eta$, even under $V_j\mapsto V_{-j}$. Theorem 5.6 of \cite{CvS} then applies, after rescaling the interval $[0,1]$ of \cite{CvS} to the window $[-\frac L2,\frac L2]$ in $x=\log u$, where $\cM\eta(it)=\int\eta(e^x)e^{itx}dx$. The hypothesis is checked at finite $N$: the lowest even eigenvalue is separated from the next even one and from the lowest odd one by factors comparable to $\delta_8/\delta_4$ and $\delta_6/\delta_4$ (Tables \eqref{eq:table2}, \eqref{eq:table3}). If moreover the lowest eigenvalue of $QW_\lambda$ is simple and isolated, the ground states of $QW^N_\lambda$ converge in $L^2$ to that of $QW_\lambda$ as $N\to\infty$, their Mellin transforms converge uniformly on compact sets, and by Hurwitz' theorem the zeros of $\cM\eta_\lambda$ are real as well.

\begin{rem}[counting the zeros]\label{rem:countzeros}
Since $P$ is even, of degree $2N$, with $P(0)\neq0$, its roots come in pairs $\pm s$ and the function $t\mapsto\cM\eta(it)$ has \emph{at most} $N$ zeros in $t>0$, counted with multiplicity, apart from the points $2\pi j/L$, $j>N$, with equality if and only if all the roots of $P$ are real. This bound does not use the positivity: it only uses the degree of $P$, so it also holds for the finite function defined by approximate (exported) coefficients. Consequently, if the sign of $\cM\eta(it)$ is certified at points $0<t_0<t_1<\cdots$ (beyond $2\pi N/L$, the sign of $\cM\eta(it)/\sin(tL/2)$, which does not vanish at the points $2\pi j/L$) and $N$ sign changes are found, each interval of sign change contains exactly one zero, which is simple, the count of positive real zeros is exactly $N$, and there is no other zero, real or complex: the list of zeros is complete, without using Theorem \ref{thm:realzeros} and without any bound on derivatives. The leading coefficient of $P$ is $\sum_k\xi_k=\langle\eta\vert\sum_kV_k\rangle$, which is very small for the ground state (of the order of $\epsilon^N_\lambda(1)^{1/2}$ in our computations), and so are the next moments; correspondingly, part of the $N$ zeros lie far beyond $2\pi N/L$ (\S\ref{ss:zeroscomp}).
\end{rem}

By Theorem \ref{thm:mellintheta}, $\cM\tau(s)=\xi(s+\frac12)/(4\theta(1))$, whose zeros in the variable $s$ are the $\rho-\frac12$, $\rho$ the non-trivial zeros of $\zeta$. Hurwitz' theorem on the zeros of uniform limits of holomorphic functions gives:

\begin{prop}\label{prop:hurwitz}
Let $\lambda_j\to\infty$ and let $\tilde\eta_{\lambda_j}$ be the ground states normalized by $P(\tilde\eta_{\lambda_j})=1$. If
$$\Vert\tilde\eta_{\lambda_j}-\tau\Vert_1+\Vert(\tilde\eta_{\lambda_j}-\tau)\circ\iota\Vert_1\longrightarrow0 ,$$
then $\cM\tilde\eta_{\lambda_j}\to\xi(\cdot+\frac12)/(4\theta(1))$ uniformly on the closed strip $\vert\Re s\vert\le\frac12$, every zero of $\xi(s+\frac12)$ in the open strip is a limit of zeros of $\cM\tilde\eta_{\lambda_j}$, and if these are real (Theorem \ref{thm:realzeros}) then all zeros of $\zeta$ in the open critical strip are on the critical line.
\end{prop}
\begin{proof}
The uniform convergence is \eqref{eq:stripbound} applied to $g=\tilde\eta_{\lambda_j}-\tau\mathbf 1_{[\lambda_j^{-1},\lambda_j]}$ (the tail of $\tau$ outside the window is negligible in both norms). Hurwitz' theorem applies on every disc inside the open strip since the limit is not identically zero; a limit of real points is real.
\end{proof}

\begin{rem}[on the logic]
The convergence required in Proposition \ref{prop:hurwitz} is what Section \ref{sec:ground} estimates, but under the assumption $QW_{\lambda}\ge0$ for the $\lambda$ considered, which for all $\lambda$ is equivalent to the conclusion. The conditional results therefore do not give a proof; what they show is that the picture is self-consistent -- positivity and the Fuchs ladder force the ground state onto the target, whose Mellin transform is $\xi$ -- and that the distance to be controlled decays like $\lambda^{-2}\log\lambda$, entirely because of the prolate-to-Hermite convergence. An unconditional treatment would have to replace the positivity by information on the spectral measure of $QW_\lambda$, which is exactly the lower bound direction of Hypothesis \ref{hyp:LA}.
\end{rem}

\subsection{Computations}\label{ss:zeroscomp}

\begin{figure}[ht]\centering
\includegraphics[width=0.95\textwidth]{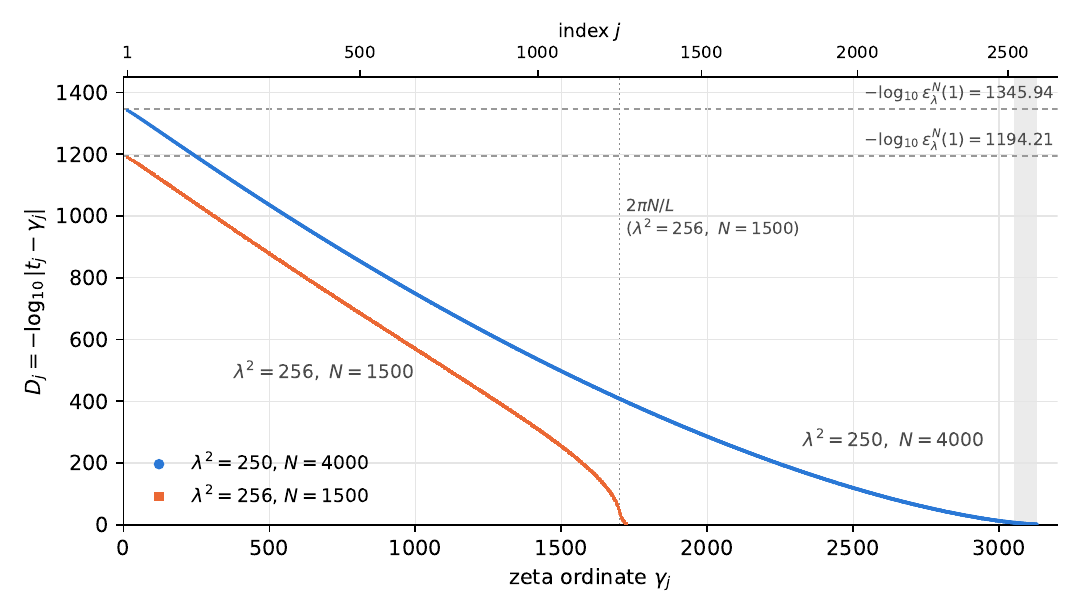}
\caption{The number of digits of agreement $D_j=-\log_{10}\vert t_j-\gamma_j\vert$ between the zeros $t_j$ of $t\mapsto\cM\eta(it)$ and the ordinates $\gamma_j$ of the zeros of $\zeta$, plotted against $\gamma_j$ (lower axis) and $j$ (upper axis), one dot per zero (the dots merge into curves): $\lambda^2=250$, $N=4000$, $j=1,\dots,2596$, and $\lambda^2=256$, $N=1500$, $j=1,\dots,1268$ (ground states computed by A.L.G.). Dashed: the levels $-\log_{10}\epsilon^N_\lambda(1)$. Dotted: $2\pi N/L$ for $\lambda^2=256$, $N=1500$, where that curve ends; shaded: $t\approx3050$--$3130$, where the agreement disappears for $\lambda^2=250$. Each $t_j$ is obtained by Newton's method from $\gamma_j$ with the exported coefficients of the ground state, at a precision exceeding $D_j$ by at least $150$ digits; for $j\le200$ the values agree with the interval enclosures of A.L.G. to $5\cdot10^{-9}$.}\label{fig:mellindigits}
\end{figure}

The phenomenon was discovered in \cite[\S6]{CCM}, where the zeros of $\cM\eta$ appear as the spectrum of a self-adjoint rank-one perturbation of the scaling operator: with $N=120$, the first $50$ zeros of $\zeta$ are reproduced for $\lambda^2=12$, $13$, $14$, the first one to within $3.4\cdot10^{-50}$, $2.4\cdot10^{-55}$ and $1.1\cdot10^{-60}$ respectively, and the fiftieth to within about $10^{-3}$ for $\lambda^2=13$. At $\lambda=16$ ($\lambda^2=256$) it was then observed by R.A., with a much larger accuracy; it has been reproduced by A.L.G. with an independent computation of the ground state at $N=1500$ (even sector of dimension $1501$, $3400$ decimal digits, quadrature order $18816$), for which $\epsilon^N_\lambda(1)=6.142\cdot10^{-1195}$, i.e. $-\log_{10}\epsilon^N_\lambda(1)=1194.21$. The zeros $t_j$ of $t\mapsto\cM\eta(it)$ near the first $200$ ordinates $\gamma_j$ of the zeros of $\zeta$ were isolated as described below, and we have recomputed them independently from the exported coefficients, with identical results. With $D_j=-\log_{10}\vert t_j-\gamma_j\vert$:
\begin{equation}\label{eq:zeros256}
\begin{array}{r|r|c}
j & \gamma_j & D_j\ (\lambda^2=256,\ N=1500)\\ \hline
1 & 14.13472514 & 1190.12 \\
2 & 21.02203964 & 1186.54 \\
10 & 49.77383248 & 1168.90 \\
50 & 143.11184581 & 1108.98 \\
100 & 236.52422967 & 1047.95 \\
150 & 318.85310426 & 994.49 \\
200 & 396.38185422 & 944.15 \\
\end{array}
\end{equation}
All the differences $t_j-\gamma_j$ are positive, and $D_j$ decreases from $1190$ to $944$, comparable to the values (about $1150$ to $950$) of the original computation of R.A., whose truncation is not recorded. Note that $-\log_{10}\epsilon^N_\lambda(1)-D_1=4.090$, and $4.086$ at $\lambda^2=250$ below, so that $\vert t_1-\gamma_1\vert\approx1.2\cdot10^4\,\epsilon^N_\lambda(1)$ in both cases, although $N=1500\approx6\lambda^2$ is far from the converged regime of \S\ref{sec:numspec}: the accuracy of the zeros is governed by the eigenvalue of the truncated form, not by its distance to the limit. Beyond the first $200$ zeros, $D_j$ decreases to $570$ at $t=1000$ and $254$ at $t=1500$, and then drops abruptly to $0$ between $t=1680$ and $t=1730$ (Figure \ref{fig:mellindigits}): at this truncation the agreement ends at the highest log-frequency $2\pi N/L=1699.6$ of the basis, whereas at $\lambda^2=250$ with $N=3000$ and $4000$ ($2\pi N/L=3414$ and $4553$) it ends earlier, near $t\approx3100$, independently of $N$. 

An independent computation by A.L.G. at $\lambda^2=250$ confirms the phenomenon, with the truncations $N=3000$ and $N=4000$ of \S\ref{sec:numspec} (even sectors of dimension $3001$ and $4001$, ground states computed with $2900$ decimal digits). In the variable $s=it$, $t$ real, and for $\eta$ in the even sector with coefficients $v_j$ in the basis $b_j$, the transform is the finite sum
\begin{equation}\label{eq:mellinfinite}
\cM\eta(it)=2\sin\big(\tfrac{tL}{2}\big)\Big(\frac{v_0}{\sqrt L\,t}+\sqrt{\tfrac2L}\;t\sum_{j=1}^N\frac{v_j}{t^2-a_j^2}\Big),\qquad a_j=\frac{2\pi j}{L},\quad L=\log\lambda^2 ,
\end{equation}
whose zeros other than those of $\sin(tL/2)$ were located near the first $200$ ordinates $\gamma_j$ of the zeros of $\zeta$ by Newton's method followed by an interval-Newton inclusion, at $1500$ and again at $1800$ decimal digits, with identical results. All $200$ roots $t_j$ are simple and isolated, and satisfy $t_j>\gamma_j$. With $D_j=-\log_{10}\vert t_j-\gamma_j\vert$:
\begin{equation}\label{eq:zeros250}
\begin{array}{r|r|cc}
j & \gamma_j & D_j\ (N=3000) & D_j\ (N=4000)\\ \hline
1 & 14.13472514 & 1341.68 & 1341.86 \\
2 & 21.02203964 & 1338.10 & 1338.28 \\
10 & 49.77383248 & 1320.51 & 1320.69 \\
50 & 143.11184581 & 1260.99 & 1261.17 \\
100 & 236.52422967 & 1200.76 & 1200.94 \\
150 & 318.85310426 & 1148.34 & 1148.52 \\
200 & 396.38185422 & 1099.30 & 1099.48 \\
\end{array}
\end{equation}
We have checked these values independently, evaluating \eqref{eq:mellinfinite} directly with the exported coefficients of the $N=4000$ ground state: the differences $t_1-\gamma_1=1.38999579112714\cdot10^{-1342}$ and $t_{200}-\gamma_{200}=3.31178230707278\cdot10^{-1100}$ are reproduced to all digits shown, and the reference ordinate $\gamma_1$ agrees with an independent evaluation to $1400$ digits. The accuracy is governed by the lowest eigenvalue: between $N=3000$ and $N=4000$, $-\log_{10}\epsilon_\lambda(1)$ increases by $0.1799$, and every $D_j$, $1\le j\le200$, increases by the same amount, between $0.1799$ and $0.1804$, and again from $N=4000$ to $N=5000$ (A.L.G.): $0.06520$ for $-\log_{10}\epsilon_\lambda(1)$, between $0.06508$ and $0.06520$ for the $D_j$ ($D_1=1341.92$, $D_{200}=1099.55$ at $N=5000$); thus $\vert t_j-\gamma_j\vert$ is proportional to $\epsilon_\lambda(1)$, uniformly in $j$, and $\vert t_1-\gamma_1\vert\approx1.2\cdot10^4\,\epsilon_\lambda(1)$ (at $N=4000$, $-\log_{10}\epsilon_\lambda(1)=1345.94$ and $-\log_{10}\delta_4(\lambda)=1347.21$).

Beyond the first $200$ zeros, we located all the sign changes of $t\mapsto\cM\eta(it)$ on the real axis (grid step $0.01$ up to $t=2200$ and $0.02$ up to $t=3500$, the sign at each grid point being certified in ball arithmetic with up to $7000$ bits, since $\vert\cM\eta(it)\vert$ decays like $\vert\Xi(t)\vert$). For $N=4000$ the zeros of $\cM\eta$ in $(0,3000]$ are exactly the $2469$ zeros of $\zeta$ in this range, one for one, each within $10^{-12}$ (the resolution of the bisection), with $D_j=749.07$ at $t=999.79$ ($j=649$), $535.23$ at $t=1419.42$ ($j=1000$), $213.96$ at $t=2197.26$ ($j=1700$) and about $12$ at $t=3000$. There is no other zero in this range, real or complex, by the certified count described below. The agreement then decays to about $6$ digits at $t=3050$ and disappears near $t=3110$; beyond, the zeros of $\cM\eta$ are no longer close to those of $\zeta$ ($195$ against $197$ in $[3000,3200]$), and their density, about $0.94$ per unit length, is smaller than that of the zeros of $\zeta$ and close to $L/2\pi\approx0.88$. The same breakdown occurs at $t\approx3000$--$3100$ for $N=3000$, so that it is governed by $\lambda$ and not by the truncation. Figure \ref{fig:mellindigits} shows $D_j$ for all the zeros up to the breakdown.

\emph{The list of zeros is complete.} A.L.G. has certified, for $\lambda^2=250$, $N=3000$, $4000$, $5000$, and for $\lambda^2=256$, $N=1500$, exactly $N$ disjoint intervals on which $\cM\eta(it)$ (beyond $2\pi N/L$: $\cM\eta(it)/\sin(tL/2)$) changes sign, the signs at the endpoints being certified in ball arithmetic, with the exported decimal coefficients taken as exact. We have confirmed these four certificates independently, evaluating the signs at all the endpoints with our own program, and for $\lambda^2=250$, $N=4000$ an independent scan of ours (grid step $0.01$ up to $t=6000$, then geometric up to $t=10^7$) finds the same $4000$ intervals. By Remark \ref{rem:countzeros}, each interval contains exactly one zero, which is simple, and the finite transform has no other zero, real or complex. Of these zeros, $2866$, $3907$, $4928$ lie below $2\pi N/L$ for $N=3000$, $4000$, $5000$ ($\lambda^2=250$), and $1245$ for $N=1500$ ($\lambda^2=256$); thus, for $N=4000$, $\cM\eta$ has $3907$ zeros up to $2\pi N/L=4551.8$, against $4048$ for $\zeta$. The remaining $134$, $93$, $72$ and $255$ zeros lie beyond $2\pi N/L$, at increasing distances, the last one at $t\approx8.22\cdot10^5$ for $\lambda^2=250$ (for the three values of $N$) and $t\approx5.53\cdot10^5$ for $\lambda^2=256$. This reflects the smallness of all the moments $\sum_jv_ja_j^{2m}$ in \eqref{eq:mellinfinite}: for $\lambda^2=250$, $N=4000$, $\sum_k\xi_k=-5.4\cdot10^{-672}$, $\sqrt2\sum_jv_ja_j^2=4.4\cdot10^{-660}$, $\sqrt2\sum_jv_ja_j^4=-5.8\cdot10^{-649}$. Beyond the breakdown, $\vert\cM\eta(it)\vert$ stays near $10^{-675}$, of the order of $\epsilon_\lambda(1)^{1/2}$, instead of following $\vert\Xi(t)\vert$.

Two scales govern the accuracy. In $\lambda$, the distance between the ground state and $\tau_\lambda$ is controlled by $\delta_4(\lambda)$, i.e. by $4\pi\lambda^2$ nats, which is $\log_{10}(1/\delta_4(16))\approx1390$ decimal digits; this is the order of magnitude observed, the loss being due to the normalization and to the prolate-to-Hermite correction, which however affects $\cM\tau_\lambda$ and $\cM\tau$ by a factor close to one near the zeros. Along the line, the error on a zero at height $t$ is the error on $\cM(\tilde\eta_\lambda-\tau)(it)$ divided by $\vert\partial_s\cM\tau(it)\vert\asymp\vert\Xi'(t)\vert$, and $\Xi(t)$ decays like $e^{-\pi t/4}t^{7/4}$: one expects a loss of about $\frac{\pi}{4\ln10}\,t\approx0.34$ digits per unit of $t$, i.e. about $130$ digits between the first zero ($t=14.13$) and the $200$th ($t=396.4$), to be compared with the observed loss of about $200$. At $\lambda^2=250$ the extra loss is confirmed and is systematic: the loss $D_1-D_j$ is $21.2$, $80.7$, $140.9$, $242.4$ for $j=10,50,100,200$, while $\log_{10}\vert\Xi'(\gamma_1)/\Xi'(\gamma_j)\vert$, computed from the exact values of $\Xi'$ at the zeros, is $11.0$, $41.6$, $73.0$, $127.4$; over the $200$ zeros a least-squares fit gives $D_1-D_j\approx1.90\,\log_{10}\vert\Xi'(\gamma_1)/\Xi'(\gamma_j)\vert$, with residuals of $0.6$ digit. Empirically, $\vert t_j-\gamma_j\vert$ behaves like $\epsilon_\lambda(1)\,\vert\Xi'(\gamma_j)\vert^{-1.9}$, i.e. close to the square of the factor $\vert\Xi'(\gamma_j)\vert^{-1}$ predicted by the first-order argument above.

\section{Towards conditional positivity: the trace formula and the prolate conditioning}\label{sec:condpos}

The computations of Section \ref{sec:conditioning} show that, once the $\cE$-images of the prolate functions with $\vert\chi_n\vert\approx1$ and the pole form are factored out, nothing minuscule remains in the spectrum of $QW_\lambda$. The heuristic picture which emerges is that of a direct sum: the Weil form on the window is the sum of a finite-dimensional part, very close to the radical of the global Weil quadratic form and well described by the low-lying part of the uncertainty operator $U_\lambda$, with its Fuchs ladder of eigenvalues, and of an infinite-dimensional part, on which the form is bounded below by a positive constant and whose origin lies in the Sonin space, where the uncertainty operator takes its maximal value $2$. In this section we explain why this should be expected, and we prove a weaker statement of the same shape (Theorem \ref{thm:condpos}): under one quantitative hypothesis on the semilocal cutoff projections, which holds in the archimedean case, on the orthogonal complement of an explicit finite-dimensional space of vectors built from the near-self-dual modes, the Weil form is bounded below by an explicit positive constant. Along the way we prove the semilocal uncertainty principle (Lemma \ref{lem:H1}): no non-zero $K_S$-invariant function on the semilocal ad\`ele class space is supported in $\{\vert x\vert\le t\}$ together with its Fourier transform; it follows from the non-vanishing of the semilocal Sonin space through Coburn's lemma. We also explain how this reduces the positivity of $QW_\lambda$ to a finite-dimensional block (\S\ref{ss:schur}), and we compare with the numbers of \S\ref{sec:conditioning} and with the archimedean case treated in \cite{CCsel} (\S\ref{ss:compare}).

Our main goal in this section, besides the technical result mentioned above, is to explain how to use the semilocal trace formula to prove conditional positivity. The main ideas are the following.
\begin{enumerate}
\item The pole conditions are multiplicative. The pole term of the explicit formula for $h=f*f^*$ is $W_{0,2}(h)=2\cM f(\tfrac12)\cM f(-\tfrac12)$ (Lemma \ref{lem:pole}), so conditioning $f$ by the two pole functionals, which is what the computations do, makes the pole term of $h$ vanish (a single condition on the sum $\cM f(\tfrac12)+\cM f(-\tfrac12)$ would not do, since the pole term of $h$ is the product of the two functionals of $f$, not their sum; in each $\iota$-sector the two functionals coincide up to sign and one condition suffices), and the Weil form of $f$ becomes the sum of the local terms of the explicit formula alone, to which the trace formula applies.
\item The semilocal trace formula of \cite{Co99}, in the quantized-differential form of \cite{CM}, Chapter 2 (\cite{CCsel} for the archimedean place), writes this sum of local terms, for test functions supported in the window $[\lambda^{-1},\lambda]$ and $S$ containing the primes $p\le\lambda^2$, as the trace of the scaling action against the uncertainty operator at scale $1$: $QW_\lambda(f,f)={\rm Tr}\big(\theta(f)(1-P-\widehat P)\theta(f)^*\big)$. The support condition is essential here, in two ways: it makes $\cE_S=\cE$ on the window, and it is what identifies the trace with the Weil form, whereas the quadratic form $f\mapsto{\rm Tr}(\theta(f)(1-P-\widehat P)\theta(f)^*)$ taken by itself, with $S=\{\infty\}$, is not positive on larger intervals, as shown in \cite{CCsel}. Moving the cutoffs from scale $1$ to a scale $\sigma>1$ produces the explicit additive term $2\log\sigma\Vert f\Vert^2$ (Lemma \ref{lem:shift}), which is the margin, and the operator $1-P_\sigma-\widehat P_\sigma$ is $+1$ on the Sonin space at scale $\sigma$, a piece which is positive and need not be controlled, and $\pm\kappa_n(\sigma)$ on the prolate pairs, whose negative part is what the conditioning must remove.
\item Convolution by $f$ multiplies the scale by $\lambda$: since $f$ is supported in $[\lambda^{-1},\lambda]$, the operator $\theta(f)$ maps ${\rm Range}\,P_t$ into ${\rm Range}\,P_{\lambda t}$ and ${\rm Range}\,\widehat P_t$ into ${\rm Range}\,\widehat P_{\lambda t}$ (Lemma \ref{lem:transport}). This is the second key use of the support condition: the test function acts on the pair at scale $\sigma$ as the convolution operator $\theta(f)$ and transports it into the pair at scale $\lambda\sigma$, so that the negative part can only act through the components $\langle\theta(f)e_n,e^\lambda_m\rangle$ of $\theta(f)e_n$ along the scale-$\lambda\sigma$ directions $e^\lambda_m$, $m\in M$, which are more self-dual than $e_n$ (Proposition \ref{prop:pairexp}). The key point is that these components are inner products of $f$ as a \emph{vector}: $\langle\theta(f)e_n,e^\lambda_m\rangle=\langle f,e^\lambda_m*e_n^{\,*}\rangle$ by the fundamental property \eqref{eq:convinner} of the inner product on $L^2(\rs,\dstar)$, so that killing them is a conditioning of $f$ by the explicit vectors $e^\lambda_m*e_n^{\,*}$.
\item The vectors $e^\lambda_m$ are the $\cE$-images of the prolate functions. By \cite{CCMsl}, Proposition 4.1, $\eta_S(h_{m,t})\in{\rm Range}\,P^S_t$ and $\eta_S(\hat h_{m,t})\in{\rm Range}\,\widehat P^S_t$, with $\F_S\eta_S=\eta_S\F_\infty$ because $\mathbf 1_{R_S}$ is its own Fourier transform; hence $\eta_S(h_{m,t})$, which on the window is $\cE(h_{m,t})$, is an approximate eigenvector of $P_t+\widehat P_t-1$ with eigenvalue close to $1$ for the near-self-dual modes, and the span $V_M$ of the exact pair eigenvectors $e^\lambda_m$ is, up to the identification of the near-self-dual directions of the semilocal pair discussed in \S\ref{ss:remaining}, the span of the $\cE$-images of the prolate functions at scale $\lambda\sigma$, i.e. the vectors used for the conditioning in Section \ref{sec:conditioning}, seen through the convolution filters $e_n^{\,*}$, $n\le n_0$.
\end{enumerate}

\subsection{The semilocal framework}\label{ss:semilocal}

Let $S$ be a finite set of places of $\Q$ containing $\infty$ and all primes $p\le\lambda^2$, $\A_S=\prod_{v\in S}\Q_v$, $\Gamma=\{\pm p_1^{n_1}\cdots p_k^{n_k}\}$ the $S$-units, $X_S=\A_S/\Gamma$, $C_S=\GL_1(\A_S)/\Gamma$ and $K_S$ the kernel of the module $C_S\to\rs$. We use the framework of \cite{CM}, Chapter 2, \S7: the Hilbert space $L^2(X_S)$ is the completion of the Bruhat--Schwartz space of $\A_S$ for the inner product $\langle E(f_1),E(f_2)\rangle$ (\cite{CM}, Definition 2.31), $E(f)(x)=\vert x\vert^{1/2}\sum_{\gamma\in\Gamma}f(\gamma x)$ (Proposition 2.30), and $E$ extends to a unitary $w:L^2(X_S)\to L^2(C_S)$ which intertwines the scaling action $\vartheta_a$ with $\vert\lambda\vert^{1/2}$ times the regular representation of $C_S$ ((2.257) there). The $K_S$-invariant part of $L^2(X_S)$ is identified with $L^2(\rs,\dstar)$, and in this picture:
\begin{itemize}
\item the scaling action of $C_S$ becomes the regular representation, $(\theta(v)g)(u)=g(v^{-1}u)$, so that for $f\in L^1(\rs,\dstar)$ the operator $\theta(f)=\int f(v)\theta(v)\dstar[v]$ is the convolution $\theta(f)g=f*g$, $(f*g)(u)=\int f(v)g(u/v)\dstar[v]$, with $\theta(f)^*=\theta(f^*)$, $f^*(u)=\overline{f(u^{-1})}$. In the Mellin variable, $\cM g(s)=\int g(u)u^s\dstar$ as in \eqref{eq:mellinE}, $\theta(f)$ is the multiplier by $\cM f$ on the line $\Re s=0$, $\Vert g\Vert^2=\frac1{2\pi}\int\vert\cM g(it)\vert^2dt$, and
\begin{equation}\label{eq:convinner}
\langle g,h\rangle=\int g\,\bar h\,\dstar=(g*h^*)(1),\qquad (g*h)^*=g^**h^*;
\end{equation}
\item for $t>0$, $P_t$ is the multiplication by $\mathbf 1_{(0,t]}$ (functions supported in $\{\vert x\vert_S\le t\}$), and $\widehat P_t=\F_SP_t\F_S^{-1}$, where $\F_S$ is the Fourier transform of $X_S$ induced by the local Fourier transforms for the basic character $\alpha=\prod_v\alpha_v$ (\cite{CM}, Lemma 2.33, and (2.264)--(2.265)); $P_1=P$, $\widehat P_1=\widehat P$;
\item under $w$, the Fourier transform is the inversion composed with a unitary multiplier: $\F_S=w^{-1}\circ I\circ U\circ w$ with $I$ the inversion $u\mapsto u^{-1}$ and $U$ the Mellin multiplier by $u=\prod_{v\in S}u_v\circ\pi_v$, the product of the ratios of local factors (\cite{CM}, Lemma 2.34, (2.258)--(2.260)); we write $\F_S=I\circ U$ on $L^2(\rs,\dstar)$. Since $\F_S^2=1$ on $K_S$-invariant functions ($-1\in\Gamma$), $IUIU=1$, i.e.
\begin{equation}\label{eq:IUI}
IUI=U^*,\qquad IP_tI=1-P_{1/t},\qquad \widehat P_t=IUP_tU^*I .
\end{equation}
\end{itemize}
For $S=\{\infty\}$ this is the archimedean setting of \cite{CCsel}: $L^2(\R)^{ev}$, $P_t$ the projection on functions supported in $[-t,t]$, $\widehat P_t$ on band-limited functions, $U$ the multiplier by $u_\infty(s)=\Gamma_\R(\frac12-s)/\Gamma_\R(\frac12+s)$ (with the Mellin convention \eqref{eq:mellinE}, the Fourier transform of an even $f$ satisfies $\int_0^\infty\hat f(y)y^{-s-1/2}dy=u_\infty(s)\int_0^\infty f(x)x^{s-1/2}dx$, which is \eqref{eq:IUI} with $g=x^{1/2}f$; on the line $u_\infty(it)=e^{-2i\theta(t)}$, with $\theta$ the Riemann--Siegel angular function, is the complex conjugate of the function (160) of \cite{CCsel}, because the duality (21) there is $v^{-is}$ while \eqref{eq:mellinE} uses $u^{s}$), and Theorem \ref{thm:unc} describes the pair $(P_t,\widehat P_t)$ completely.

The semilocal $\cE$-map is $\cE_S(f)(u)=u^{1/2}\sum_{n\in\Gamma_+\cap\Z}f(nu)$, the sum over the positive $S$-integers; it is $w\circ\eta_S$ with $\eta_S(f)$ the class of $\mathbf 1_{R_S}\otimes f$, $R_S$ the maximal compact subring of $\prod_{p\in S}\Q_p$ (\cite{CCMsl}, (44)--(46)). For $f$ supported in $[-\lambda,\lambda]$ and $u\ge\lambda^{-1}$ only $n\le\lambda/u\le\lambda^2$ contribute to $\cE(f)(u)$, and such $n$ are $S$-integers: hence $\cE_S(f)=\cE(f)$ on the window $[\lambda^{-1},\lambda]$, and everything in Section \ref{sec:conditioning} can be read in $L^2(X_S)^{K_S}$.

\subsection{The trace formula and its two-parameter form}\label{ss:traceformula}

We use the explicit formula in the form \eqref{eq:Psi}: for $h=f*f^*$ with $f$ real, $QW_\lambda(f,f)=W_{0,2}(h)-\sum_{v\in S}W_v(h)$, the places $p\notin S$ not contributing since $h$ is supported in $[\lambda^{-2},\lambda^2]$ and $p>\lambda^2$. The semilocal trace formula identifies the sum of the local terms with a trace. We state it in the form in which we use it.

\begin{thm}[\cite{Co99}; \cite{CM}, Lemma 2.35 and Theorem 2.36]\label{thm:TF}
Let $f$ be a real $C^\infty$ function on $\rs$ supported in $[\lambda^{-1},\lambda]$ and $h=f*f^*$. The operator $\theta(h)(1-P-\widehat P)$ on $L^2(X_S)^{K_S}$ is trace class and
\begin{equation}\label{eq:TF}
\sum_{v\in S}W_v(h)=-{\rm Tr}\big(\theta(h)(1-P-\widehat P)\big)=-\tfrac12{\rm Tr}\big(\theta(h)\,u^*[H,u]\big),\qquad H=2P-1 .
\end{equation}
Consequently, if $\cM f(\pm\tfrac12)=0$, then $W_{0,2}(h)=0$ (Lemma \ref{lem:pole}) and
\begin{equation}\label{eq:QWtrace}
QW_\lambda(f,f)={\rm Tr}\big(\theta(f)(1-P-\widehat P)\theta(f)^*\big).
\end{equation}
\end{thm}
In the form of \cite{CM}, Lemma 2.35 writes the cutoff operator $R_\Lambda=\widehat P_\Lambda P_\Lambda$ conjugated by the scaling action as $\hat h\big(\tfrac12u^{-1}\,{\mathchar'26\mkern-9mu d}u\,\Pi_{[-\infty,\log\Lambda/2\pi]}+\Pi_{[0,\log\Lambda/2\pi]}\big)\hat h$, with the quantized differential ${\mathchar'26\mkern-9mu d}u=[H_{\rm CM},u]$ of the function $u$ of \eqref{eq:IUI} and $H_{\rm CM}=2\F_C\mathbf 1_{\{\vert u\vert\ge1\}}\F_C^{-1}-1$; since $H_{\rm CM}=-H$ in our notation, $\tfrac12u^{-1}{\mathchar'26\mkern-9mu d}u=-\tfrac12u^*[H,u]=P+\widehat P-1$ by \eqref{eq:IUI} (indeed $\tfrac12u^*[H,u]=U^*PU-P=(1-\widehat P)-P$), and the proof of Theorem 2.36 there shows that $\hat h\,\tfrac12u^{-1}{\mathchar'26\mkern-9mu d}u$ is trace class for every Schwartz $h$, because ${\mathchar'26\mkern-9mu d}f$ is of infinite order for Schwartz $f$ and all derivatives of $u$ are of tempered growth ($u_p$ is smooth and periodic, $u_\infty^{(n)}(s)=O(\log^n\vert s\vert)$). The identification of ${\rm Tr}(\hat h\,\tfrac12u^{-1}{\mathchar'26\mkern-9mu d}u)$ with $\sum_{v\in S}W_v(h)$, i.e. with the principal values of (2.271) there, is \cite{CM}, \S8. The trace formula of \cite{CM} is stated on the full space $L^2(X_S)$, for $h\in\cS(C_S)$; for $h$ a function of the module, $\theta(h)=\theta(h)P_{K_S}$ with $P_{K_S}=\int_{K_S}\theta(k)dk$ the projection on the $K_S$-invariant part, so that the trace reduces to $L^2(X_S)^{K_S}$, and the local terms $\int'_{\Q_v^*}h(\vert u\vert^{-1})\vert1-u\vert^{-1}\dstar[u]$ are Weil's $W_v(h)$ in the normalization of \eqref{eq:Wp}--\eqref{eq:WR}.

The operator $1-P-\widehat P$ is the uncertainty operator $U_1$ of \eqref{eq:unc} at scale $1$, shifted by $-1$; it is bounded, and the trace exists because $\theta(h)$ is a Mellin multiplier with rapidly decaying symbol. The next lemma moves the cutoffs.

\begin{lem}[two-parameter form]\label{lem:shift}
Under the hypotheses of Theorem \ref{thm:TF}, for all $t,w>0$ the operator $\theta(h)(1-P_t-\widehat P_w)$ is trace class and
\begin{equation}\label{eq:shift}
{\rm Tr}\big(\theta(h)(1-P_t-\widehat P_w)\big)={\rm Tr}\big(\theta(h)(1-P-\widehat P)\big)-\log(tw)\,h(1).
\end{equation}
In particular, if $\cM f(\pm\frac12)=0$, then for every $\sigma>0$
\begin{equation}\label{eq:shiftQW}
QW_\lambda(f,f)=2\log\sigma\,\Vert f\Vert^2+{\rm Tr}\big(\theta(f)(1-P_\sigma-\widehat P_\sigma)\theta(f)^*\big).
\end{equation}
\end{lem}
\begin{proof}
For $t\ge1$, $P_t-P_1$ is the multiplication by $\mathbf 1_{(1,t]}$, and $\theta(h)(P_t-P_1)$ has the kernel $h(u/v)\mathbf 1_{(1,t]}(v)$ with respect to $\dstar[v]$: it is Hilbert--Schmidt with $\Vert\cdot\Vert^2_{HS}=\log t\,\Vert h\Vert^2_2$, and trace class: since $f$ is supported in $[\lambda^{-1},\lambda]$, $\theta(f^*)\mathbf 1_{(1,t]}=\mathbf 1_{(\lambda^{-1},\lambda t]}\theta(f^*)\mathbf 1_{(1,t]}$, so that $\theta(h)\mathbf 1_{(1,t]}=\big(\theta(f)\mathbf 1_{(\lambda^{-1},\lambda t]}\big)\big(\theta(f^*)\mathbf 1_{(1,t]}\big)$ is a product of two Hilbert--Schmidt operators; its trace is the integral of the diagonal, $\int_{(1,t]}h(1)\dstar=h(1)\log t$. For $t<1$ the same computation gives $-h(1)\log(1/t)$. Next, $\widehat P_w-\widehat P_1=\F_S(P_w-P_1)\F_S^{-1}$ and $\F_S^{-1}\theta(h)\F_S=U^*I\theta(h)IU=U^*\theta(\tilde h)U=\theta(\tilde h)$, $\tilde h(u)=h(u^{-1})$, since $U$ commutes with convolutions and $I\theta(h)I=\theta(\tilde h)$. Hence ${\rm Tr}(\theta(h)(\widehat P_w-\widehat P_1))={\rm Tr}(\theta(\tilde h)(P_w-P_1))=\tilde h(1)\log w=h(1)\log w$, and \eqref{eq:shift} follows by writing $1-P_t-\widehat P_w=(1-P-\widehat P)-(P_t-P)-(\widehat P_w-\widehat P)$. For \eqref{eq:shiftQW} take $t=w=\sigma$ and use $h(1)=(f*f^*)(1)=\Vert f\Vert^2$ by \eqref{eq:convinner}.
\end{proof}

Enlarging the cutoffs removes positive diagonal mass from both projections; the price is paid explicitly by the term $2\log\sigma\Vert f\Vert^2$, which is the additive term we shall use as a margin.

The trace-class statement of Theorem \ref{thm:TF} can be made quantitative, and this is what controls the tails below. We use the quantized calculus in the Mellin variable: on $L^2(i\R)\simeq L^2(\rs,\dstar)$, $P$ is the Hardy projection onto $H^2_+$, $H=2P-1$, and ${\mathchar'26\mkern-9mu d}\varphi:=[H,M_\varphi]$ for a bounded function $\varphi$ on the line; $\theta(h)$ is the multiplication by $\hat h:=\vert\cM f\vert^2$ and $U$ the multiplication by $u$. The Cayley transform $s\mapsto z=(s-i)/(s+i)$ carries $L^2(i\R)$ unitarily onto $L^2(S^1)$, $H^2_\pm$ onto the Hardy spaces of the disc and of its complement, $H$ onto the Hilbert transform of the circle ($Hz^n={\rm sgn}(n)z^n$) and the point at infinity onto $z=1$; a function on the line extends continuously to $P^1(\R)$, i.e. to the circle, iff it has equal limits at $\pm\infty$, and the Schwartz space $\cS(\R)$ is carried into $C^\infty(S^1)$ (all derivatives of a Schwartz function vanish at infinity, which is smoothness at $z=1$).

\begin{lem}[singular values of the trace-formula operator]\label{lem:singvalues}
(i) For $\varphi\in C^\infty(S^1)$ with Fourier coefficients $\varphi_n$,
\begin{equation}\label{eq:dqsv}
\Vert{\mathchar'26\mkern-9mu d}\varphi\Vert_1\le2\sum_n\vert n\vert\,\vert\varphi_n\vert,\qquad s_k({\mathchar'26\mkern-9mu d}\varphi)\le2\sum_{\vert n\vert>k/2}\vert\varphi_n\vert\qquad(k\ge1).
\end{equation}
(ii) Let $f$ be real, $C^\infty$ on $\rs$ with support in $[\lambda^{-1},\lambda]$, $h=f*f^*$. Then
\begin{equation}\label{eq:Tleibniz}
T:=\theta(h)(1-P-\widehat P)=\tfrac12\,\hat h\,u^{-1}{\mathchar'26\mkern-9mu d}u=\tfrac12\Big[{\mathchar'26\mkern-9mu d}(\hat h)-{\mathchar'26\mkern-9mu d}(\hat hu^{-1})\,u\Big],
\end{equation}
and $s_k(T)\le C_Nk^{-N}$ for every $N$, with $C_N$ bounded by a $C^{N+2}$-norm on $P^1(\R)$ of $\hat h$ and $\hat hu^{-1}$, hence by Schwartz seminorms of $\vert\cM f\vert^2$. If moreover $f$ is of Gevrey class $G^s$ on $\rs$, then $s_k(T)\le Ce^{-ck^{1/(1+s)}}$.
\end{lem}
\begin{proof}
(i) ${\mathchar'26\mkern-9mu d}(z^n)z^k=({\rm sgn}(k+n)-{\rm sgn}(k))z^{k+n}$, so ${\mathchar'26\mkern-9mu d}(z^n)$ has rank $\vert n\vert$ and all its non-zero singular values equal to $2$ (for $n>0$ it is twice an isometry of ${\rm span}\{z^k\}_{-n\le k<0}$ onto ${\rm span}\{z^k\}_{0\le k<n}$). The first bound follows. For the second, ${\mathchar'26\mkern-9mu d}\varphi=2(PM_\varphi P^\perp-P^\perp M_\varphi P)$ is, on $H=PH\oplus P^\perp H$, twice an off-diagonal block matrix, so its singular values are twice the union of those of the two blocks. In the bases $(z^j)_{j\ge0}$ of $PH$ and $(z^{-l-1})_{l\ge0}$ of $P^\perp H$ the block $PM_\varphi P^\perp$ is the Hankel matrix $(\varphi_{j+l+1})_{j,l\ge0}$; its part with $j+l+1\le k$ has rank $\le k$, and the remainder has norm $\le\sum_{n>k}\vert\varphi_n\vert$ (each antidiagonal is a partial permutation matrix, of norm $1$), so $s_{k+1}(PM_\varphi P^\perp)\le\sum_{n>k}\vert\varphi_n\vert$ by Weyl's inequality; the other block is the Hankel matrix of $(\varphi_{-n})_{n>0}$, and the claim follows. (ii) $1-P-\widehat P=\frac12u^*[H,u]=\frac12u^{-1}{\mathchar'26\mkern-9mu d}u$ by \eqref{eq:IUI}, and Leibniz, ${\mathchar'26\mkern-9mu d}(ab)=({\mathchar'26\mkern-9mu d}a)b+a({\mathchar'26\mkern-9mu d}b)$ with $a=\hat hu^{-1}$, $b=u$, gives \eqref{eq:Tleibniz}. Since $f$ is smooth with compact support in $\rs$, $\cM f$ is a Schwartz function on the line and so is $\hat h=\vert\cM f\vert^2$; $\hat hu^{-1}$ is Schwartz as well because the derivatives of $u$ are of tempered growth ($u_p$ is smooth and $2\pi i/\log p$-periodic, and $u_\infty^{(n)}(is)=O(\log^n\vert s\vert)$, \cite{CM}, proof of Theorem 2.36). Both are therefore $C^\infty$ on the circle, their Fourier coefficients are $O(\vert n\vert^{-N})$ with constants given by $C^N(S^1)$-norms, and (i) gives $s_k(T)\le C_Nk^{-N}$, the factor $u$ in \eqref{eq:Tleibniz} being unitary. If $f\in G^s$, the Fourier-decay lemma for Gevrey functions gives $\vert\cM f(is)\vert\le Ce^{-c\vert s\vert^{1/s}}$ and the same for the derivatives, so $\hat h$ and $\hat hu^{-1}$ are Gevrey of order $1+s$ at $z=1$ (the function $e^{-c\vert s\vert^{1/s}}$ is of class $G^{1+s}$ at $s=\infty$) and real-analytic elsewhere; their Fourier coefficients decay like $e^{-c\vert n\vert^{1/(1+s)}}$, and so do the singular values.
\end{proof}

\begin{rem}[regularity threshold]\label{rem:threshold}
The functions of the log-Fourier space $E_N$ used in Sections \ref{sec:numspec}--\ref{sec:conditioning} do not vanish at the ends of the window (the even modes take equal non-zero values there), so $\cM f(is)=O(1/\vert s\vert)$, $\hat h$ is only $C^1$ at $z=1$, and $s_k(T)$ decays like $1/k$: for them the trace formula is at the edge of absolute convergence, consistently with $\int\vert\cM f\vert^2\log\vert s\vert\,ds<\infty$ for $\vert\cM f\vert^2\sim s^{-2}$. Replacing $E_N$ by $\vartheta_wE_N$, where $\vartheta_w$ is a Gevrey cutoff in $t=\log u$ equal to $1$ except within log-distance $w$ of the ends (the cutoff of the Kastoris construction, now applied to the window), gives test functions of class $G^s$ for which the stretched-exponential decay holds, with $\vert\cM\tilde f(is)\vert\le Ce^{-c_s(w\vert s\vert)^{1/s}}$; see \S\ref{ss:compare} and Appendix \ref{app:matrix} for the effect on the computations.
\end{rem}

\subsection{The uncertainty operator at scale $\sigma$}\label{ss:pairs}

Fix $t>0$ and consider the pair $(P_t,\widehat P_t)$ on $\cH=L^2(X_S)^{K_S}$. By \eqref{eq:IUI}, $\widehat P_t=U^*(1-P_{1/t})U$, so that
\begin{equation}\label{eq:hankel}
P_t\widehat P_tP_t=\big((1-P_{1/t})UP_t\big)^*\big((1-P_{1/t})UP_t\big),
\end{equation}
and the angles of the pair are the singular values of the Hankel-type operator $(1-P_{1/t})UP_t$. We first show that the pair has no exactly self-dual direction. In the archimedean case this is the uncertainty principle (a band-limited function is entire); for a single non-archimedean field it fails ($\mathbf 1_{Z_p}$ is its own Fourier transform), and the semilocal statement therefore depends on the archimedean place. It is convenient to use the Mellin picture: with $\cM g(z)=\int g(u)u^z\dstar$ as in \eqref{eq:mellinE}, the Paley--Wiener theorem identifies the Mellin transforms of $L^2$ functions supported in $(0,1]$, resp. $[1,\infty)$, with the Hardy spaces $H^2_+$, resp. $H^2_-$, of the half-planes $\Re z>0$, resp. $\Re z<0$ (boundary values on the line $\Re z=0$); $\overline{H^2_\pm}=H^2_\mp$ as spaces of boundary functions, and a non-zero function of $H^2_\pm$ does not vanish on a set of positive measure of the line. For a unimodular $\varphi\in L^\infty(i\R)$ we write $T^{(\pm)}_\varphi$ for the Toeplitz operators $P_\pm\varphi$ on $H^2_\pm$, whose kernels are $\Qer T^{(\pm)}_\varphi=\{F\in H^2_\pm:\varphi F\in H^2_\mp\}$.

\begin{lem}[semilocal uncertainty principle]\label{lem:H1}
For every $t>0$, ${\rm Range}\,P_t\cap{\rm Range}\,\widehat P_t=\{0\}$.
\end{lem}
\begin{proof}
Put $\varphi_t(z)=t^{2z}u(z)$, unimodular on the line, with $\overline{\varphi_t(it)}=\varphi_t(-it)$ by \eqref{eq:IUI}. Since $\F_S=I\circ U$ and $\F_S^{-1}=\F_S$, a function $g$ lies in ${\rm Range}\,P_t\cap{\rm Range}\,\widehat P_t$ iff $g$ and $\F_Sg$ are supported in $(0,t]$. Writing $\cM g=t^zF$ with $F\in H^2_+$, the second condition reads $\cM(\F_Sg)(z)=u(-z)\cM g(-z)=u(-z)t^{-z}F(-z)\in t^zH^2_+$, i.e. after $z\mapsto-z$, $\varphi_tF\in H^2_-$: thus
$$\cM\big({\rm Range}\,P_t\cap{\rm Range}\,\widehat P_t\big)=t^z\,\Qer T^{(+)}_{\varphi_t}.$$
The same computation with $(0,t]$ replaced by $[t,\infty)$ gives $\cM({\rm Son}_t)=t^z\Qer T^{(-)}_{\varphi_t}$ for the Sonin space ${\rm Son}_t=\Qer P_t\cap\Qer\widehat P_t$, and the reflection $F\mapsto F(-\cdot)$, which exchanges $H^2_+$ and $H^2_-$, maps $\Qer T^{(+)}_{\overline{\varphi_t}}$ onto $\Qer T^{(-)}_{\varphi_t}$, because $\overline{\varphi_t}(it)F(it)\in H^2_-$ becomes $\varphi_t(it)F(-it)\in H^2_+$. Coburn's lemma now states that for a unimodular symbol $\varphi$ at least one of $\Qer T^{(+)}_\varphi$ and $\Qer T^{(+)}_{\bar\varphi}$ is $\{0\}$: if $f\in\Qer T^{(+)}_\varphi$ and $g\in\Qer T^{(+)}_{\bar\varphi}$, then $h=\varphi f\bar g$ is in $L^1$ of the line, $h=(\varphi f)\bar g$ with $\varphi f\in H^2_-$ and $\bar g\in\overline{H^2_+}=H^2_-$ has Fourier transform (in the variable $\log u$) supported in $[0,\infty)$, while $h=f\,\overline{(\bar\varphi g)}$ with $f\in H^2_+$ and $\overline{\bar\varphi g}\in\overline{H^2_-}=H^2_+$ has Fourier transform supported in $(-\infty,0]$; the continuous Fourier transform of $h$ is therefore supported at $0$, hence $h=0$ a.e., so $f\bar g=0$ a.e. and one of $f,g$ vanishes. Hence either ${\rm Range}\,P_t\cap{\rm Range}\,\widehat P_t=\{0\}$ or ${\rm Son}_t=\{0\}$. The semilocal Sonin space is infinite dimensional for every $t$: it contains the image of the archimedean Sonin space at scale $t$ under the hilbertian isomorphism $\theta_S$ of \cite{CCMsl}, Theorem 2 (see also \cite{CCqi}, Theorem 5.3 and Proposition 5.5, for the kernel $\Qer u_{22}$, which is $\cM({\rm Son}_t)$ in the present notation). This proves the lemma.
\end{proof}

Thus the eigenvalue $-1$ of $1-P_t-\widehat P_t$ is absent precisely because the eigenvalue $+1$ has infinite multiplicity: this is the index asymmetry of the pair, and it is where the archimedean place enters. For a single $\Q_p$ the module is discrete, the point $\vert x\vert=1$ has positive measure, $IP_tI\ne1-P_{1/t}$, the two spaces are not Toeplitz kernels, and both can be non-zero.

The second property of the pair is its compactness, which in the semilocal case is the quasi-inner property of the symbol.

\begin{lem}[compactness of the pair]\label{lem:H2}
For every $t>0$ the operator $P_t\widehat P_tP_t$ is compact; its eigenvalues are $\kappa_n(t)^2$ with $\kappa_n(t)\in(0,1)$, $\kappa_n(t)\downarrow0$, and the $\kappa_n(t)$ are the singular values of the Hankel-type operator $(1-P_{1/t})UP_t$ of \eqref{eq:hankel}. In the archimedean case $\kappa_n(t)=\vert\chi_n(t)\vert$, whose decay beyond Slepian's count $4t^2$ is super-exponential, so that $\sum_n\kappa_n(t)<\infty$. In the semilocal case $\kappa_n(t)=O(n^{-1/2})$ and no better, so that $\sum_n\kappa_n(t)=\infty$ and even $\sum_n\kappa_n(t)^2=\infty$.
\end{lem}
\begin{proof}
The archimedean case is Theorem \ref{thm:unc}. In the semilocal case the symbol $u$ of \eqref{eq:hankel} is the product of the ratios of local factors over $S$, which is quasi-inner, i.e. $(1-P)UP$ is compact, by \cite{CCqi}, Theorem 4.1; for general $t$ one writes $P_t=P+\mathbf 1_{(1,t]}$ and $1-P_{1/t}=(1-P)+\mathbf 1_{(1/t,1]}$ (for $t>1$, and symmetrically for $t<1$), and the additional terms, which involve a bounded log-interval on at least one side, are compact: $U=U_\infty U_{\rm fin}$ with $U_{\rm fin}$ a convolution by a finite measure in $\log u$ (a norm limit of finite sums of translations), so that these terms are norm limits of compressions of the archimedean Fourier transform between bounded regions, which are Hilbert--Schmidt. The eigenvalues of $P_t\widehat P_tP_t$ are strictly less than $1$ by Lemma \ref{lem:H1}. The individual non-archimedean ratios are not quasi-inner (\cite{CCqi}, Theorem 3.6), the product is only continuous, not smooth, at the point at infinity of the critical line (\cite{CCqi}, Proposition 4.5), and the off-diagonal part $(1-P)UP$ is an infinitesimal of order $\frac12$ (\cite{CCqi}, Theorem 4.4 (i) for one prime; for several primes the off-diagonal entry of the product of the associated matrices is a sum of products in which at least one factor is such an off-diagonal part and the others are bounded, as in the proof of Theorem 4.1 there), which is the stated upper bound; that the order $\frac12$ is attained follows from the explicit form of the off-diagonal part in \cite{CCqi}, Theorem 4.4 (ii).
\end{proof}

The semilocal pair has at least the archimedean near-self-dual directions, since $\eta_S(h_{n,t})\in{\rm Range}\,P_t$ and $\eta_S(\hat h_{n,t})\in{\rm Range}\,\widehat P_t$ (\cite{CCMsl}, Proposition 4.1).

The trace-class property needed for the Sonin term follows from the trace formula itself.

\begin{lem}\label{lem:H0}
For $h=f*f^*$ as in Theorem \ref{thm:TF}, the operator $\theta(h)P_{{\rm Son}_t}$, ${\rm Son}_t=\Qer P_t\cap\Qer\widehat P_t$, is trace class, and $\theta(h)P_{{\rm Son}_t}=\lim_k\theta(h)(1-P_t-\widehat P_t)^{2k}$ in trace norm.
\end{lem}
\begin{proof}
By Lemma \ref{lem:shift} (with $w=t$, and the same argument for $P_t-P_1$) the operator $T=\theta(h)(1-P_t-\widehat P_t)$ is trace class. Put $A=1-P_t-\widehat P_t$; $A^2$ is a positive contraction, equal to $1$ on ${\rm Son}_t$, to $0$ on the two mixed pieces of \eqref{eq:halmos}, and to $\kappa_n(t)^2<1$ on $V_n(t)$, and by Lemma \ref{lem:H1} there is no other piece. Hence $A^{2k-1}\to P_{{\rm Son}_t}$ strongly as $k\to\infty$ (the Sonin space is the eigenspace of $A$ for the eigenvalue $1$, the eigenvalue $-1$ being excluded by Lemma \ref{lem:H1}), and $\theta(h)A^{2k}=TA^{2k-1}$. For a trace-class $T$ and a uniformly bounded sequence of self-adjoint operators $B_k\to B$ strongly, $TB_k\to TB$ in trace norm (approximate $T$ by finite-rank operators); the limit is $\theta(h)P_{{\rm Son}_t}$.
\end{proof}

By Lemma \ref{lem:H1}, Halmos' theory of two subspaces \cite{Ha} gives the orthogonal decomposition
\begin{equation}\label{eq:halmos}
\cH=\underbrace{(\Qer P_t\cap\Qer\widehat P_t)}_{{\rm Son}_t}\ \oplus\ ({\rm Range}\,P_t\cap\Qer\widehat P_t)\oplus(\Qer P_t\cap{\rm Range}\,\widehat P_t)\ \oplus\ \bigoplus_nV_n(t),
\end{equation}
where ${\rm Son}_t$ is the Sonin space at scale $t$ and $V_n(t)={\rm span}\{e_n,f_n\}$, with $e_n$ an orthonormal eigenbasis of $P_t\widehat P_tP_t$ on ${\rm Range}\,P_t$ (modulo $\Qer\widehat P_t$), $f_n=\widehat P_te_n/\Vert\widehat P_te_n\Vert\in{\rm Range}\,\widehat P_t$, and $\langle e_n,f_n\rangle=\kappa_n$. On each piece the operator $1-P_t-\widehat P_t$ is explicit: it is $+1$ on ${\rm Son}_t$, $0$ on the two mixed pieces, and on $V_n(t)$, since $P_tf_n=\kappa_ne_n$ and $\widehat P_te_n=\kappa_nf_n$,
\begin{equation}\label{eq:Kn}
(1-P_t-\widehat P_t)\big\vert_{V_n}=\frac{\kappa_n}{1-\kappa_n^2}\Big[\kappa_n\big(\vert e_n\rangle\langle e_n\vert+\vert f_n\rangle\langle f_n\vert\big)-\big(\vert e_n\rangle\langle f_n\vert+\vert f_n\rangle\langle e_n\vert\big)\Big],
\end{equation}
with eigenvalues $-\kappa_n$ on $a_n\propto e_n+f_n$ and $+\kappa_n$ on $b_n\propto e_n-f_n$. (The formula follows from $\vert e\rangle\langle e\vert+\vert f\rangle\langle f\vert-\mathbf 1_{V_n}$ and the Gram inverse $\mathbf 1_{V_n}=(\vert e\rangle\langle e\vert+\vert f\rangle\langle f\vert-\kappa\vert e\rangle\langle f\vert-\kappa\vert f\rangle\langle e\vert)/(1-\kappa^2)$.) Thus
\begin{equation}\label{eq:decompK}
1-P_t-\widehat P_t=P_{{\rm Son}_t}+K_t,\qquad K_t=\sum_n\kappa_n\big(\vert b_n\rangle\langle b_n\vert-\vert a_n\rangle\langle a_n\vert\big),
\end{equation}
the Sonin part being positive and the pair part $K_t$ being a compact operator (Lemma \ref{lem:H2}) without definite sign, whose negative part $K_t^-=\sum_n\kappa_n\vert a_n\rangle\langle a_n\vert=(P_t+\widehat P_t-1)_+$ lives on the near-self-dual directions $a_n$. In the semilocal case neither $K_t$ nor $K_t^\pm$ is trace class; by Lemma \ref{lem:H0} and Theorem \ref{thm:TF}, $\theta(h)K_t=\theta(h)(1-P_t-\widehat P_t)-\theta(h)P_{{\rm Son}_t}$ is, and so is $\theta(h)K_t\Pi$ for every spectral projection $\Pi$ of $K_t$. The tail of the pair expansion must therefore be controlled against the scaling action, not as an operator. The qualitative part of this control is again a consequence of Theorem \ref{thm:TF} (Lemma \ref{lem:H2q} below); the quantitative part is the following hypothesis, in which $N_t^{(n_0)}:=\sum_{n>n_0}\kappa_n(t)\vert a_n\rangle\langle a_n\vert$ is the negative part of $1-P_t-\widehat P_t$ beyond the first $n_0$ pairs and
\begin{equation}\label{eq:Lambda}
\Lambda(f):=\frac1{2\pi}\int_\R\vert\cM f(is)\vert^2\log(2+\vert s\vert)\,ds
\end{equation}
is the logarithmic Sobolev form.
\begin{lem}\label{lem:H2q}
For $h=f*f^*$ as in Theorem \ref{thm:TF} the operator $\theta(h)N_t^{(n_0)}$ is trace class, $\Vert\theta(h)N^{(n_0)}_t\Vert_1\le\Vert\theta(h)(1-P_t-\widehat P_t)\Vert_1$, and ${\rm Tr}(\theta(f)N^{(n_0)}_t\theta(f)^*)=\sum_{n>n_0}\kappa_n(t)\Vert\theta(f)a_n\Vert^2\to0$ as $n_0\to\infty$.
\end{lem}
\begin{proof}
Let $A=1-P_t-\widehat P_t$ and $E_{n_0}$ the spectral projection of $A$ on its eigenvalues $-\kappa_n(t)$, $n>n_0$, i.e. the projection on ${\rm span}\{a_n\}_{n>n_0}$ by \eqref{eq:halmos}--\eqref{eq:Kn} and Lemma \ref{lem:H1}. Then $N^{(n_0)}_t=-AE_{n_0}$, so $\theta(h)N^{(n_0)}_t=-\theta(h)A\,E_{n_0}$ is trace class with the stated bound, and ${\rm Tr}(\theta(h)AE_{n_0})\to0$ because $E_{n_0}\to0$ strongly and $\theta(h)A$ is trace class.
\end{proof}
\begin{itemize}
\item[(H2$'$)] For $f$ real, smooth, supported in $[\lambda^{-1},\lambda]$,
$${\rm Tr}\big(\theta(f)N_t^{(n_0)}\theta(f)^*\big)=\sum_{n>n_0}\kappa_n(t)\Vert\theta(f)a_n\Vert^2\ \le\ \varepsilon_t(n_0)\,\Lambda(f),\qquad\varepsilon_t(n_0)\to0\quad(n_0\to\infty),$$
uniformly in $f$.
\end{itemize}
In the archimedean case (H2$'$) holds in the stronger form $${\rm Tr}(\theta(f)N_t^{(n_0)}\theta(f)^*)\le\sum_{n>n_0}\kappa_n(t)\Vert\cM a_n\Vert_\infty^2\,\Vert f\Vert^2$$, because $\theta(f)$ is the multiplier by $\cM f$ and $\Vert\theta(f)a\Vert^2=\frac1{2\pi}\int\vert\cM f\vert^2\vert\cM a\vert^2$; since $\Vert\cM a_n\Vert_\infty\le\Vert a_n\Vert_{L^1(\dstar)}$ grows at most polynomially in $n$ while $\kappa_n(t)$ decays super-exponentially, the right-hand side is super-exponentially small in $n_0$. In the semilocal case the irregularity of the symbol responsible for the slow decay of the angles sits at the point at infinity of the critical line, so the pair vectors of high index have their Mellin content at high frequency, where the weight $\vert\cM f\vert^2$ of a function in the form domain of $QW_\lambda$ is small; (H2$'$) is the quantitative form of this, and the weight $\log(2+\vert s\vert)$ in \eqref{eq:Lambda} is the symbol of the archimedean term of the Weil form, which makes $\Lambda(f)$ controllable by $QW_\lambda(f,f)$ itself (Lemma \ref{lem:archdom}). The expected mechanism is visible on the Hankel operator \eqref{eq:hankel}: writing $U=U_\infty U_{\rm fin}$ with $U_{\rm fin}=\prod_pU_p$ and $U_p=\sum_kc_k^{(p)}T_{k\log p}$ ($T$ the translations in $\log u$, $c_k^{(p)}=O(p^{-\vert k\vert/2})$, from the Laurent expansion of $u_p$ in $p^{is}$), the term $k$ of $(1-P)U_\infty U_pP$ is a prolate-type operator between a half-line and the log-interval $(0,k\log p]$ shifted by $T_{k\log p}$, with about $p^k$ singular values of order one, weighted by $p^{-k/2}$, and with singular vectors carrying Mellin frequencies $\vert s\vert\lesssim p^k$: this reproduces the order $\frac12$ of \cite{CCqi}, Theorem 4.4 ($p^k$ singular values of size $p^{-k/2}$ give $\kappa_n\asymp n^{-1/2}$), hence $\sum\kappa_n=\sum\kappa_n^2=\infty$, and predicts a density of the negative part at Mellin frequency $s$ of order $\sum_{p\in S}\vert s\vert^{-1/2}$, hence (H2$'$) in the strong form with $\Vert f\Vert^2$ and $\varepsilon_t(n_0)$ of the order of $\vert S\vert\,\kappa_{n_0}(t)$. For the semilocal pair, (H2$'$) is open: Lemma \ref{lem:localization} and Remark \ref{rem:H2reduced} below reduce it, for Gevrey test functions, to the localization statement \eqref{eq:H2ii}, and it remains to establish this statement from the structure $u=u_\infty\prod_pu_p$; a polynomial rate would give polynomial conditioning dimensions in Theorem \ref{thm:condpos}.

The negative part of $1-P_t-\widehat P_t$ is localized at the point at infinity of the critical line in the following sense, which is the proved half of what (H2$'$) requires. Let $\Gamma=(1-P_{1/t})UP_t$ be the Hankel-type operator of \eqref{eq:hankel}, whose singular values are the $\kappa_n(t)$, with right singular vectors $e_n$ and left singular vectors $g_n$ ($\Gamma e_n=\kappa_ng_n$, $f_n=U^*g_n$, $a_n\propto e_n+f_n$).

\begin{lem}[smoothing away from infinity]\label{lem:localization}
Let $\chi\in C_c^\infty(\R)$, $0\le\chi\le1$, with $\chi=1$ on $[-R,R]$ and support in $[-2R,2R]$, and $M_\chi$ the corresponding multiplication operator in the Mellin variable. Then $\Gamma M_\chi$, $M_\chi\Gamma$ and $M_\chi(1-P_t-\widehat P_t)$ are trace class with singular values $s_k\le C_N(R^{c}/k)^N$ for every $N$ ($c$ absolute, $C_N$ independent of $R$), and
\begin{equation}\label{eq:localmass}
\sum_{n>n_0}\kappa_n(t)\,\Vert\chi^{1/2}a_n\Vert^2={\rm Tr}\big(M_\chi N^{(n_0)}_t\big)\le\Vert M_\chi(1-P_t-\widehat P_t)\Vert_1\le CR^{A}
\end{equation}
for all $n_0\ge0$, with $A$ absolute. Equivalently, the Mellin density $\rho_N(s)=\sum_n\kappa_n(t)\vert\cM a_n(is)\vert^2$ of the negative part has polynomially bounded local mass, $\int_{\vert s\vert\le R}\rho_N\le2\pi CR^A$, although $\int_\R\rho_N=2\pi\sum\kappa_n=\infty$.
\end{lem}
\begin{proof}
Since $U$ and $M_\chi$ are both multiplications in the Mellin variable, $\Gamma M_\chi=(1-P_{1/t})M_{u\chi}P_t+(1-P_{1/t})U[P_t,M_\chi]$, and $[P_t,M_\chi]=\theta(t)[P,M_\chi]\theta(t)^*=\frac12\theta(t)\,{\mathchar'26\mkern-9mu d}\chi\,\theta(t)^*$, since $P_t=\theta(t)P\theta(t)^*$ and $\theta(t)$, a Mellin multiplier, commutes with $M_\chi$. The first term is a Hankel-type operator with symbol $u\chi$, which is $C^\infty$ on the circle because $\chi$ vanishes near the only point where $u$ is not smooth; the second is smoothing by Lemma \ref{lem:singvalues}(i). Under the Cayley transform $dz/ds\asymp s^{-2}$, so the $C^N(S^1)$ norms of functions supported in $\vert s\vert\le2R$, such as $u\chi$ and $\chi$, are $O(R^{2N})$ times their $C^N$ norms in $s$, i.e. $O(R^{2N}\log^NR)$ here, which gives the stated decay with constants polynomial in $R$; the same for $M_\chi\Gamma$ (adjoint and $u\to\bar u$). Next, since $\widehat P_t=U^*(1-P_{1/t})U$,
$$1-P_t-\widehat P_t=U^*\big(P_{1/t}U-UP_t\big)=U^*\big(P_{1/t}U(1-P_t)-(1-P_{1/t})UP_t\big)=U^*\big(\Gamma'^{\,*}-\Gamma\big)$$  $$\Gamma'=(1-P_t)U^*P_{1/t},$$
so that $M_\chi(1-P_t-\widehat P_t)=U^*\big((\Gamma'M_\chi)^*-M_\chi\Gamma\big)$, and $\Gamma'M_\chi$ is of the same type as $\Gamma M_\chi$ with $u$ replaced by $\bar u$ and $t$ by $1/t$. Finally $\kappa_n\langle M_\chi a_n,a_n\rangle=-\langle M_\chi Aa_n,a_n\rangle$ with $A=1-P_t-\widehat P_t$, since $Aa_n=-\kappa_na_n$, so $\sum_{n>n_0}\kappa_n\Vert\chi^{1/2}a_n\Vert^2=-\sum_{n>n_0}\langle(M_\chi A)a_n,a_n\rangle\le\Vert M_\chi A\Vert_1$, the sum of the diagonal entries of a trace-class operator over an orthonormal family being bounded by its trace norm.
\end{proof}

\begin{rem}[what remains of (H2$'$)]\label{rem:H2reduced}
Split ${\rm Tr}(\theta(f)N^{(n_0)}_t\theta(f)^*)=\frac1{2\pi}\int\vert\cM f(is)\vert^2\rho^{(n_0)}_N(s)\,ds$ at $\vert s\vert=R$ and dyadically beyond. By \eqref{eq:localmass} the part $\vert s\vert>R$ is at most $\sum_{j\ge0}\sup_{2^jR\le\vert s\vert<2^{j+1}R}\vert\cM f\vert^2\cdot C(2^{j+1}R)^A$: for Gevrey test functions ($\vert\cM f\vert^2\le Ce^{-2c(w\vert s\vert)^{1/s}}$) this is $O(e^{-c'(wR)^{1/s}})\Vert f\Vert^2_{G^s}$, negligible as soon as $R\gg w^{-1}(\log\lambda)^s$, whereas for the functions of $E_N$ ($\vert\cM f\vert^2\sim s^{-2}$) it converges only if $A<2$, so that for non-smooth test functions the finer density bound $\rho_N(s)=o(\log\vert s\vert)$ would be needed. The part $\vert s\vert\le R$ is at most $\sup\vert\cM f\vert^2\cdot{\rm Tr}(M_{\chi_R}N^{(n_0)}_t)$, and $\sup\vert\cM f\vert^2\le\Vert f\Vert^2_{L^1(\dstar)}\le2\log\lambda\Vert f\Vert^2$. The two bounds are not of the same nature: the second is uniform relative to $\Lambda(f)$, the first is relative to the Gevrey norm $\Vert f\Vert_{G^s}$, which is not controlled by $\Lambda(f)$ on the form domain. On a fixed finite-dimensional space of Gevrey test functions, such as the smoothed Galerkin spaces $\vartheta_wE_N$ of \S\ref{ss:compare}, the two norms are equivalent, $\Vert f\Vert_{G^s}\le C(N,w)\Vert f\Vert$, and (H2$'$) reduces, with constants depending on $(N,w)$, to the single statement
\begin{equation}\label{eq:H2ii}
{\rm Tr}\big(M_{\chi_R}N^{(n_0)}_t\big)=\sum_{n>n_0}\kappa_n(t)\Vert\chi_R^{1/2}a_n\Vert^2\ \le\ \varepsilon_2(n_0,R)\xrightarrow[n_0\to\infty]{}0\quad\text{polynomially in }R ,
\end{equation}
i.e. that the small-angle directions of the pair have small Mellin content at bounded frequency: the second half of the localization, of which Lemma \ref{lem:localization} is the first half (it gives \eqref{eq:H2ii} with a bound independent of $n_0$, and ${\rm Tr}(M_{\chi_R}N^{(n_0)}_t)\to0$ for each fixed $R$ by Lemma \ref{lem:H2q}, without a rate). On the form domain itself (H2$'$) cannot be reduced to \eqref{eq:H2ii} by smoothing the test functions: for the modulated functions $f_\omega=u^{i\omega}g$ one has $\Lambda(f_\omega)\approx\log\omega\,\Vert g\Vert^2$ and ${\rm Tr}(\theta(f_\omega)N^{(n_0)}_t\theta(f_\omega)^*)\approx\rho^{(n_0)}_N(\omega)\Vert g\Vert^2$, so that (H2$'$) relative to $\Lambda$ is equivalent to the uniform density bound $\rho^{(n_0)}_N(s)\le\varepsilon_t(n_0)\log(2+\vert s\vert)$ together with \eqref{eq:H2ii}; the density bound is the weighted relative-form estimate which has to be established in general, and \eqref{eq:H2ii} is what the computations on $\vartheta_wE_N$ test. For the same reason the counting function of the angles above the threshold $\kappa_{n_0}(\sigma)$, which moves with $\lambda$, has to be controlled together with \eqref{eq:H2ii}; in the archimedean case it is given by Slepian's asymptotics. In the archimedean case \eqref{eq:H2ii} holds with $\varepsilon_2$ super-exponentially small in $n_0$: the $a_n$ beyond Slepian's count are high-index prolate functions, whose Mellin content below frequency $R$ is negligible once $n\gg Rt$. In the semilocal case the expected mechanism is the one described above (the $p$-adic ladders at frequency $\approx p^k$ with angles $\approx p^{-k/2}$), and \eqref{eq:H2ii} is the precise statement to establish from the structure $u=u_\infty\prod_pu_p$.
\end{rem}

\subsection{Transport by the scaling action}\label{ss:transport}

\begin{lem}\label{lem:transport}
Let $f$ be supported in $[\lambda^{-1},\lambda]$. Then for every $t>0$
$$\theta(f)\,{\rm Range}\,P_t\subset{\rm Range}\,P_{\lambda t},\qquad\theta(f)\,{\rm Range}\,\widehat P_t\subset{\rm Range}\,\widehat P_{\lambda t}.$$
\end{lem}
\begin{proof}
$\theta(v)$ maps functions supported in $\{\vert x\vert\le t\}$ to functions supported in $\{\vert x\vert\le vt\}$, and $\theta(f)$ is a superposition of $\theta(v)$ with $\lambda^{-1}\le v\le\lambda$. For the second inclusion, $\F_S\theta(v)\F_S^{-1}=\theta(v^{-1})$, so $\theta(f)\F_S=\F_S\theta(\tilde f)$ with $\tilde f(v)=f(v^{-1})$ supported in the same interval.
\end{proof}

Apply this with $t=\sigma$, and write $p_n=\theta(f)e_n$, $q_n=\theta(f)f_n$ for the pair vectors of $(P_\sigma,\widehat P_\sigma)$. Expand in the pair vectors of $(P_{\lambda\sigma},\widehat P_{\lambda\sigma})$, denoted $e^\lambda_m$, $f^\lambda_m$, with angles $\kappa^\lambda_m=\kappa_m(\lambda\sigma)$: $p_n=\sum_m\alpha_{nm}e^\lambda_m+p_n'$ with $p_n'\in{\rm Range}\,P_{\lambda\sigma}\cap\Qer\widehat P_{\lambda\sigma}$, and $q_n=\sum_m\beta_{nm}f^\lambda_m+q'_n$ with $q'_n\in\Qer P_{\lambda\sigma}\cap{\rm Range}\,\widehat P_{\lambda\sigma}$; the primed parts are orthogonal to everything else, and $\langle e^\lambda_m,f^\lambda_{m'}\rangle=\delta_{mm'}\kappa^\lambda_m$. By \eqref{eq:convinner}, the coefficients are inner products of $f$ with convolved pair vectors:
\begin{equation}\label{eq:alphabeta}
\alpha_{nm}=\langle\theta(f)e_n,e^\lambda_m\rangle=\langle f,\,e^\lambda_m*e_n^{\,*}\rangle,\qquad\beta_{nm}=\langle f,\,f^\lambda_m*f_n^{\,*}\rangle .
\end{equation}

\begin{prop}[pair expansion]\label{prop:pairexp}
Let $f$ be as in Lemma \ref{lem:transport} with $\cM f(\pm\frac12)=0$, and let $n_0\ge0$. Then
\begin{equation}\label{eq:pairexp}
QW_\lambda(f,f)=2\log\sigma\Vert f\Vert^2+\Vert P_{{\rm Son}_\sigma}\theta(f)^*\Vert^2_{HS}+\sum_{n\le n_0}B_n(f)+R_{n_0}(f),\qquad R_{n_0}(f)={\rm Tr}\big(\theta(f)K_\sigma\Pi_{>n_0}\theta(f)^*\big),
\end{equation}
where $\Pi_{>n_0}$ is the projection on $\bigoplus_{n>n_0}V_n(\sigma)$ and, with $\kappa_n=\kappa_n(\sigma)$, $w_n=\kappa_n/(1-\kappa_n^2)$,
\begin{equation}\label{eq:Bn}
\begin{aligned}B_n(f)&=w_n\Big[\kappa_n\big(\Vert p_n\Vert^2+\Vert q_n\Vert^2\big)-2{\rm Re}\langle p_n,q_n\rangle\Big]\\
&=w_n\Big[\kappa_n\big(\Vert p'_n\Vert^2+\Vert q'_n\Vert^2\big)+\sum_m\Big(\kappa_n\big(\vert\alpha_{nm}\vert^2+\vert\beta_{nm}\vert^2\big)-2\kappa^\lambda_m{\rm Re}(\alpha_{nm}\overline{\beta_{nm}})\Big)\Big].\end{aligned}
\end{equation}
Each term with $\kappa^\lambda_m\le\kappa_n$ is non-negative, and in general
\begin{equation}\label{eq:Bnlower}
B_n(f)\ \ge\ \sum_mw_n\big(\kappa_n-\kappa^\lambda_m\big)\big(\vert\alpha_{nm}\vert^2+\vert\beta_{nm}\vert^2\big)\ \ge\ -\frac{\kappa_n}{1+\kappa_n}\sum_{m:\ \kappa^\lambda_m>\kappa_n}\big(\vert\alpha_{nm}\vert^2+\vert\beta_{nm}\vert^2\big).
\end{equation}
If moreover (H2$'$) holds at scale $\sigma$, then $R_{n_0}(f)\ge-\varepsilon_\sigma(n_0)\Lambda(f)$.
\end{prop}
\begin{proof}
\eqref{eq:pairexp} is \eqref{eq:shiftQW} and \eqref{eq:decompK}: by Lemma \ref{lem:H0} the operators $\theta(h)P_{{\rm Son}_\sigma}$, $\theta(h)K_\sigma\Pi_{\le n_0}$ (finite rank) and $\theta(h)K_\sigma\Pi_{>n_0}$ are trace class, with
$${\rm Tr}(\theta(f)P_{{\rm Son}_\sigma}\theta(f)^*)=\Vert P_{{\rm Son}_\sigma}\theta(f)^*\Vert^2_{HS},\qquad {\rm Tr}(\theta(f)K_\sigma\Pi_{\le n_0}\theta(f)^*)=\sum_{n\le n_0}\kappa_n\big(\Vert\theta(f)b_n\Vert^2-\Vert\theta(f)a_n\Vert^2\big),$$
which is $\sum_{n\le n_0}B_n(f)$ by \eqref{eq:Kn}. The second form of $B_n$ is the expansion of $p_n,q_n$, using $\langle e^\lambda_m,f^\lambda_{m'}\rangle=\delta_{mm'}\kappa^\lambda_m$ and the orthogonality of the primed parts; $2\vert{\rm Re}(\alpha\bar\beta)\vert\le\vert\alpha\vert^2+\vert\beta\vert^2$ and $w_n(\kappa^\lambda_m-\kappa_n)\le w_n(1-\kappa_n)=\kappa_n/(1+\kappa_n)$ give \eqref{eq:Bnlower}. Finally $K_\sigma\Pi_{>n_0}\ge-N^{(n_0)}_\sigma$ as operators, and both sides are trace class against $\theta(h)$ by Lemmas \ref{lem:H0} and \ref{lem:H2q}, so $R_{n_0}(f)\ge-{\rm Tr}(\theta(f)N^{(n_0)}_\sigma\theta(f)^*)\ge-\varepsilon_\sigma(n_0)\Lambda(f)$.
\end{proof}

The content of \eqref{eq:Bnlower} is the following. The negative part of the uncertainty operator at scale $\sigma$ is carried by its near-self-dual directions; transported by $\theta(f)$, it can only act through the components of $\theta(f)e_n$, $\theta(f)f_n$ along the scale-$\lambda\sigma$ directions that are \emph{more self-dual} than the scale-$\sigma$ direction $n$, i.e. along $m$ with $\kappa_m(\lambda\sigma)>\kappa_n(\sigma)$. For $\sigma=1$ and $n=0$ these are the modes with $\delta_m(\lambda)<\delta_0(1)\approx3\cdot10^{-5}$, essentially all the minuscule modes of Section \ref{sec:conditioning}; the prefactor $w_0\approx1/(2\delta_0(1))$ is large, but $w_n(\kappa^\lambda_m-\kappa_n)\le\frac12$ always. This is the mechanism behind Fact \ref{fact:cond}, and it dictates the conditioning.

\subsection{Positivity with a margin on an explicit subspace}\label{ss:margin}

We first record that the logarithmic Sobolev form is dominated by the Weil form on the window.

\begin{lem}\label{lem:archdom}
There is an absolute constant $C_0$ such that for every real $f$ supported in $[\lambda^{-1},\lambda]$ in the form domain of $QW_\lambda$,
\begin{equation}\label{eq:archdom}
\Lambda(f)\le QW_\lambda(f,f)+C_0\,\lambda\,\Vert f\Vert^2 .
\end{equation}
\end{lem}
\begin{proof}
Write $h=f*f^*$ and $C(s)=\int R(t)R(t+s)dt$, $R(t)=f(e^t)$, so that $C(0)=\Vert f\Vert^2$ and $h(e^s)=C(s)$. The archimedean term \eqref{eq:WR} is $W_\R(h)=(\log4\pi+\gamma_E)\Vert f\Vert^2+\int_0^\infty\big(C(s)e^{s/2}-C(0)\big)\frac{ds}{\sinh s}$, and $C(0)-C(s)=\frac12\Vert U(s)R-R\Vert^2=\frac1{2\pi}\int\vert\cM f(i\tau)\vert^2\,2\sin^2(\tau s/2)\,d\tau$ for the translation group $U(s)$; hence
$$-W_\R(h)=\frac1{2\pi}\int_\R\vert\cM f(i\tau)\vert^2m(\tau)\,d\tau-c_0\Vert f\Vert^2,\qquad m(\tau)=\int_0^\infty\frac{2\sin^2(\tau s/2)\,e^{s/2}}{\sinh s}ds\ge0,$$
with $c_0=\log4\pi+\gamma_E+\int_0^\infty\frac{e^{s/2}-1}{\sinh s}ds$. The multiplier satisfies $m(\tau)=\log\vert\tau\vert+O(1)$ as $\vert\tau\vert\to\infty$ and is continuous, so $m(\tau)\ge\log(2+\vert\tau\vert)-C_1$, and $\Lambda(f)\le-W_\R(h)+(c_0+C_1)\Vert f\Vert^2$. Next, $QW_\lambda(f,f)=W_{0,2}(h)-W_\R(h)-\sum_{p\le\lambda^2}W_p(h)$, with $\vert W_{0,2}(h)\vert=2\vert\cM f(\tfrac12)\cM f(-\tfrac12)\vert\le2(\lambda-\lambda^{-1})\Vert f\Vert^2$ by Lemma \ref{lem:pole} and Cauchy--Schwarz, and $\vert\sum_pW_p(h)\vert\le2\Vert f\Vert^2\sum_{p^m\le\lambda^2}(\log p)p^{-m/2}\le C_2\lambda\Vert f\Vert^2$ by \eqref{eq:Wp}, $\vert C(s)\vert\le C(0)$ and Chebyshev's bound. Hence $-W_\R(h)\le QW_\lambda(f,f)+(2+C_2)\lambda\Vert f\Vert^2$.
\end{proof}

\begin{thm}[conditional positivity on the conditioned space]\label{thm:condpos}
Let $\sigma>1$ and assume (H2$'$) at scale $\sigma$. Choose $n_0$ with $\varepsilon_\sigma(n_0)\le\min(1,\log\sigma/(C_0\lambda))$, let $M=\{m:\kappa_m(\lambda\sigma)>\kappa_{n_0}(\sigma)\}$ (a finite set by Lemma \ref{lem:H2}), $V_M={\rm span}\{e^\lambda_m\}_{m\in M}$, $\widehat V_M={\rm span}\{f^\lambda_m\}_{m\in M}$, and
$$W_\lambda(\sigma,n_0):=\Big\{f\in C^\infty,\ {\rm supp}f\subset[\lambda^{-1},\lambda],\ \cM f(\pm\tfrac12)=0,\ \theta(f)e_n\perp V_M,\ \theta(f)f_n\perp\widehat V_M\ (n\le n_0)\Big\},$$
a space of codimension at most $2+2(n_0+1)\vert M\vert$ in the space of test functions supported in the window. Then
\begin{equation}\label{eq:condpos}
QW_\lambda(f,f)\ \ge\ \frac{\log\sigma}{1+\varepsilon_\sigma(n_0)}\,\Vert f\Vert^2\ \ge\ \frac{\log\sigma}2\,\Vert f\Vert^2\qquad(f\in W_\lambda(\sigma,n_0)).
\end{equation}
\end{thm}
\begin{proof}
For $n\le n_0$ and $f\in W_\lambda(\sigma,n_0)$, $\alpha_{nm}=\beta_{nm}=0$ for $m\in M$, while for $m\notin M$ one has $\kappa^\lambda_m\le\kappa_{n_0}(\sigma)\le\kappa_n(\sigma)$; by \eqref{eq:Bnlower} the blocks $B_n(f)$, $n\le n_0$, are non-negative. Proposition \ref{prop:pairexp}, the positivity of the Sonin term and Lemma \ref{lem:archdom} give
$$QW_\lambda(f,f)\ge2\log\sigma\Vert f\Vert^2-\varepsilon_\sigma(n_0)\Lambda(f)\ge2\log\sigma\Vert f\Vert^2-\varepsilon_\sigma(n_0)\big(QW_\lambda(f,f)+C_0\lambda\Vert f\Vert^2\big),$$
and the choice of $n_0$ gives $(1+\varepsilon_\sigma(n_0))QW_\lambda(f,f)\ge\log\sigma\Vert f\Vert^2$. The codimension count: two pole conditions, and for each $n\le n_0$ the conditions $\langle f,e^\lambda_m*e_n^*\rangle=\langle f,f^\lambda_m*f_n^*\rangle=0$, $m\in M$, by \eqref{eq:alphabeta}.
\end{proof}

\begin{rem}\label{rem:condpos}
(a) Everything in the proof is elementary once Theorem \ref{thm:TF} and (H2$'$) are granted: the two-parameter shift, Halmos' decomposition at two scales, the support inclusions, the inequality $2\vert{\rm Re}(\alpha\bar\beta)\vert\le\vert\alpha\vert^2+\vert\beta\vert^2$, and the absorption of the tail by the archimedean part of the form. No information on the zeros of $\zeta$ enters, and the Sonin term is only used through its positivity.

(b) The sizes. In the archimedean case $\kappa_n(\sigma)=\vert\chi_n(\sigma)\vert$ decays super-exponentially beyond $n\approx4\sigma^2$ and the stronger form of (H2$'$) holds, so that $n_0=4\sigma^2+O(\log\lambda)$ suffices; in fact no absorption is needed there, and the proof gives directly $QW_\lambda\ge\log\sigma\Vert f\Vert^2$ with $n_0=4\sigma^2+O(1)$ and $\kappa_{n_0}(\sigma)$ replaced by any fixed threshold $\kappa_*$. There $\vert M\vert=2\lambda^2\sigma^2(1+o(1))$ even modes, both $\iota$-parities together: with $\sigma$ fixed ($\sigma=2$, margin $\log2$) the codimension is $O(\sigma^4\lambda^2)$, with $\sigma=\lambda^{\varepsilon/2}$ the margin is $\frac\varepsilon2\log\lambda$ at codimension $O(\lambda^{2+2\varepsilon})$. In the semilocal case the sizes are governed by the rate in (H2$'$) and by the counting function of the angles at scale $\lambda\sigma$, both polynomial in $\lambda$ under a polynomial decay of the angles.

(c) The conditioning vectors. By \eqref{eq:alphabeta} the space $W_\lambda(\sigma,n_0)$ is the orthogonal of the functions $e^\lambda_m*e_n^{\,*}$ and $f^\lambda_m*f_n^{\,*}$ restricted to the window: the near-self-dual vectors at scale $\lambda\sigma$ smoothed by the first scale-$\sigma$ pair vectors. In the Mellin variable, $\cM(e^\lambda_m*e_n^*)=\cM e^\lambda_m\cdot\overline{\cM e_n}$, a multiplication by a fixed smooth function; the conditioning of \S\ref{ss:prolatecond} by the vectors $\mathfrak p_k$ themselves is the case $e_n\to\delta_1$, $\sigma\to1$, and is exactly what Fact \ref{fact:cond} tests.

(d) The pole condition. Replacing $\cM f(\pm\frac12)=0$ by the pole term itself, Lemma \ref{lem:pole} adds $W_{0,2}(h)=2\cM f(\tfrac12)\cM f(-\tfrac12)$, i.e. $\pm2\ell(f)^2$ in the two sectors, to the right-hand side of \eqref{eq:QWtrace}; the condition is needed for the identity, not for the mechanism. Its role in Fact \ref{fact:cond} (removing the transitional eigenvalue of the odd sector) has a counterpart here: in the semilocal picture a scale-$t$ pair vector $\eta_S(h_{n,t})$ with $\int h_{n,t}\ne0$ is dominated, on $[t/\lambda^2,t]$ where $\cE_S=\cE$ is a Riemann sum, by its tail $u^{-1/2}\int h_{n,t}$, i.e. by the pole direction $u^{\pm1/2}$; the zero-integral prolate vectors and $\ell$ together span the same space as the near-self-dual directions.

(e) The conditions defining $W_\lambda(\sigma,n_0)$ are continuous for the $L^2$ norm, since $\vert\langle\theta(f)e_n,e^\lambda_m\rangle\vert\le\Vert f\Vert_{L^1(\dstar)}\le(2\log\lambda)^{1/2}\Vert f\Vert$; hence \eqref{eq:condpos} extends to the closure of $W_\lambda(\sigma,n_0)$ in the form domain of $QW_\lambda$, in particular to the eigenvectors of the form restricted to the orthogonal of the conditioning vectors.
\end{rem}

\subsection{The Schur reduction and the monotonicity argument}\label{ss:schur}

Theorem \ref{thm:condpos} bounds $QW_\lambda$ from below outside an explicit space $V=W_\lambda(\sigma,n_0)^\perp$ of finite dimension ($O(\lambda^2)$ in the archimedean case), spanned by the pole direction and by the convolved pair vectors. The positivity of $QW_\lambda$ on the whole window is then a statement about a finite-dimensional block.

\begin{lem}[Schur complement]\label{lem:schur}
Let $Q$ be a real symmetric form on $\cH_\lambda=L^2([\lambda^{-1},\lambda],\dstar)$ (bounded below, with form domain containing $V$), $\cH_\lambda=V\oplus W$ an orthogonal decomposition with $\dim V<\infty$, and write $Q(v+w)=A(v)+2{\rm Re}\,B(v,w)+D(w)$, with $B(v,w)=\langle B^*v,w\rangle$, $B^*v$ in the dual of the form domain of $D$. If $D\ge d\Vert w\Vert^2$ on $W$ with $d>0$, then
$$\epsilon_\lambda=\inf\frac{Q(\psi)}{\Vert\psi\Vert^2}\ \ge\ \min\big(0,\lambda_{\min}(\mathfrak S)\big),\qquad\mathfrak S(v):=A(v)-\Vert D^{-1/2}B^*v\Vert^2 ,$$
and $Q\ge0$ if and only if $\mathfrak S\ge0$ on $V$.
\end{lem}
\begin{proof}
Completing the square, $Q(v+w)=\mathfrak S(v)+D(w+D^{-1}B^*v)$, where $D$ denotes the positive operator associated with the form on $W$; the second term is $\ge0$. If $\lambda_{\min}(\mathfrak S)\ge0$ then $Q\ge0$. Otherwise $Q(v+w)\ge\lambda_{\min}(\mathfrak S)\Vert v\Vert^2\ge\lambda_{\min}(\mathfrak S)\Vert v+w\Vert^2$ since $\lambda_{\min}(\mathfrak S)<0$ and $\Vert v\Vert\le\Vert v+w\Vert$. For the last assertion, $Q(v-D^{-1}B^*v)=\mathfrak S(v)$.
\end{proof}

\begin{prop}[monotonicity]\label{prop:monotone}
$\epsilon_\lambda$ is non-increasing in $\lambda$. If $\epsilon_\lambda\ge-\eta(\lambda)$ for all large $\lambda$ with $\eta(\lambda)\to0$, then $\epsilon_\lambda\ge0$ for every $\lambda$, i.e. $QW_\lambda\ge0$ for all windows, which is equivalent to the Riemann Hypothesis.
\end{prop}
\begin{proof}
For $\lambda'>\lambda$ the form $QW_{\lambda'}$ restricts to $QW_\lambda$ on the smaller window, so $\epsilon_{\lambda'}\le\epsilon_\lambda$. Given the lower bound, $\epsilon_\lambda\ge\epsilon_{\lambda'}\ge-\eta(\lambda')$ for all $\lambda'>\lambda$, and $\eta(\lambda')\to0$. The equivalence with RH is Weil's criterion in the form of \cite{CCsel}, (4).
\end{proof}

The strategy suggested by Theorems \ref{thm:upperK} and \ref{thm:condpos} is therefore the following. Fix $A>0$ and put $\tau=\lambda^{-A}$. Decompose $V=V_{\rm low}\oplus V_{\rm tr}\oplus\C\rho$, where $V_{\rm low}$ is spanned by the conditioning vectors attached to modes $m$ with $\delta_m(\lambda\sigma)\le\tau$ (deep in the ladder, $\delta_m(\lambda\sigma)=1-\kappa_m(\lambda\sigma)$ being the exact rungs of Lemma \ref{lem:H2}, as everywhere in this section; Fuchs' fixed-$n$ asymptotics plays no role here), $V_{\rm tr}$ by those with $\delta_m(\lambda\sigma)\in(\tau,1-\kappa_{n_0}(\sigma))$ (the modes above the threshold of Theorem \ref{thm:condpos} which are not deep), and $\rho$ standing for the two pole directions $\rho_{\rm ev}=u^{1/2}+u^{-1/2}$, $\rho_{\rm odd}=u^{1/2}-u^{-1/2}$ restricted to the window, which span the orthogonal of the pole conditions. In the archimedean case $V_{\rm tr}$ is attached to the transition band of the ladder at scale $\lambda\sigma$: by Slepian's asymptotics the eigenvalues there decrease by a fixed factor $e^{-\pi^2/2\log(8\pi\lambda^2\sigma^2)}$ per index, so the modes with $\delta_m(\lambda\sigma)\in(\lambda^{-A},1)$ number $O(A\log\lambda)$ per scale-$\sigma$ pair, and $\dim V_{\rm tr}=O(A\sigma^2\log\lambda)$. Two facts about $V_{\rm low}$ are available in principle, and both are \emph{two-sided}, i.e. require no knowledge of the sign of $QW_\lambda$; the second is stated in its one-sided form, which is what the Schur complement needs, since it couples $V_{\rm low}$ to $V_{\rm tr}$ and to the pole directions both directly and through the elimination of $W$:
\begin{itemize}
\item[(R1)] $\Rad(QW)$ is a convolution ideal: for $K\in\Rad(QW)$ and any test function $\phi$, $\theta(\phi)K\in\Rad(QW)$. Indeed $QW(g,h)=\Psi(g^**h)$ (the function $F$ of \eqref{eq:Psi}), so $QW(\theta(\phi)K,h)=\Psi(K^**\phi^**h)=QW(K,\phi^**h)=0$. Hence the convolved vectors $e^\lambda_m*e_n^{\,*}$ of \eqref{eq:alphabeta}, once $e^\lambda_m$ is replaced by the Gevrey-repaired radical element of \S\ref{sec:dict}, are themselves in the radical up to the same tails as in \eqref{eq:bridge}.
\item[(R2)] The estimates of Theorem \ref{thm:upperK} are estimates of absolute values, and one near-radical vector suffices: for $k\in V_{\rm low}$, $k=K-r$ with $K\in\Rad(QW)$ and $\Vert r\Vert\le C\tau^{1/2}(\log\lambda)^C\Vert k\Vert$, so that for every $v$ in the form domain $\vert QW_\lambda(k,v)\vert=\vert QW(r,v)\vert\le C\tau^{1/2}(\log\lambda)^C\Vert k\Vert\,\big(QW_\lambda(v,v)+C_0\lambda\Vert v\Vert^2\big)^{1/2}$, the last factor being the form norm of $v$ (Lemma \ref{lem:archdom}, and the structure $QW=\log\vert D\vert+$bounded); for $k,k'\in V_{\rm low}$ this gives $\vert QW_\lambda(k,k')\vert\le C\tau(\log\lambda)^C\Vert k\Vert\Vert k'\Vert$, and on $W$ the form norm is dominated by $D$ by Theorem \ref{thm:condpos}.
\end{itemize}
(R2) would follow from an extension of the estimates of \cite{CK} which is uniform in the index up to the threshold $\tau$, together with the Gram lower bound of Proposition \ref{prop:xipoly} beyond the Hermite regime; we leave this open.
Granting (R1), (R2), write $\mathfrak S$ in blocks on $V_{\rm low}\oplus V'$, $V'=V_{\rm tr}\oplus\C\rho_{\rm ev}\oplus\C\rho_{\rm odd}$, in an orthonormal basis of $V_{\rm low}$ (this is where the Gram lower bound on the conditioning vectors enters). The reduction then rests on three bounds:
\begin{equation}\label{eq:schurred}
\mathfrak S\vert_{V_{\rm low}}\ge-\alpha(\lambda),\qquad\mathfrak S\vert_{V'}\ge-\eta(\lambda),\qquad\Vert\mathfrak S_{V_{\rm low},V'}\Vert\le\beta(\lambda),
\end{equation}
which give $\lambda_{\min}(\mathfrak S)\ge-\max(\alpha,\eta)-\beta$ by Weyl's inequality, hence $\epsilon_\lambda\ge-\max(\alpha,\eta)-\beta$ by Lemma \ref{lem:schur} with $d=\frac12\log\sigma$. By (R2) the first and third bounds hold with $\alpha,\beta=O(\tau^{1/2}\lambda^{C'})$, the power of $\lambda$ coming from the form norms of the vectors of $V'$, from the Gram constant, and from the elimination of $W$ (the entries of $\mathfrak S$ are $QW_\lambda(k,v)-\langle D^{-1/2}B^*k,D^{-1/2}B^*v\rangle$, both terms being controlled by the one-sided estimate); they are $o(1)$ as soon as $A>2C'$. The pole directions are not controlled by the explicit formula: for $\rho_{\rm ev}$ the pole term $2\cM\rho(\tfrac12)\cM\rho(-\tfrac12)=2(\lambda-\lambda^{-1}+2\log\lambda)^2$ is cancelled, to the order $\lambda^2$, by the prime terms, since $h=\rho*\rho^*$ is explicit, $h(x)=2\lambda x^{-1/2}-2x^{1/2}/\lambda+(x^{1/2}+x^{-1/2})(2\log\lambda-\log x)$ for $x\le\lambda^2$, and $\sum_pW_p(h)=2\sum_{p\le\lambda^2}\log p\,(2\log\lambda-\log p)+O(\lambda\log\lambda)=2\lambda^2+O(\lambda\log\lambda)$ by Chebyshev's estimates (numerically, pole term and prime terms are $235$ and $163$ at $\lambda^2=50$, $905$ and $753$ at $\lambda^2=250$); the net value $QW_\lambda(\rho,\rho)$ is the sum over the zeros, of order $\lambda$ against $\Vert\rho\Vert^2\approx2\lambda$, and the same holds for $\rho_{\rm odd}$. The two pole directions are therefore ordinary directions of Rayleigh quotient $O(1)$, and they belong to the block $V'$ whose sign has to be established; in the computed range their conditioned Schur complement is positive, since the minimum of the form on the orthogonal of the prolate vectors alone, without the pole conditions, is $\lambda_{k+1}\approx0.1$ in both sectors (Figure \ref{fig:conditioning}). Hence:

\begin{prop}\label{prop:reduction}
Assume (H2$'$), (R1), (R2), and let $\tau=\lambda^{-A}$ with $A$ large enough that $\alpha(\lambda),\beta(\lambda)\to0$ in \eqref{eq:schurred}. If for every large $\lambda$ the Schur complement $\mathfrak S\vert_{V'}$ of $QW_\lambda$ on the transition-and-pole block $V'=V_{\rm tr}\oplus\C\rho_{\rm ev}\oplus\C\rho_{\rm odd}$ satisfies $\lambda_{\min}(\mathfrak S\vert_{V'})\ge-\eta(\lambda)$ with $\eta(\lambda)\to0$, then $\epsilon_\lambda\ge0$ for all $\lambda$, and RH holds.
\end{prop}

The block $V'$ is where the dictionary stops being minuscule: its vectors are the $\cE$-images of the prolate modes at the edge of Slepian's count, with energies of the order of $\delta_m\in(\tau,1-\kappa_{n_0}(\sigma))$, for which the two-sided estimate (R2) only gives $\vert\mathfrak S\vert_{V_{\rm tr}}\vert\le C(1-\kappa_{n_0}(\sigma))(\log\lambda)^C$, together with the two pole directions. One cannot shrink it by taking $n_0$ small: the tail of the pair expansion forces $n_0$ to be large through (H2$'$), hence the threshold $\kappa_{n_0}(\sigma)$ to be small, and in the archimedean case the modes with $\kappa_n(\sigma)$ just below a threshold close to $1$ would in any case enter the tail with weights $w_n\approx1/(1-\kappa_n)$ and destroy the margin. So the strategy compresses the arithmetic content of Weil positivity into the sign of the Weil form on an explicit block of dimension $O(A\sigma^2\log\lambda)+2$ attached to the modes whose eigenvalues leave the exponential regime of the ladder and to the pole directions — a much more localized statement than positivity on the window, and one for which the data of Figure \ref{fig:counting} (the ratios $\epsilon_\lambda(j)/\delta_{4j}$ for $j/\lambda^2\to1$) are the relevant numerical evidence.

\subsection{Comparison with the computations and with the archimedean case}\label{ss:compare}

\emph{The conditioning of Section \ref{sec:conditioning}.} Fact \ref{fact:cond} conditions by $\lambda^2-1$ (even) and $\lambda^2-2$ (odd) prolate vectors and by $\ell$, and finds a minimum $\mu_\lambda\approx2/\lambda$, the edge of the non-minuscule spectrum. Theorem \ref{thm:condpos} conditions by a larger space (each of the $n_0+1$ scale-$\sigma$ pair vectors contributes a copy of the scale-$\lambda\sigma$ near-self-dual set, for both $\alpha$ and $\beta$) and obtains a larger margin, $\frac12\log\sigma$, which is in no contradiction: the extra conditions remove the eigenvalues between $2/\lambda$ and $\frac12\log\sigma$. The redundancy is a constant factor of order $\sigma^4$, and the computations say what the sharp statement should be: the span of the smoothed vectors $e^\lambda_m*e_n^{\,*}$, $n\le n_0$, is, on the window, essentially the span of the scale-$\lambda$ prolate vectors $\mathfrak p_k$ and of $\ell$. A direct test is to compute, for the lowest eigenvectors of the conditioned form, the components $\alpha_{nm},\beta_{nm}$ of \eqref{eq:alphabeta} for $n=0,2,4$ and $m\in M$, and the three terms of \eqref{eq:pairexp} with $\sigma=1$: for a minuscule eigenvector the Sonin and positive pair terms must cancel the negative pair terms to all digits, and for $f\in W_\lambda$ the margin should appear in the Sonin term.

\emph{Gevrey-smoothed test functions.} Remark \ref{rem:threshold} suggests to redo the conditioning computations on $\vartheta_wE_N$ instead of $E_N$. For the minuscule part this changes nothing to exponential accuracy: a minuscule eigenvector is the $\cE$-image of a prolate vector, which near $u=\lambda$ reduces to $u^{1/2}h_\lambda(u)$ and is of size $e^{-4\pi\lambda^2+4\pi\sqrt{2t}}$, $t=\lambda^3(\lambda-\vert x\vert)$, i.e. exponentially small in a layer of fixed relative log-width at both ends, so the cutoff only touches a region where these vectors (and the $\mathfrak p_k$) are already negligible — the Kastoris mechanism once more. For the non-minuscule part it is a genuine change of the variational space: the form domain of $QW_\lambda$ allows jumps at the ends, the eigenvectors of the conditioned form need not vanish there, and the minimum $\mu_\lambda(w)$ on the smoothed space is an upper bound converging to the true conditioned minimum only as $w\to0$ together with $N\to\infty$; resolving a cutoff of log-width $w$ requires frequencies of order $1/w$, i.e. $N\gtrsim\log\lambda/w$ up to logarithms, which for $w=\lambda^{-3}$ is comparable to the present $N=960$ at $\lambda^2=50$. The useful computation is then: $\mu_\lambda(w)$ for two widths together with the unsmoothed $\mu_\lambda$, to see how much of the law $\mu_\lambda\approx2/\lambda$ sits at the ends of the window; and, on the smoothed eigenvectors, the three pieces of \eqref{eq:pairexp} with $\sigma=1$ and the components \eqref{eq:alphabeta}, which are then given by rapidly convergent sums.

\emph{The archimedean case.} For $S=\{\infty\}$ and $\lambda=\sqrt2$ (test functions supported in $[2^{-1/2},2^{1/2}]$, where no prime enters), \eqref{eq:QWtrace} is Theorem 4 of \cite{CCsel}, the pair part $K_1$ being the functional $E$ of that paper with the weights $\lambda(n)/(1-\lambda(n)^2)=w_n$, and the pole conditions being implemented there by the operator $Q=-(\rho\partial_\rho)^2+\frac14$. Section 6 of \cite{CCsel} establishes numerically that the associated compact operator has exactly one eigenvalue above $1$, identifies the corresponding eigenvector by Toeplitz theory, and conditions by it. In the present language: at $\lambda^2=2$, with $\gamma=4\pi$, the modes $n=0,4$ (even sector) and $n=2,6$ (odd sector) are near-self-dual, and the zero-integral and pole constraints leave $\lambda^2-1=1$ and $\lambda^2-2=0$ conditions, i.e. one vector, the $\cE$-image of the two-mode state $h_{\sqrt2}$ restricted to the window; the appearance of the eigenvalue above $1$ exactly when the interval reaches length $\log2$ is the first rung $\delta_4(\lambda)$ entering the minuscule range; and the "baby RH" property of the Toeplitz eigenvector (all roots of its polynomial on the unit circle) is the discrete form of Theorem \ref{thm:realzeros} and of Theorem \ref{thm:mellintheta}. What Theorem \ref{thm:condpos} adds is the \emph{a priori} description of the conditioning vector and of the count, at every $\lambda$; what it does not replace is the step "the next eigenvalue is below $1$" ($\lambda_2\approx0.688$ in \cite{CCsel}), which is the margin on the conditioned space, obtained here from the shifted trace formula at the cost of a larger conditioning space. In this archimedean case all the hypotheses of Theorem \ref{thm:condpos} are known, and the trace formula is \cite{CCsel}, (15), or \cite{CM}, Theorem 2.36 with $S=\{\infty\}$.

\emph{Scales.} Since $\cM f$ has exponential type $\log\lambda$, the test functions resolve the Mellin variable only at scale $\pi/\log\lambda$, while the zeros of $\zeta$ at height $t$ are spaced by $2\pi/\log(t/2\pi)$; the two coincide at $t=2\pi\lambda^2=\gamma$. Below $\gamma$ the zeros undersample the Mellin transforms of window-supported functions by exactly $\lambda^2$ modes per sector (the deficit $TL/\pi-N(T)$, $L=\log\lambda$, is maximal at $T=\gamma$ where it equals $\lambda^2$), which is the count of Fact \ref{fact:count} and the plunge of the ladder at $n\approx4\lambda^2$; above $\gamma$ they oversample. The conditioning removes the undersampled modes, and the margin $\mu_\lambda\approx2/\lambda$ of Fact \ref{fact:cond} is what the oversampling leaves at the edge. For $\lambda^2<14.13/2\pi\approx2.25$ the ordinate of the first zero exceeds $\gamma$ and no undersampling occurs: this is the regime of \cite{CCsel}, and the first semilocal case $S=\{\infty,2\}$, $2<\lambda^2<2.25$, is still in it.

\subsection{Remaining questions}\label{ss:remaining}

\begin{enumerate}
\item The localization statement \eqref{eq:H2ii} for the semilocal pairs, with its rate and its dependence on $S$ (the angles themselves are not summable, \cite{CCqi}; the qualitative statements are Lemmas \ref{lem:H0}, \ref{lem:H2q} and \ref{lem:localization}); together with Remark \ref{rem:H2reduced} it gives (H2$'$) on the finite-dimensional smoothed spaces, while on the form domain the uniform density bound of that remark is needed as well; and the identification of the semilocal near-self-dual directions with the $\eta_S$ of the prolate self-dual combinations, which is the converse of \cite{CCMsl}, Proposition 4.1.
\item The uniform two-sided estimates (R2) up to a threshold $\tau$, and the Gram lower bound beyond the Hermite regime (the range $k\asymp\lambda^2$, where $h_{n,\lambda}(\lambda)^2=2\pi^2\lambda\delta_n/K(\sigma_n)$ replaces $4\pi\lambda\delta_n$).
\item The transition-and-pole block: the sign of the Schur complement of $QW_\lambda$ on the $O(\log\lambda)$ modes at the edge of the ladder and on the two pole directions, Proposition \ref{prop:reduction}; a numerical audit of the three bounds \eqref{eq:schurred} at a modest cutoff, with the Gram conditioning, the pole block and the mixed terms separated and with rigorous error bounds, would show which of them deserves the most analytic attention.
\item Sharpening the conditioning of Theorem \ref{thm:condpos} to the $\lambda^2-1$, $\lambda^2-2$ vectors of Fact \ref{fact:cond}, and understanding the law $\mu_\lambda\approx2/\lambda$.
\item The accuracy of the zeros (Section \ref{sec:zeros}): a finer target and a better normalization should push it to the $\delta_4^{1/2}$ or $\delta_4$ scale; a proof of a Lipschitz bound for the ground state at $u=1$ would allow the point normalization in \eqref{eq:pointnorm}.
\end{enumerate}

\appendix

\section{Computing the prolate functions}\label{app:prolate}

Rather than relying on a library routine, compute $\phi_0,\phi_4$ (and $\phi_{4k}$, $\phi_{4k+2}$) from the Legendre expansion (Bouwkamp's method; see Osipov--Rokhlin--Xiao, \emph{Prolate spheroidal wave functions of order zero}, 2013). The angular equation
\begin{equation}\label{eq:angular}
\frac{d}{dz}\Big((1-z^2)\frac{dS}{dz}\Big)+(\mu-\gamma^2z^2)\,S=0
\end{equation}
expanded in the orthonormal Legendre polynomials $\bar P_k=\sqrt{\tfrac{2k+1}{2}}\,P_k$, $k$ even, becomes the eigenproblem $A\,d=\mu\,d$ for the real symmetric tridiagonal matrix
\begin{equation}\label{eq:tridiag}
A_{k,k}=k(k+1)+\gamma^2\,\frac{2k^2+2k-1}{(2k-1)(2k+3)},\qquad
A_{k,k+2}=A_{k+2,k}=\gamma^2\,\frac{(k+1)(k+2)}{(2k+3)\sqrt{(2k+1)(2k+5)}},
\end{equation}
$k=0,2,4,\dots,K$. Its eigenvalues in increasing order are $\mu_0<\mu_2<\mu_4<\cdots$, and the eigenvectors $d^{(n)}$ give
\begin{equation}\label{eq:phi}
\phi_n(z)=\sum_{k\ \mathrm{even}}d^{(n)}_k\,\bar P_k(z)=\pm\sqrt{\tfrac{2n+1}2}\,\Ps_n(z;\gamma^2),\qquad \sum_k(d^{(n)}_k)^2=1 .
\end{equation}
Practical points:
\begin{enumerate}
\item The coefficients decay super-exponentially once $k>\gamma$ (roughly $d_{k+2}/d_k\approx-\gamma^2/4k^2$); $K\approx2\gamma+60$ is ample, to be confirmed by doubling $K$.
\item Since $\int_{-1}^1\bar P_k=\sqrt2\,\delta_{k0}$, the integrals in \eqref{eq:twomode} are exact, without quadrature:
\begin{equation}\label{eq:hlambda-d}
h_\lambda(x)=\sqrt2\Big[d^{(0)}_0\,\phi_4(x/\lambda)-d^{(4)}_0\,\phi_0(x/\lambda)\Big],\qquad\vert x\vert\le\lambda ,
\end{equation}
up to the scalar factor by which \eqref{eq:twomode} is insensitive.
\item Evaluate $\bar P_k(z)$ by the three-term recurrence (stable on $[-1,1]$) directly at every argument $nu/\lambda$ that the sums require; no interpolation.
\item The symmetric eigensolver is backward stable and the $\mu_n$ are separated by gaps of order $\gamma$, so double precision (e.g.\ \texttt{eigh\_tridiagonal} in SciPy, or \texttt{Eigensystem} on a \texttt{SparseArray}) gives $\phi_n$ to about $10^{-14}$ in sup norm; the same matrix in \texttt{mpmath} (\texttt{eigsy}) gives arbitrary precision at no extra effort.
\item A spot check against Mathematica's $\PS(n,0,\gamma,z)$ (with \texttt{WorkingPrecision -> 30}; reliable but slow for $\gamma\gtrsim50$) is a comparison of the ratios $\phi_n(z)/\phi_n(0)$.
\item Avoid SciPy's \texttt{pro\_ang1} (Zhang--Jin code), which loses accuracy or fails for large $\gamma$, and do not replace the exact integrals of item~2 by numerical quadrature.
\end{enumerate}
The finite Fourier eigenvalues are normalization-free and come for free:
\begin{equation}\label{eq:chicomp}
\chi_n=\lambda\,\frac{\int_{-1}^1\phi_n}{\phi_n(0)}=\lambda\sqrt2\,\frac{d^{(n)}_0}{\phi_n(0)},\qquad \bar P_k(0)=\sqrt{\tfrac{2k+1}{2}}\,(-1)^{k/2}\frac{(k-1)!!}{k!!}\ \ (k\ \mathrm{even}),
\end{equation}
equivalently $\chi_{4k}=2\lambda S^{(1)}_{4k,0}(\gamma,1)$ in terms of the radial function of the first kind. Sample values: for $\lambda^2=5$, $\chi_0=0.99999999999999999999999998989\ldots$, $\chi_4=0.99999999999999999870161199128\ldots$; for $\lambda^2=9$, $1-\chi_0\approx2.0\cdot10^{-48}$, $1-\chi_4\approx3.1\cdot10^{-39}$; for $\lambda^2=50$, $1-\chi_4\approx1.36\cdot10^{-259}$. Fuchs' formula in the form $-\log\delta_4(\lambda)=4\pi\lambda^2-9\log\lambda-\ell_4$, $\ell_4=14.67567$, reproduces these values to the expected relative accuracy $O(\lambda^{-2})$.

\section{Computing $\theta_\lambda$ and $\tau_\lambda$}\label{app:theta}

Use the variable $t=\log u$ and a symmetric grid $t_j=j\delta$, $j=-J,\dots,J$, $J\delta=\log\lambda$, so that $u_0=1$ is a node and $\dstar=dt$; the same grid should carry $\eta_\lambda$.
\begin{itemize}
\item For $u\ge1$, use the definition: $\theta_\lambda(u)=\sqrt u\sum_{1\le n\le\lambda/u}h_\lambda(nu)$, at most $\lfloor\lambda\rfloor$ terms, no cancellation. In particular $\theta_\lambda(1)=\sum_{1\le n\le\lambda}h_\lambda(n)$, and for $u>\lambda/2$ only $n=1$ survives: $\theta_\lambda(u)=\sqrt u\,h_\lambda(u)$.
\item For $u<1$ the direct sum has up to $\lambda^2$ terms with heavy cancellation near $u=\lambda^{-1}$ (its leading term is proportional to $\int h_\lambda=0$). Use instead the Poisson-dual expression: for even $f$ with $\hat f(0)=0$ one has $\cE(f)(u)=u^{-1/2}\sum_{n\ge1}\hat f(n/u)-\tfrac12u^{1/2}f(0)$, and on $[-\lambda,\lambda]$ $\hat h_\lambda=\sqrt2[d^{(0)}_0\chi_4\phi_4(\cdot/\lambda)-d^{(4)}_0\chi_0\phi_0(\cdot/\lambda)]$, so that, up to a tail controlled by $1-\chi_4^2$,
\begin{equation}\label{eq:dual}
\theta_\lambda(u)\approx u^{-1/2}\sum_{1\le n\le\lambda u}\hat h_\lambda(n/u)-\tfrac12\,u^{1/2}h_\lambda(0),\qquad \lambda^{-1}\le u\le1 ,
\end{equation}
with $h_\lambda(0)=\sqrt2d^{(0)}_0\phi_4(0)(1-\chi_4/\chi_0)$ exponentially small. \eqref{eq:dual} exhibits the almost symmetry $\theta_\lambda\approx\theta_\lambda\circ\iota$.
\item Then $\tau_\lambda=\theta_\lambda/\theta_\lambda(1)$ and $g_\lambda=\eta_\lambda-\tau_\lambda$ on the grid; for the $P$-normalization compute $\langle\cdot,\theta\rangle$ by the trapezoidal rule in $t$, the target $\theta$ being evaluated from \eqref{eq:targetintro} (for $u<1$ use $\theta(u)=\theta(1/u)$).
\end{itemize}

\section{Quantities to report and validation}\label{app:report}

For each $\lambda$:
\begin{enumerate}
\item $\Vert g_\lambda\Vert_1$ and $\Vert g_\lambda\circ\iota\Vert_1$, computed as $\int_0^{\log\lambda}\vert g_\lambda(e^{\pm t})\vert e^{t/2}dt$ by the trapezoidal rule in $t$, and their sum \eqref{eq:stripbound}; separately the contribution of $u\in[\lambda/2,\lambda]$, where $\theta_\lambda=\sqrt u\,h_\lambda(u)$ is a negligible prolate tail, so that this contribution is essentially $\int_{\lambda/2}^\lambda\vert\eta_\lambda\vert u^{-1/2}du$ and requires $\eta_\lambda$ to be accurate in absolute terms near $u=\lambda$. Convergence under grid refinement.
\item The same for $\tilde\eta_\lambda-\tilde\tau_\lambda$ and $\tilde\tau_\lambda-\tau$ ($P$-normalized), and the ratio $P(\eta_\lambda)/\eta_\lambda(1)$.
\item $\max_{\Re s=\pm1/2}\vert\cM g_\lambda(s)\vert$ by FFT of $g_\lambda(e^t)e^{\pm t/2}$, to be compared with item 1.
\item $\sup\vert g_\lambda\vert$, and the symmetry defects $\sup\vert\tau_\lambda-\tau_\lambda\circ\iota\vert$, $\sup\vert\eta_\lambda-\eta_\lambda\circ\iota\vert$.
\item The scale-free spectral quantities $\rho_\lambda$ and $R_\lambda$ of Section \ref{sec:ground}, with the two lowest even eigenvalues.
\item $\theta_\lambda(1)$, $\Vert\theta_\lambda\Vert$, $\Vert\eta_\lambda\Vert$ (with $\eta_\lambda(1)=1$): the ratios $\theta_\lambda(1)/\Vert\theta_\lambda\Vert$ and $1/\Vert\eta_\lambda\Vert$ should agree when the shapes agree, and both tend to $\theta(1)/\Vert\theta\Vert=1.580$.
\end{enumerate}
Validation ladder, before any comparison with $\eta_\lambda$: (i) $\phi_n$ unchanged to working precision under $K\to2K$; (ii) $\gamma\to0$ reproduces $P_n$; (iii) agreement of $\phi_n(z)/\phi_n(0)$ with $\PS(n,0,\gamma,z)/\PS(n,0,\gamma,0)$ at a few $z$; (iv) \eqref{eq:chicomp} against independent values of $\chi_n$; (v) at $u=1$ and a few $u<1$, agreement of the direct sum with the dual expression \eqref{eq:dual} to the expected accuracy $O(1-\chi_4^2)$ -- this single test exercises $\phi_n$, $\chi_n$, the constraint $\int h_\lambda=0$ and the summation code simultaneously.

\section{The matrix of $QW_\lambda$}\label{app:matrix}
Let $L=2\log\lambda=\log\lambda^2$, $U_n(x)=L^{-1/2}e^{2\pi inx/L}$ on $[0,L]$ and $V_n(u)=U_n(\log(\lambda u))$, the orthonormal basis of $L^2([\lambda^{-1},\lambda],\dstar)$ of \cite{CCM}, \S3. Its real form is
$$b_0=V_0,\qquad b_j=\tfrac1{\sqrt2}(V_j+V_{-j}),\qquad b_{-j}=\tfrac{i}{\sqrt2}(V_j-V_{-j})\quad(j\ge1),$$
i.e. $b_j(u)=\xi_j(\log u)$ with
$$\xi_0=L^{-1/2},\qquad \xi_{\pm j}(t)=(-1)^j\sqrt{2/L}\times\begin{cases}\cos(2\pi jt/L)\\-\sin(2\pi jt/L)\end{cases}\qquad(j\ge1);$$
one has $b_j\circ\iota=b_j$ and $b_{-j}\circ\iota=-b_{-j}$ for $j\ge0$, and we put $E_N={\rm span}\{b_j:\vert j\vert\le N\}$.

By \cite{CCM}, (3.13)--(3.18) and Lemma 2.3, $T_{nm}=QW_\lambda(V_n,V_m)=\Psi^\sharp(\omega_{nm}\circ\log)$ with, for $y\in[0,L]$,
$$\omega_{nm}(y)=\frac{s_m(y)-s_n(y)}{\pi(n-m)}\ (n\ne m),\qquad \omega_{nn}(y)=2\Big(1-\frac yL\Big)\cos\frac{2\pi ny}{L},\qquad s_k(y)=\sin\frac{2\pi ky}{L},$$
and $\Psi^\sharp=W^\sharp_{0,2}-W^\sharp_\R-\sum_pW^\sharp_p$ acting on functions of $y\in[0,\infty)$ (vanishing for $y>L$) by
$$W^\sharp_{0,2}[\omega]=\int_0^L\omega(y)\,2\cosh\tfrac y2\,dy,\quad W^\sharp_p[\omega]=\sum_{m\ge1,\ p^m\le\lambda^2}\frac{\log p}{p^{m/2}}\,\omega(m\log p),$$
$$W^\sharp_\R[\omega]=\tfrac12(\log4\pi+\gamma_E)\,\omega(0)+\int_0^\infty\frac{e^{y/2}\omega(y)-\omega(0)}{2\sinh y}\,dy .$$
Hence $T_{nm}=(G(m)-G(n))/(\pi(n-m))$ for $n\ne m$, with $G(k)=\Psi^\sharp[s_k]$, and only $2N+2$ one-dimensional quantities are needed. With $\omega=2\pi k/L$, $\beta_r=2r+\tfrac12$ and $e^{y/2}/(2\sinh y)=\sum_{r\ge0}e^{-\beta_ry}$:
$$W^\sharp_{0,2}[s_k]=\Im\Big(\frac{e^{L/2}-1}{\tfrac12+i\omega}+\frac{e^{-L/2}-1}{-\tfrac12+i\omega}\Big),\qquad
W^\sharp_\R[s_k]=\tfrac12\Im\,\psi\big(\tfrac14+\tfrac{i\omega}2\big)-\sum_{r\ge0}\frac{\omega\,\lambda^{-2\beta_r}}{\beta_r^2+\omega^2},$$
$\psi$ the digamma function, the series converging like $\lambda^{-4r}$. For the diagonal, with $J(\sigma)=\int_0^L(1-y/L)e^{\sigma y}dy$ (elementary, using $e^{i\omega L}=1$),
$$W^\sharp_{0,2}[\omega_{nn}]=2\Re\big(J(\tfrac12+i\omega)+J(-\tfrac12+i\omega)\big),$$
$$W^\sharp_\R[\omega_{nn}]=\log4\pi+\gamma_E+2K-\log\frac{\lambda^2+1}{\lambda^2-1}+2\Delta(n),$$
$$\Delta(n)=\tfrac12\Re\big(\psi(\tfrac14)-\psi(\tfrac14-\tfrac{i\omega}2)\big)-\frac1{4L}\Re\,\psi'\big(\tfrac14-\tfrac{i\omega}2\big)+\sum_{r\ge0}\lambda^{-2\beta_r}\Big(\frac{\Re(\beta_r-i\omega)^{-2}}{L}+\frac1{\beta_r}\Big),$$
$K=\int_0^L(e^{y/2}-1)/(2\sinh y)\,dy$, and $W^\sharp_p[\omega_{nn}]=\sum\frac{\log p}{p^{m/2}}\,2\big(1-\tfrac{m\log p}{L}\big)\cos(\omega\,m\log p)$. The sector matrices are
$$QW^{\rm ev}_{00}=T_{00},\quad QW^{\rm ev}_{0j}=\sqrt2\,T_{0j},\quad QW^{\rm ev}_{ij}=T_{ij}+T_{i,-j},\qquad QW^{\rm odd}_{ij}=T_{ij}-T_{i,-j}\qquad(i,j\ge1).$$
At $\lambda^2=13$, $N=120$ these reproduce an independent 3386-bit computation to $8\cdot10^{-109}$ in every entry, and the restriction of the $W_{0,2}$ part to each sector is $\pm2\,\ell\ell^T$ with $\ell$ as in \eqref{eq:ellbasis}.

\emph{The components of $\cE(f)$.} For $f$ supported in $[-\lambda,\lambda]$ and $t=\log u\in[-L/2,L/2]$, $\langle\cE(f),b_j\rangle=\int_{-L/2}^{L/2}e^{t/2}\sum_{n<\lambda^2,\ ne^t\le\lambda}f(ne^t)\,\xi_j(t)\,dt$. The integrand is analytic on each interval between consecutive points $\log(\lambda/n)$, and we use Gauss--Legendre on each of them, with a number of nodes proportional to the phase $2\pi N\Delta t/L+\gamma\,n\,\Delta(e^t)/\lambda$; the prolate functions are evaluated from their Legendre expansion (Appendix \ref{app:prolate}) at every point $ne^t/\lambda$. Doubling the number of nodes changes the vectors by less than $2\cdot10^{-14}$, and single components agree with adaptive quadrature to $4\cdot10^{-16}$.

\emph{Smoothed basis.} For the trace-formula diagnostics of Section \ref{sec:condpos} the basis $E_N$ is replaced by $\vartheta_wE_N$, $\vartheta_w(u)=\vartheta\big((L/2-\vert\log u\vert-w)/w\big)$ in the notation $t=\log u\in[-L/2,L/2]$, with $\vartheta$ a Gevrey cutoff of class $G^s$ ($\vartheta=0$ on $(-\infty,0]$, $\vartheta=1$ on $[1,\infty)$), so that $\vartheta_w=1$ except within log-distance $2w$ of the ends. The matrix of $QW_\lambda$ on $\vartheta_wE_N$ is obtained from the Fourier coefficients of $\vartheta_w$ on $[-L/2,L/2]$ (a fixed one-dimensional computation, $O(N\log N)$ coefficients needed) by a congruence of the matrix on $E_{N'}$, $N'\gtrsim N+\log\lambda/w$; the Gram matrix of $\vartheta_wE_N$ is no longer the identity and the eigenvalue problem becomes generalized. Quantities to report: $\mu_\lambda(w)$ for $w=\lambda^{-3}$ and $w=\lambda^{-1}$ against the unsmoothed $\mu_\lambda$; the angles $\sin\angle(\mathfrak p_k,\eta_\lambda(k))$ recomputed on the smoothed space, which should be unchanged to exponential accuracy; and the decay of $\vert\cM\tilde f(is)\vert$ for the smoothed eigenvectors, which should follow $e^{-c(w\vert s\vert)^{1/s}}$.

\vspace{20pt}
\scriptsize
\noindent
\begin{tabular}{@{}l@{\hspace{2.2cm}}l@{\hspace{2.2cm}}l@{}}
Ronnie Andrews, Jr. & Alain Connes & Antonio La Gatta\\
ORCID \href{https://orcid.org/0009-0003-9724-3104}{0009-0003-9724-3104}
  & Coll\`ege de France & Dipartimento per la Trasformazione Digitale\\
585 Hopewell Downs Dr. & 3 Rue d'Ulm & Presidenza del Consiglio dei Ministri\\
Milton, GA 30004 & 75005 Paris & Roma\\
United States of America & France & Italy\\
\href{mailto:randrewsmath@gmail.com}{randrewsmath@gmail.com}
  & \href{mailto:alain@connes.org}{alain@connes.org}
  & \href{mailto:prof.lagatta@gmail.com}{prof.lagatta@gmail.com}
\end{tabular}
\end{document}